\documentclass[11pt,leqno]{article}
\usepackage{amssymb}
\usepackage[mathscr]{eucal}
\usepackage{amsmath,amssymb,latexsym,theorem,bbm}
\usepackage{color}
\usepackage{appendix}

\newcommand{\CC}{\mathbb{C}}

\newcommand{\NN}{\mathbb{N}}

\newcommand{\RR}{\mathbb{R}}

\newcommand{\ZZ}{\mathbb{Z}}

\newcommand{\bA}{{\boldsymbol{A}}}

\newcommand{\Bf}{{\boldsymbol{f}}}

\newcommand{\bG}{{\boldsymbol{G}}}
\newcommand{\tbG}{\widetilde{\bG}}
\newcommand{\bh}{{\boldsymbol{h}}}
\newcommand{\bI}{{\boldsymbol{I}}}
\newcommand{\bJ}{{\boldsymbol{J}}}

\newcommand{\bM}{{\boldsymbol{M}}}

\newcommand{\bQ}{{\boldsymbol{Q}}}

\newcommand{\bx}{{\boldsymbol{x}}}

\newcommand{\bZ}{{\boldsymbol{Z}}}

\newcommand{\bV}{{\boldsymbol{V}}}

\newcommand{\Beta}{{\boldsymbol{\eta}}}

\newcommand{\bGamma}{{\boldsymbol{\Gamma}}}
\newcommand{\btheta}{{\boldsymbol{\theta}}}

\newcommand{\bfeta}{{\boldsymbol{\eta}}}

\newcommand{\hbtheta}{{\widehat{\btheta}}}

\newcommand{\bvarphi}{{\boldsymbol{\varphi}}}

\newcommand{\bxi}{{\boldsymbol{\xi}}}
\newcommand{\hsigma}{\widehat{\sigma}}
\newcommand{\bzero}{{\boldsymbol{0}}}

\newcommand{\cA}{{\mathcal A}}

\newcommand{\cC}{{\mathcal C}}
\newcommand{\cD}{{\mathcal D}}

\newcommand{\cF}{{\mathcal F}}

\newcommand{\cN}{{\mathcal N}}

\newcommand{\cX}{{\mathcal X}}
\newcommand{\cY}{{\mathcal Y}}

\newcommand{\cW}{{\mathcal W}}

\newcommand{\tcX}{\widetilde{\cX}}
\newcommand{\tcY}{\widetilde{\cY}}

\newcommand{\cc}{\mathrm{c}}
\newcommand{\dd}{\mathrm{d}}
\newcommand{\ee}{\mathrm{e}}
\newcommand{\ii}{\mathrm{i}}

\DeclareMathOperator*{\argmin}{arg\,min}

\newcommand{\EE}{\operatorname{\mathbb{E}}}
\newcommand{\PP}{\operatorname{\mathbb{P}}}
\newcommand{\OO}{\operatorname{O}}

\newcommand{\diag}{\operatorname{diag}}
\renewcommand{\Re}{\operatorname{Re}}

\newcommand{\ha}{\widehat{a}}
\newcommand{\hb}{\widehat{b}}
\newcommand{\hc}{\widehat{c}}
\newcommand{\hd}{\widehat{d}}

\newcommand{\halpha}{\widehat{\alpha}}
\newcommand{\hbeta}{\widehat{\beta}}
\newcommand{\hgamma}{\widehat{\gamma}}
\newcommand{\hdelta}{\widehat{\delta}}
\newcommand{\hvare}{\widehat{\vare}}
\newcommand{\hzeta}{\widehat{\zeta}}
\newcommand{\hvarrho}{\widehat{\varrho}}

\newcommand{\tW}{\widetilde{W}}

\newcommand{\vare}{\varepsilon}

\renewcommand{\mid}{\,|\,}
\newcommand{\bmid}{\,\big|\,}

\renewcommand{\leq}{\leqslant}
\renewcommand{\geq}{\geqslant}

\newcommand{\stoch}{\stackrel{\PP}{\longrightarrow}}
\newcommand{\distr}{\stackrel{\cD}{\longrightarrow}}
\newcommand{\distre}{\stackrel{\cD}{=}}

\newcommand{\as}{\stackrel{{\mathrm{a.s.}}}{\longrightarrow}}

\newcommand{\bbone}{\mathbbm{1}}

\newcommand{\nT}{{\lfloor nT\rfloor}}

\newcommand{\proofend}{\hfill\mbox{$\Box$}}

\numberwithin{equation}{section}

\numberwithin{equation}{section}

\theoremstyle{change} \theorembodyfont{\em}
\newtheorem{Lem}{Lemma.}[section]
\newtheorem{Thm}[Lem]{Theorem.}
\newtheorem{Pro}[Lem]{Proposition.}

\newtheorem{Def}[Lem]{Definition.}

\theorembodyfont{\rm}
\newtheorem{Rem}[Lem]{Remark.}

\begin{document}

\begin{center}
 {\bfseries\Large
  On conditional least squares estimation for affine diffusions} \\[1mm]
 {\bfseries\Large
  based on continuous time observations} \\[5mm]
 {\sc\large
  Be\'ata $\text{Bolyog}^*$,
  \ Gyula $\text{Pap}$}
\end{center}

\vskip0.2cm

\noindent
 Bolyai Institute, University of Szeged,
 Aradi v\'ertan\'uk tere 1, H--6720 Szeged, Hungary

\noindent e--mails: bbeata@math.u-szeged.hu (B. Bolyog),
                    papgy@math.u-szeged.hu (G. Pap).

\noindent * Corresponding author.

\renewcommand{\thefootnote}{}
\footnote{\textit{2010 Mathematics Subject Classifications\/}:
          60J60, 62F12.}
\footnote{\textit{Key words and phrases\/}:
 affine processes, conditional least squares estimators.}

\vspace*{-10mm}

\begin{abstract}
We study asymptotic properties of conditional least squares estimators for the drift
 parameters of two-factor affine diffusions based on continuous time observations. 
We distinguish three cases: subcritical, critical and supercritical. 
For all the drift parameters, in the subcritical and supercritical cases, asymptotic
 normality and asymptotic mixed normality is proved, while in the critical case, non-standard
 asymptotic behavior is described.
\end{abstract}

\section{Introduction}

Affine processes are applied in mathematical finance in several models including interest
 rate models (e.g.\ the Cox--Ingersoll--Ross, Vasi\v{c}ek or general affine term structure
 short rate models), option pricing (e.g.\ the Heston model) and credit risk models, see
 e.g.\ Duffie, Filipovi\'{c} and Schachermayer \cite{DufFilSch}, Filipovi\'{c} \cite{Fil},
 Baldeaux and Platen \cite{BalPla}, and Alfonsi \cite{Alf}.
In this paper we consider two-factor affine processes, i.e.\ affine processes with state-space
 \ $[0, \infty) \times \RR$.
\ Dawson and Li \cite{DawLi} derived a jump-type stochastic differential equation (SDE) for
 such processes.
Specializing this result to the diffusion case, i.e.\ two-factor affine processes without
 jumps, we obtain that for every \ $a \in [0, \infty)$, \ $b, \alpha, \beta, \gamma \in \RR$,
 \ $\sigma_1, \sigma_2, \sigma_3 \in [0, \infty)$ \ and \ $\varrho \in [-1, 1]$, \ the SDE
 \begin{align}\label{2dim_affine}
  \begin{cases}
   \dd Y_t = (a - b Y_t) \, \dd t + \sigma_1 \sqrt{Y_t} \, \dd W_t , \\
   \dd X_t = (\alpha - \beta Y_t - \gamma X_t) \, \dd t
             + \sigma_2 \sqrt{Y_t} \, (\varrho \, \dd W_t + \sqrt{1 - \varrho^2} \, \dd B_t)
             + \sigma_3 \, \dd L_t ,
  \end{cases}
  \qquad t \in [0, \infty) ,
 \end{align}
 with an arbitrary initial value \ $(Y_0, X_0)$ \ with \ $\PP(Y_0 \in [0, \infty)) = 1$ \ and
 independent of a 3-dimensional standard Wiener process \ $(W_t, B_t, L_t)_{t\in[0, \infty)}$,
 \ has a pathwise unique strong solution being a two-factor affine diffusion process, and
 conversely, every two-factor affine diffusion process is a pathwise unique strong solution of
 a SDE \eqref{2dim_affine} with appropriate parameters \ $a \in [0, \infty)$,
 \ $b, \alpha, \beta, \gamma \in \RR$, \ $\sigma_1, \sigma_2, \sigma_3 \in [0, \infty)$
 \ and \ $\varrho \in [-1, 1]$.
 
In this paper we study asymptotic properties of conditional least squares estimators (CLSE)
 \ $(\ha_T, \hb_T, \halpha_T, \hbeta_T, \hgamma_T)$ \ of the drift parameters
 \ $(a, b, \alpha, \beta, \gamma)$ \ based on continuous time observations
 \ $(Y_t, X_t)_{t\in[0,T]}$ \ with \ $T > 0$.
\ This estimator is the high frequency limit in probability as \ $n \to \infty$ \ of the CLSE
 based on discrete time observations \ $(Y_{k/n}, X_{k/n})_{k\in\{0,\ldots,\nT\}}$,
 \ $n \in \NN$. 
\ We do not estimate the parameters \ $\sigma_1$, $\sigma_2$, $\sigma_3$ \ and \ $\varrho$,
 \ since for all \ $T \in (0, \infty)$, \ they are measurable functions (i.e., statistics) of
 \ $(Y_t, X_t)_{t\in[0,T]}$, \ see Appendix \ref{sigma_varrho}.
It will turn out that for the calculation of
 \ $(\ha_T, \hb_T, \halpha_T, \hbeta_T, \hgamma_T)$ \ one does not need to know the values of
 the diffusion coefficients \ $\sigma_1$, $\sigma_2$, $\sigma_3$ \ and \ $\varrho$, \ see
 \eqref{LSE_cont}.

The first coordinate process \ $Y$ \ in \eqref{2dim_affine} is called a
 Cox--Ingersoll--Ross (CIR) process (see Cox, Ingersoll and Ross \cite{CoxIngRos}). 
In the submodel consisting only of the process \ $Y$, \ Overbeck and Ryd\'en
 \cite[Theorems 3.4, 3.5 and 3.6]{OveRyd} derived the CLSE of \ $(a, b)$ \ based on continuous 
 time observations \ $(Y_t)_{t\in[0,T]}$ \ with \ $T > 0$, \ i.e., the limit in probability as
 \ $n \to \infty$ \ of the CLSE based on discrete time observations
 \ $(Y_{k/n})_{k\in\{0,\ldots,\nT\}}$, \ $n \in \NN$, \ which turns to be the same as the CLSE
 \ $(\ha_T, \hb_T)$ \ of \ $(a, b)$ \ based on continuous time observations
 \ $(Y_t, X_t)_{t\in[0,T]}$, \ and they proved strong consistency and asymptotic
 normality in case of a subcritical CIR process \ $Y$, \ i.e., when \ $b > 0$ \ and the
 initial distribution is the unique stationary distribution of the model.  
 
Barczy at al.\ \cite{BarDorLiPap3} considered a submodel of \eqref{2dim_affine} with
 \ $a \in (0, \infty)$, \ $\beta = 0$, \ $\sigma_1 = 1$, \ $\sigma_2 = 1$, \ $\varrho = 0$
 \ and \ $\sigma_3 = 0$. 
\ The estimator of the parameters \ $(\alpha, \gamma)$ \ based on continuous time observations
 \ $(X_t)_{t\in[0,T]}$ \ with \ $T > 0$ \ (which they call a least square estimator) is in
 fact the CLSE, i.e., the limit in probability as \ $n \to \infty$ \ of the CLSE based on 
 discrete time observations \ $(X_{k/n})_{k\in\{0,\ldots,\nT\}}$, \ $n \in \NN$, \ which can
 be shown by the method of the proof of Lemma \ref{LEMMA_LSE_exist}.
They proved strong consistency and asymptotic normality in case of a subcritical process
 \ $(Y, X)$, \ i.e., when \ $b > 0$ \ and \ $\gamma > 0$.
 
Barczy at al.\ \cite{BarNyuPap} considered the so-called Heston model, which is a submodel of
 \eqref{2dim_affine} with \ $a, \sigma_1, \sigma_2 \in (0, \infty)$, \ $\gamma = 0$,
 \ $\varrho \in (-1, 1)$ \ and \ $\sigma_3 = 0$. 
\ The estimator of the parameters \ $(a, b, \alpha, \beta)$ \ based on continuous time
 observations \ $(Y_t, X_t)_{t\in[0,T]}$ \ with \ $T > 0$ \ (which they call least square
 estimator) is in fact the CLSE, i.e., the limit in probability as \ $n \to \infty$ \ of the
 CLSE based on discrete time observations \ $(Y_{k/n}, X_{k/n})_{k\in\{0,\ldots,\nT\}}$,
 \ $n \in \NN$ \ which can be shown by the method of the proof of Lemma \ref{LEMMA_LSE_exist}. 
They proved strong consistency and asymptotic normality in case of a subcritical process
 \ $(Y, X)$, \ i.e., when \ $b > 0$.
\ Note that Barczy and Pap \cite{BarPap2} studied the maximum likelihood estimator (MLE) of
 the parameters \ $(a, b, \alpha, \beta)$ \ in this Heston model under the additional
 assumption \ $a \geq \frac{\sigma_1^2}{2}$.   
\ In the subcritical case, i.e., when \ $b > 0$, \ they proved strong consistency and
 asymptotic normality of the MLE of \ $(a, b, \alpha, \beta)$ \ under the additional
 assumption \ $a > \frac{\sigma_1^2}{2}$.
\ In the critical case, namely, if \ $b = 0$, \ they showed weak consistency of the MLE of
 \ $(a, b, \alpha, \beta)$ \ and determined the asymptotic behavior of the MLE under the
 additional assumption \ $a > \frac{\sigma_1^2}{2}$.
\ In a special supercritical case, namely, when \ $b < 0$, \ they showed strong consistency of
 the MLE of \ $b$, \ weak consistency of the MLE of \ $\beta$ \ and proved asymptotic mixed
 normality of the MLE of \ $(a, b, \alpha, \beta)$. 
\ Barczy at al.\ \cite{BarBenKebPap1, BarBenKebPap2} studied the asymptotic
 behavior of maximum likelihood estimators for a jump-type Heston model and for the growth
 rate of a jump-type CIR process, respectively, based on continuous time observations.
 
We consider general two-factor affine diffusions \eqref{2dim_affine}.
In the subcritical case, i.e., when \ $b > 0$ \ and \ $\gamma > 0$, \ we prove strong
 consistency and asymptotic normality of \ $(\ha_T, \hb_T, \halpha_T, \hbeta_T, \hgamma_T)$
 \ under the additional assumptions \ $a > 0$, \ $\sigma_1 > 0$ \ and
 \ $(1 - \varrho^2) \sigma_2^2 + \sigma_3^2 > 0$.
\ In a special critical case, namely if \ $b = 0$ \ and \ $\gamma = 0$, \ we show weak
 consistency of \ $(\hb_T, \hbeta_T, \hgamma_T)$ \ and determine the asymptotic behavior of
 \ $(\ha_T, \hb_T, \halpha_T, \hbeta_T, \hgamma_T)$ \ under the additional assumptions
 \ $\beta = 0$ \ and \ $(1 - \varrho^2) \sigma_2^2 + \sigma_3^2 > 0$.
\ In a special supercritical case, namely, when \ $\gamma < b < 0$, \ we show strong
 consistency of \ $\hb_T$, \ weak consistency of \ $(\hbeta_T, \hgamma_T)$ \ and prove
 asymptotic mixed normality of \ $(\ha_T, \hb_T, \halpha_T, \hbeta_T, \hgamma_T)$ \ under the
 additional assumptions \ $\alpha \beta \leq 0$, \ $\sigma_1 > 0$, \ and either
 \ $\sigma_3 > 0$, \ or
 \ $\bigl(a - \frac{\sigma_1^2}{2}\bigr) (1 - \varrho^2) \sigma_2^2 > 0$.
\ Note that we decided to deal with the CLSE of \ $(a, b, \alpha, \beta, \gamma)$, \ since
 the MLE of \ $(a, b, \alpha, \beta, \gamma)$ \ contains, for example,
 \ $\int_0^T \frac{X_t}{(1-\varrho^2)\sigma_2^2Y_t+\sigma_3^2} \, \dd t$, \ and the question
 of the asymptotic behavior of this integral as \ $T \to \infty$ \ is still open in the
 critical and supercritical cases.

\section{The affine two-factor model}\label{Section_preliminaires}

Let \ $\NN$, \ $\ZZ_+$, \ $\RR$, \ $\RR_+$, \ $\RR_{++}$, \ $\RR_-$, \ $\RR_{--}$ \ and
 \ $\CC$ \ denote the sets of positive integers, non-negative integers, real numbers,
 non-negative real numbers, positive real numbers, non-positive real numbers, negative real
 numbers and complex numbers, respectively.
For \ $x , y \in \RR$, \ we will use the notations \ $x \land y := \min(x, y)$ \ and
 \ $x \lor y := \max(x, y)$.
\ By \ $C^2_\cc(\RR_+ \times \RR, \RR)$, \ we denote the set of twice continuously
 differentiable real-valued functions on \ $\RR_+ \times \RR$ \ with compact support.
\ Let \ $(\Omega, \cF, \PP)$ \ be a probability space equipped with the augmented filtration
 \ $(\cF_t)_{t\in\RR_+}$ \ corresponding to \ $(W_t, B_t, L_t)_{t\in\RR_+}$ \ and a given
 initial value \ $(\eta_0, \xi_0)$ \ being independent of \ $(W_t, B_t, L_t)_{t\in\RR_+}$
 \ such that \ $\PP(\eta_0 \in \RR_+) = 1$, \ constructed as in Karatzas and Shreve
 \cite[Section 5.2]{KarShr}.
Note that \ $(\cF_t)_{t\in\RR_+}$ \ satisfies the usual conditions, i.e., the filtration
 \ $(\cF_t)_{t\in\RR_+}$ \ is right-continuous and \ $\cF_0$ \ contains all the $\PP$-null
 sets in \ $\cF$.
\ We will denote the convergence in distribution, convergence in probability, almost surely
 convergence and equality in distribution by \ $\distr$, \ $\stoch$, \ $\as$ \ and
 \ $\distre$, \ respectively.
By \ $\|\bx\|$ \ and \ $\|\bA\|$, \ we denote the Euclidean norm of a vector \ $\bx \in \RR^d$
 \ and the induced matrix norm of a matrix \ $\bA \in \RR^{d \times d}$, \ respectively.
By \ $\bI_d \in \RR^{d \times d}$, \ we denote the \ $d\times d$ unit matrix.
For quadratic matrices \ $\bA_1, \ldots, \bA_k$, \ $\diag(\bA_1, \ldots, \bA_k)$ \ will
 denote the quadratic block matrix containing the matrices \ $\bA_1, \ldots, \bA_k$ \ in its
 diagonal.
 
The next proposition is about the existence and uniqueness of a strong solution of the SDE
 \eqref{2dim_affine}, see Bolyog and Pap \cite[Proposition 2.2]{BolPap}.

\begin{Pro}\label{Pro_affine}
Let \ $(\eta_0, \xi_0)$ \ be a random vector independent of the process
 \ $(W_t, B_t, L_t)_{t\in\RR_+}$ \ satisfying \ $\PP(\eta_0 \in \RR_+) = 1$.
\ Then for all \ $a \in \RR_+$, \ $b, \alpha, \beta, \gamma \in \RR$,
 \ $\sigma_1, \sigma_2, \sigma_3 \in \RR_+$, \ $\varrho \in [-1, 1]$, \ there is a (pathwise)
 unique strong solution \ $(Y_t, X_t)_{t\in\RR_+}$ \ of the SDE \eqref{2dim_affine} such that
 \ $\PP((Y_0, X_0) = (\eta_0, \xi_0)) = 1$ \ and \ $\PP(\text{$Y_t \in \RR_+$ \ for all
 \ $t \in \RR_+$}) = 1$.
\ Further, for all \ $s, t \in \RR_+$ \ with \ $s \leq t$, \ we have
 \begin{equation}\label{SolutionY}
  Y_t = \ee^{-b(t-s)} Y_s + a \int_s^t \ee^{-b(t-u)} \, \dd u
        + \sigma_1 \int_s^t \ee^{-b(t-u)} \sqrt{Y_u} \, \dd W_u
 \end{equation}
 and
 \begin{equation}\label{SolutionX}
  \begin{aligned}
   X_t &= \ee^{-\gamma(t-s)} X_s
          + \int_s^t \ee^{-\gamma(t-u)} (\alpha - \beta Y_u) \, \dd u \\
       &\quad
          + \sigma_2
            \int_s^t
             \ee^{-\gamma(t-u)}\sqrt{Y_u}
             \, (\varrho \, \dd W_u + \sqrt{1 - \varrho^2} \, \dd B_u)
          + \sigma_3 \int_s^t \ee^{-\gamma(t-u)} \, \dd L_u .
  \end{aligned}
 \end{equation}
Moreover, \ $(Y_t, X_t)_{t\in\RR_+}$ \ is a two-factor affine process with infinitesimal
 generator
 \begin{equation}\label{infgen}
  \begin{aligned}
   (\cA_{(Y,X)} f)(y, x)
   &= (a - b y) f_1'(y, x) + (\alpha - \beta y - \gamma x) f_2'(y, x) \\
   &\quad 
      + \frac{1}{2} y
        \bigl[ \sigma_1^2 f_{1,1}''(y, x)
               + 2 \varrho \sigma_1 \sigma_2 f_{1,2}''(y, x)
               + \sigma_2^2 f_{2, 2}''(y,x) \bigr]
      + \frac{1}{2} \sigma_3^2 f_{2, 2}''(y,x) ,
  \end{aligned}
 \end{equation}
 where \ $(y ,x) \in \RR_+ \times \RR$, \ $f \in \cC^2_c(\RR_+ \times \RR, \RR)$, \ and
 \ $f_i'$, \ $i \in \{1, 2\}$, \ and \ $f_{i,j}''$, \ $i, j \in\{1, 2\}$, \ denote the first
 and second order partial derivatives of $f$ with respect to its \ $i$-th and \ $i$-th and
 \ $j$-th variables.
 
Conversely, every two-factor affine diffusion process is a (pathwise) unique strong solution
 of a SDE \eqref{2dim_affine} with suitable parameters \ $a \in \RR_+$,
 \ $b, \alpha, \beta, \gamma \in \RR$, \ $\sigma_1, \sigma_2, \sigma_3 \in \RR_+$ \ and
 \ $\varrho \in [-1, 1]$.
\end{Pro}

The next proposition gives the asymptotic behavior of the first moment of the process
 \ $(Y_t, X_t)_{t\in\RR_+}$ \ as \ $t \to \infty$, \ see Bolyog and Pap
 \cite[Proposition 2.3]{BolPap}.

\begin{Pro}\label{Pro_moments_asymptotics}
Let us consider the two-factor affine diffusion model \eqref{2dim_affine} with
 \ $a \in \RR_+$, \ $b, \alpha, \beta, \gamma \in \RR$,
 \ $\sigma_1, \sigma_2, \sigma_3 \in \RR_+$, \ $\varrho \in [-1, 1]$.
\ Suppose that \ $\EE(Y_0 |X_0|) < \infty$.
\ In case of \ $b \in \RR_{++}$ \ we have \ $\EE(Y_t) = \frac{a}{b} + \OO(\ee^{-bt})$ \ and
 \[
   \EE(X_t)
   = \begin{cases}
      \frac{\alpha}{\gamma} - \frac{a\beta}{b\gamma} + \OO(\ee^{-(b\land\gamma)t}) ,
       &\text{$\gamma \in \RR_{++}$,} \\[1mm]
      \bigl(\alpha - \frac{a\beta}{b}\bigr) t + \OO(1) ,
       &\text{$\gamma = 0$,} \\[1mm]
      \bigl(\frac{\beta}{\gamma-b} \EE(Y_0) + \EE(X_0) - \frac{\alpha}{\gamma}
            + \frac{a\beta}{b\gamma} - \frac{a\beta}{(\gamma-b)b}\bigr)
      \ee^{-\gamma t}
      + \OO(1) ,
       &\text{$\gamma \in \RR_{--}$.}
     \end{cases}
 \]
In case of \ $b = 0$ \ we have \ $\EE(Y_t) = a t + \OO(1)$ \ and
 \[
   \EE(X_t)
   = \begin{cases}
      - \frac{a\beta}{\gamma} t + \OO(1) ,
       &\text{$\gamma \in \RR_{++}$,} \\[1mm]
      - \frac{1}{2} a \beta t^2 + \OO(t) ,
       &\text{$\gamma = 0$,} \\[1mm]
      \bigl(\frac{\beta}{\gamma} \EE(Y_0) + \EE(X_0) - \frac{\alpha}{\gamma} 
            - \frac{a\beta}{\gamma^2}\bigr)
      \ee^{-\gamma t}
      + \OO(t) ,
       &\text{$\gamma \in \RR_{--}$.}
     \end{cases}
 \]
In case of \ $b \in \RR_{--}$ \ we have
 \ $\EE(Y_t) = \bigl(\EE(Y_0) - \frac{a}{b}\bigr) \ee^{-bt} + \OO(1)$ \ and
 \[
   \EE(X_t)
   = \begin{cases}
      \bigl(- \frac{\beta}{\gamma-b} \EE(Y_0) + \frac{a\beta}{(\gamma-b)b}\bigr) \ee^{-bt}
       + \OO(1) ,
       &\text{$\gamma \in \RR_{++}$,} \\[1mm]
      \bigl(\frac{\beta}{b} \EE(Y_0) + \EE(X_0) - \frac{\beta a}{b^2}\bigr) \ee^{-bt}
      + \OO(t) ,
       &\text{$\gamma = 0$,} \\[1mm]
      \bigl(- \frac{\beta}{\gamma-b} \EE(Y_0) + \frac{a\beta}{(\gamma-b)b}\bigr) \ee^{-bt}
      + \OO(\ee^{-\gamma t}) ,
       &\text{$\gamma \in (b, 0)$,} \\[1mm]
      \bigl(- \beta \EE(Y_0) + \frac{a\beta}{b}\bigr) t \ee^{-bt} + \OO(\ee^{-\gamma t}) ,
       &\text{$\gamma = b$,} \\[1mm]
      \bigl(\frac{\beta}{\gamma-b} \EE(Y_0) + \EE(X_0) - \frac{\alpha}{\gamma}
            + \frac{a\beta}{b\gamma} - \frac{a\beta}{b(\gamma-b)}\bigr)
      \ee^{-\gamma t}
      + \OO(\ee^{-bt}) ,
       &\text{$\gamma \in (-\infty, b)$.}
    \end{cases}
 \] 
\end{Pro}

Based on the asymptotic behavior of the first moment of the process \ $(Y_t, X_t)_{t\in\RR_+}$
 \ as \ $t \to \infty$, \ we can classify two-factor affine diffusions in the following way. 

\begin{Def}\label{Def_criticality}
Let \ $(Y_t, X_t)_{t\in\RR_+}$ \ be the unique strong solution of the SDE
 \eqref{2dim_affine} satisfying \ $\PP(Y_0 \in \RR_+) = 1$.
\ We call \ $(Y_t, X_t)_{t\in\RR_+}$ \ subcritical, critical or supercritical if
 \ $b \land \gamma \in \RR_{++}$, \ $b \land \gamma = 0$ \ or
 \ $b \land \gamma \in \RR_{--}$, \ respectively.
\end{Def}

\section{CLSE based on continuous time observations}
\label{section_LSE_continuous}

Overbeck and Ryd\'en \cite{OveRyd} investigated the CIR process \ $Y$, \ and for each
 \ $T \in \RR_{++}$, \ they defined a CLSE \ $(\ha_T, \hb_T)$ \ of \ $(a, b)$ \ based on
 continuous time observations \ $(Y_t)_{t\in[0,T]}$ \ as the limit in probability of the CLSE
 \ $(\ha_{T,n}, \hb_{T,n})$ \ of \ $(a, b)$ \ based on discrete time observations
 \ $(Y_{\frac{iT}{n}})_{i\in\{0,1,\ldots,n\}}$ \ as \ $n \to \infty$.

We consider a two-factor affine diffusion process \ $(Y_t, X_t)_{t\in\RR_+}$ \ given in
 \eqref{2dim_affine} with known \ $\sigma_1 \in \RR_{++}$, \ $\sigma_2, \sigma_3 \in \RR_+$
 \ and \ $\varrho \in [-1, 1]$, \ and with a random initial value \ $(\eta_0, \zeta_0)$
 \ independent of \ $(W_t, B_t, L_t)_{t\in\RR_+}$ \ satisfying
 \ $\PP(\eta_0 \in \RR_+) = 1$, \ and we will consider
 \ $\btheta = (a, b, \alpha, \beta, \gamma)^\top \in \RR_+ \times \RR^4$ \ as a parameter.
The aim of the following discussion is to construct a CLSE of \ $\btheta$ \ based on
 continuous time observations \ $(Y_t, X_t)_{t\in[0,T]}$ \ with some \ $T \in \RR_{++}$. 

Let us recall the CLSE \ $\hbtheta_{T,n}$ \ of \ $\btheta$ \ based on discrete time
 observations \ $(Y_{\frac{i}{n}}, X_{\frac{i}{n}})_{i\in\{0,1,\ldots,\nT\}}$ \ with some
 \ $n \in \NN$, \ which can be obtained by solving the extremum problem
 \[
   \hbtheta_{T,n}
   := \argmin_{\btheta\in\RR^5}
       \sum_{i=1}^\nT
        \left[\Bigl(Y_{\frac{i}{n}}
                    - \EE\Bigl(Y_{\frac{i}{n}} \,\Big|\, \cF_{\frac{i-1}{n}}\Bigr)\Bigr)^2
              + \Bigl(X_{\frac{i}{n}}
                      - \EE\Bigl(X_{\frac{i}{n}}
                                 \,\Big|\, \cF_{\frac{i-1}{n}}\Bigr)\Bigr)^2\right] .
 \]
By \eqref{SolutionY} and \eqref{SolutionX}, together with Proposition 3.2.10 in Karatzas
 and Shreve \cite{KarShr}, for all \ $s, t \in \RR_+$ \ with \ $s \leq t$, \ we obtain
 \begin{gather}\label{SolutionYcond}
  \EE(Y_t \mid \cF_s)
  = \ee^{-b(t-s)} Y_s + a \int_s^t \ee^{-b(t-u)} \, \dd u , \\
  \begin{aligned}
   \EE(X_t \mid \cF_s)
   &= \ee^{-\gamma(t-s)} X_s + \alpha \int_s^t \ee^{-\gamma(t-u)} \, \dd u
      - \beta Y_s \int_s^t \ee^{-\gamma(t-u)-b(u-s)} \, \dd u \\
   &\quad
      - a \beta
        \int_s^t \ee^{-\gamma(t-u)} \biggl(\int_s^u \ee^{-b(u-v)} \, \dd v\biggr) \dd u .
  \end{aligned}\label{SolutionXcond}
 \end{gather}
Thus, for all \ $i \in \NN$, \ we have
 \[
   \EE\Bigl(Y_{\frac{i}{n}} \mid \cF_{\frac{i-1}{n}}\Bigr)
   = \ee^{-\frac{b}{n}} Y_{\frac{i-1}{n}} + a \int_0^{\frac{1}{n}} \ee^{-bw} \, \dd w
 \]
 and
 \begin{align*}
  \EE\Bigl(X_{\frac{i}{n}} \mid\cF_{\frac{i-1}{n}}\Bigr)
  &= \ee^{-\frac{\gamma}{n}} X_{\frac{i-1}{n}}
     + \alpha \int_0^{\frac{1}{n}} \ee^{-\gamma w} \, \dd w
     - \beta Y_{\frac{i-1}{n}}
       \int_0^{\frac{1}{n}} \ee^{(\gamma-b)w-\frac{\gamma}{n}} \, \dd w \\
  &\quad
     - a \beta
       \int_0^{\frac{1}{n}}
        \ee^{\gamma w-\frac{\gamma}{n}} \biggl(\int_0^w \ee^{-b(w-v)} \, \dd v\biggr) \dd w .
 \end{align*}
Consequently,
 \begin{equation}\label{CLSED}
  \begin{aligned}
   \hbtheta_{T,n}
   = \argmin_{(a,b,\alpha,\beta,\gamma)^\top\in\RR^5}
      \sum_{i=1}^\nT
       \biggl[&\Bigl(Y_{\frac{i}{n}} - Y_{\frac{i-1}{n}}
                     - \Bigl(c - d Y_{\frac{i-1}{n}}\Bigr)\Bigr)^2 \\
              &+ \Bigl(X_{\frac{i}{n}} - X_{\frac{i-1}{n}}
                       - \Bigl(\delta - \vare Y_{\frac{i-1}{n}}
                               - \zeta X_{\frac{i-1}{n}}\Bigr)\Bigr)^2\biggr] ,
  \end{aligned}
 \end{equation}
 where
 \begin{equation}\label{gn}
  (c, d, \delta, \vare, \zeta)
  := (c_n(a, b), d_n(b), \delta_n(a, b, \alpha, \beta, \gamma), \vare_n(b, \beta, \gamma),
      \zeta_n(\gamma))
  := g_n(a, b, \alpha, \beta, \gamma)
 \end{equation}
 with
 \begin{gather*}
  c := c_n(a, b) := a \int_0^{\frac{1}{n}} \ee^{-bw} \, \dd w , \qquad
  d := d_n(b) := 1 - \ee^{-\frac{b}{n}} , \\
  \delta := \delta_n(a, b, \alpha, \beta, \gamma)
  := \alpha \int_0^{\frac{1}{n}} \ee^{-\gamma w} \, \dd w
     - a \beta
       \int_0^{\frac{1}{n}}
        \ee^{\gamma w-\frac{\gamma}{n}}
        \biggl(\int_0^w \ee^{-b(w-v)} \, \dd v\biggr) \dd w , \\
  \vare := \vare_n(b, \beta, \gamma)
  := \beta \int_0^{\frac{1}{n}} \ee^{(\gamma-b)w-\frac{\gamma}{n}} \, \dd w , \qquad
  \zeta := \zeta_n(\gamma) := 1 - \ee^{-\frac{\gamma}{n}} .
 \end{gather*}
The function \ $g_n : \RR^5 \to \RR \times (-\infty, 1) \times \RR^2 \times (-\infty, 1)$ \ is
 bijective, so first we determine the CLSE
 \ $(\hc_{T,n}, \hd_{T,n}, \hdelta_{T,n}, \hvare_{T,n}, \hzeta_{T,n})$ \ of the transformed
 parameters \ $(c, d, \delta, \vare, \zeta)$ \ by minimizing the sum on the right-hand side of
 \eqref{CLSED} with respect to \ $(c, d, \delta, \vare, \zeta)$.
\ We have
 \begin{align*}
  \bigl(\hc_{T,n}, \hd_{T,n}\bigr)
  &= \argmin_{(c,d)^\top\in\RR^2}
      \sum_{i=1}^\nT
       \Bigl(Y_{\frac{i}{n}} - Y_{\frac{i-1}{n}}
             - \Bigl(c - d Y_{\frac{i-1}{n}}\Bigr)\Bigr)^2 , \\
  \bigl(\hdelta_{T,n}, \hvare_{T,n}, \hzeta_{T,n}\bigr)
  &= \argmin_{(\delta,\vare,\zeta)^\top\in\RR^3}
      \sum_{i=1}^\nT
       \Bigl(X_{\frac{i}{n}} - X_{\frac{i-1}{n}}
             - \Bigl(\delta - \vare Y_{\frac{i-1}{n}}
                     - \zeta X_{\frac{i-1}{n}}\Bigr)\Bigr)^2 ,
 \end{align*}
 hence, similarly as on page 675 in Barczy et al. \cite{BarDorLiPap}, we get
 \begin{equation}\label{CLSE_disc}
   \begin{bmatrix}
    \hc_{T,n} \\
    \hd_{T,n}
   \end{bmatrix}
   = \bigl(\bGamma_{T,n}^{(1)}\bigr)^{-1} \bvarphi_{T,n}^{(1)} , \qquad
   \begin{bmatrix}
    \hdelta_{T,n} \\
    \hvare_{T,n} \\
    \hzeta_{T,n}
   \end{bmatrix}
   = \bigl(\bGamma_{T,n}^{(2)}\bigr)^{-1} \bvarphi_{T,n}^{(2)}
 \end{equation}
 with
 \begin{gather*}
  \bGamma_{T,n}^{(1)}
  := \begin{bmatrix}
      \nT & - \sum\limits_{i=1}^\nT Y_{\frac{i-1}{n}} \\
      - \sum\limits_{i=1}^\nT Y_{\frac{i-1}{n}} & \sum\limits_{i=1}^\nT Y_{\frac{i-1}{n}}^2
     \end{bmatrix} , \qquad
  \bvarphi_{T,n}^{(1)}
  := \begin{bmatrix}
      Y_{\frac{\nT}{n}} - Y_0 \\
      - \sum\limits_{i=1}^\nT
         \Bigl(Y_{\frac{i}{n}} - Y_{\frac{i-1}{n}}\Bigr) Y_{\frac{i-1}{n}}
     \end{bmatrix} , 
 \end{gather*}
 \begin{gather*}
  \bGamma_{T,n}^{(2)}
  := \begin{bmatrix}
      \nT & - \sum\limits_{i=1}^\nT Y_{\frac{i-1}{n}}
       & - \sum\limits_{i=1}^\nT X_{\frac{i-1}{n}} \\
      - \sum\limits_{i=1}^\nT Y_{\frac{i-1}{n}} & \sum\limits_{i=1}^\nT Y_{\frac{i-1}{n}}^2
       & \sum\limits_{i=1}^\nT Y_{\frac{i-1}{n}} X_{\frac{i-1}{n}} \\
      - \sum\limits_{i=1}^\nT X_{\frac{i-1}{n}}
       & \sum\limits_{i=1}^\nT Y_{\frac{i-1}{n}} X_{\frac{i-1}{n}}
       & \sum\limits_{i=1}^\nT X_{\frac{i-1}{n}}^2
     \end{bmatrix} , \qquad
  \bvarphi_{T,n}^{(2)}
  := \begin{bmatrix}
      X_{\frac{\nT}{n}} - X_0 \\
      - \sum\limits_{i=1}^\nT
         \Bigl(X_{\frac{i}{n}} - X_{\frac{i-1}{n}}\Bigr) Y_{\frac{i-1}{n}} \\
      - \sum\limits_{i=1}^\nT
         \Bigl(X_{\frac{i}{n}} - X_{\frac{i-1}{n}}\Bigr) X_{\frac{i-1}{n}}
     \end{bmatrix}
 \end{gather*}
 on the event where the random matrices \ $\bGamma_{T,n}^{(1)}$ \ and \ $\bGamma_{T,n}^{(2)}$
 \ are invertible.
 
\begin{Lem}\label{LEMMA_CLSED_exist}
Let us consider the two-factor affine diffusion model \eqref{2dim_affine} with
 \ $a \in \RR_+$, \ $b, \alpha, \beta, \gamma \in \RR$, \ $\sigma_1 \in \RR_{++}$,
 \ $\sigma_2, \sigma_3 \in \RR_+$ \ and \ $\varrho \in [-1, 1]$ \ with a random initial value
 \ $(\eta_0, \zeta_0)$ \ independent of \ $(W_t, B_t, L_t)_{t\in\RR_+}$ \ satisfying
 \ $\PP(\eta_0 \in \RR_+) = 1$. 
\ Suppose that \ $(1 - \varrho^2) \sigma_2^2 + \sigma_3^2 > 0$.
\ Then for each \ $T \in \RR_{++}$ \ and \ $n \in \NN$, \ the random matrices
 \ $\bGamma_{T,n}^{(1)}$ \ and \ $\bGamma_{T,n}^{(2)}$ \ are invertible almost surely, and
 hence there exists a unique CLSE
 \ $\bigl(\hc_{T,n}, \hd_{T,n}, \hdelta_{T,n}, \hvare_{T,n}, \hzeta_{T,n}\bigr)$ \ of
 \ $(c, d, \delta, \vare, \zeta)$ \ taking the form given in \eqref{CLSE_disc}.
\end{Lem}

\noindent{\bf Proof.}
The aim of the following discussion is to show that the random matrix \ $\bGamma_{T,n}^{(1)}$
 \ is almost surely strictly positive definite checking that for all
 \ $\bx \in \RR^2 \setminus \{\bzero\}$, \ we have \ $\bx^\top \bGamma_{T,n}^{(1)} \bx > 0$
 \ almost surely.
Indeed, for all \ $\bx = (x_1, x_2)^\top \in \RR^2$,
 \[
   \begin{bmatrix} x_1 \\ x_2 \end{bmatrix}^\top
   \bGamma_{T,n}^{(1)}
   \begin{bmatrix} x_1 \\ x_2 \end{bmatrix}
   = \sum\limits_{i=1}^\nT
      \begin{bmatrix} x_1 \\ x_2 \end{bmatrix}^\top
      \begin{bmatrix} 1 \\ - Y_{\frac{i-1}{n}} \end{bmatrix}
      \begin{bmatrix} 1 \\ - Y_{\frac{i-1}{n}} \end{bmatrix}^\top
      \begin{bmatrix} x_1 \\ x_2 \end{bmatrix}
      \dd s
   = \sum\limits_{i=1}^\nT (x_1 - x_2 Y_{\frac{i-1}{n}})^2 \, \dd s
   \geq 0 ,
 \]
 and \ $\bx^\top \bGamma_{T,n}^{(1)} \bx = 0$ \ if and only if
 \ $x_1 - x_2 Y_{\frac{i-1}{n}} = 0$ \ for all \ $i \in \{0, 1, \ldots, \nT\}$, \ which
 happens with probability 0, since \ $\bx = (x_1, x_2)^\top \ne \bzero$ \ and, for each
 \ $i \in \{1, \ldots, \nT\}$, \ the distribution \ of \ $Y_{\frac{i-1}{n}}$ \ is absolutely
 continuous, since the conditional distribution of \ $Y_{\frac{i-1}{n}}$ \ given \ $Y_0$ \ is
 absolutely continuous, see, eg., Ben Alaya and Kebaier \cite[Proof of Proposition 2]{BenKeb2}
 and Ikeda and Watanabe \cite[page 222]{IkeWat}.

In a similar way, the stochastic matrix \ $\bGamma_{T,n}^{(2)}$ \ is almost surely strictly
 positive definite, since for all \ $\bx \in \RR^3 \setminus \{\bzero\}$, \ we have
 \ $\bx^\top \bGamma_{T,n}^{(2)} \bx > 0$ \ almost surely.
Indeed, for all \ $\bx = (x_1, x_2, x_3)^\top \in \RR^3$,
 \[
   \begin{bmatrix} x_1 \\ x_2 \\ x_3 \end{bmatrix}^\top
   \bGamma_{T,n}^{(2)}
   \begin{bmatrix} x_1 \\ x_2 \\ x_3 \end{bmatrix}
   = \sum\limits_{i=1}^\nT
      \begin{bmatrix} x_1 \\ x_2 \\ x_3 \end{bmatrix}^\top
      \begin{bmatrix} 1 \\ - Y_{\frac{i-1}{n}} \\ - X_{\frac{i-1}{n}} \end{bmatrix}
      \begin{bmatrix} 1 \\ - Y_{\frac{i-1}{n}} \\ - X_{\frac{i-1}{n}} \end{bmatrix}^\top
      \begin{bmatrix} x_1 \\ x_2 \\ x_3 \end{bmatrix}
      \dd s
   = \sum\limits_{i=1}^\nT (x_1 - x_2 Y_{\frac{i-1}{n}} - x_3 X_{\frac{i-1}{n}})^2 \, \dd s
   \geq 0 ,
 \]
 and \ $\bx^\top \bGamma_{T,n}^{(2)} \bx = 0$ \ if and only if
 \ $x_1 - x_2 Y_{\frac{i-1}{n}} - x_3 X_{\frac{i-1}{n}} = 0$ \ for all
 \ $i \in \{0, 1, \ldots, \nT\}$, \ which happens with probability 0, since
 \ $\bx = (x_1, x_2, x_3)^\top \ne \bzero$ \ and, for each
 \ $i \in \{1, \ldots, \nT\}$, \ the distribution \ of
 \ $(Y_{\frac{i-1}{n}}, X_{\frac{i-1}{n}})$ \ is absolutely continuous, because, as in the
 proof of part (b) in the proof of Theorem \ref{Thm_ergodic1} in Bolyog and Pap \cite{BolPap},
 the conditional distribution of \ $(Y_{\frac{i-1}{n}}, X_{\frac{i-1}{n}})$ \ given
 \ $(Y_0, X_0)$ \ is absolutely continuous.
\proofend

\begin{Rem}
The first order Taylor approximation of \ $g_n(a, b, \alpha, \beta, \gamma)$ \ at
 \ $(0,0,0,0,0)$ \ is \ $\frac{1}{n}(a, b, \alpha, \beta, \gamma)$, \ hence we obtain the
 first order Taylor approximations
 \begin{gather*}
  Y_{\frac{i}{n}} - \EE\Bigl(Y_{\frac{i}{n}} \mid \cF_{\frac{i-1}{n}}\Bigr)
  \approx
  Y_{\frac{i}{n}} - Y_{\frac{i-1}{n}} - \frac{1}{n} \Bigl(a - b Y_{\frac{i-1}{n}}\Bigr) , \\
  X_{\frac{i}{n}} - \EE\Bigl(X_{\frac{i}{n}} \mid \cF_{\frac{i-1}{n}}\Bigr)
  \approx
  X_{\frac{i}{n}} - X_{\frac{i-1}{n}}
  - \frac{1}{n} \Bigl(\alpha - \beta Y_{\frac{i-1}{n}} - \gamma X_{\frac{i-1}{n}}\Bigr) .
 \end{gather*}
Using these approximations, one can define an approximate CLSE
 \ $\hbtheta_{T,n}^{\mathrm{approx}}$ \ of \ $\btheta$ \ based on discrete time observations
 \ $(Y_i, X_i)_{i\in\{0,1,\ldots,\nT\}}$, \ $n \in \NN$, \ by solving the extremum problem
 \begin{align*}
  \hbtheta_{T,n}^{\mathrm{approx}}
  := \argmin_{(a,b,\alpha,\beta,\gamma)^\top\in\RR^5}
      \sum_{i=1}^\nT
       \biggl[&\Bigl(Y_{\frac{i}{n}} - Y_{\frac{i-1}{n}}
                     - \frac{1}{n} \Bigl(a - b Y_{\frac{i-1}{n}}\Bigr)\Bigr)^2 \\
              &+ \Bigl(X_{\frac{i}{n}} - X_{\frac{i-1}{n}}
                       - \frac{1}{n}
                         \Bigl(\alpha - \beta Y_{\frac{i-1}{n}}
                               - \gamma X_{\frac{i-1}{n}}\Bigr)\Bigr)^2\biggr] ,
 \end{align*}
 hence
 \ $\hbtheta_{T,n}^{\mathrm{approx}}
    = n \bigl(\hc_{T,n}, \hd_{T,n}, \hdelta_{T,n}, \hvare_{T,n}, \hzeta_{T,n}\bigr)^\top$.
\ This definition of approximate CLSE can be considered as the definition of LSE given in Hu
 and Long \cite[formula (1.2)]{HuLon2} for generalized Ornstein--Uhlenbeck processes driven by
 \ $\alpha$-stable motions, see also Hu and Long \cite[formula (3.1)]{HuLon3}.
For a heuristic motivation of the estimator \ $\hbtheta_n^{\mathrm{approx}}$ \ based on
 discrete observations, see, e.g., Hu and Long \cite[page 178]{HuLon1} (formulated for
 Langevin equations).
\proofend
\end{Rem}

We have
 \begin{gather*}
  \frac{1}{n} \bGamma_{T,n}^{(1)}
  \as
  \begin{bmatrix}
   T & -\int_0^T Y_s \, \dd s \\
   - \int_0^T Y_s \, \dd s & \int_0^T Y_s^2 \, \dd s
  \end{bmatrix}
  =: \bG_T^{(1)} , \\
  \frac{1}{n} \bGamma_{T,n}^{(2)}
  \as
  \begin{bmatrix}
   T & -\int_0^T Y_s \, \dd s & -\int_0^T X_s \, \dd s\\
   -\int_0^T Y_s \, \dd s & \int_0^T {Y_s}^2 \, \dd s & \int_0^T X_s Y_s \, \dd s\\
   -\int_0^T X_s \, \dd s & \int_0^T X_s Y_s \, \dd s & \int_0^T X_s^2 \, \dd s
  \end{bmatrix}
  =: \bG_T^{(2)}
 \end{gather*}
 as \ $n \to \infty$, \ since \ $(Y_t, X_t)_{t\in\RR_+}$ \ is almost surely continuous.
By Proposition I.4.44 in Jacod and Shiryaev \cite{JSh} with the Riemann sequence of
 deterministic subdivisions \ $\left(\frac{i}{n} \land T\right)_{i\in\NN}$, \ $n \in \NN$.,
 \ we obtain
 \[
   \bvarphi_{T,n}^{(1)}
   \stoch
   \begin{bmatrix}
    Y_T - Y_0 \\
    -\int_0^T  Y_s \, \dd Y_s
   \end{bmatrix}
   =: \Bf_T^{(1)} , \qquad
   \bvarphi_{T,n}^{(2)}
   \stoch
   \begin{bmatrix}
    X_T - X_0 \\
    -\int_0^T  Y_s\, \dd X_s \\
    -\int_0^T X_s\, \dd X_s
   \end{bmatrix}
   =: \Bf_T^{(2)} ,
 \]
 as \ $n \to \infty$.
\ By Slutsky's lemma, using also Lemma \ref{LEMMA_CLSED_exist},  we conclude
 \begin{equation}\label{LSE_cont}
  \hbtheta_{T,n}^{\mathrm{approx}}
  = n \begin{bmatrix}
       \hc_{T,n} \\
       \hd_{T,n} \\
       \hdelta_{T,n} \\
       \hvare_{T,n} \\
       \hzeta_{T,n}
      \end{bmatrix}
  \stoch
  \begin{bmatrix}
   (\bG_T^{(1)})^{-1} \Bf_T^{(1)} \\
   (\bG_T^{(2)})^{-1} \Bf_T^{(2)}
  \end{bmatrix}
  =: \begin{bmatrix}
      \ha_T \\
      \hb_T \\
      \halpha_T \\
      \hbeta_T \\
      \hgamma_T
     \end{bmatrix}
  =: \hbtheta_T \qquad \text{as \ $n \to \infty$,}
 \end{equation}
 whenever the random matrices \ $\bG_T^{(1)}$ \ and \ $\bG_T^{(2)}$ \ are invertible.

\begin{Lem}\label{LEMMA_LSE_exist}
Let us consider the two-factor affine diffusion model \eqref{2dim_affine} with
 \ $a \in \RR_+$, \ $b, \alpha, \beta, \gamma \in \RR$, \ $\sigma_1 \in \RR_{++}$,
 \ $\sigma_2, \sigma_3 \in \RR_+$ \ and \ $\varrho \in [-1, 1]$ \ with a random initial value
 \ $(\eta_0, \zeta_0)$ \ independent of \ $(W_t, B_t, L_t)_{t\in\RR_+}$ \ satisfying
 \ $\PP(\eta_0 \in \RR_+) = 1$. 
\ Suppose that \ $(1 - \varrho^2) \sigma_2^2 + \sigma_3^2 > 0$.
\ Then for each \ $T \in \RR_{++}$, \ the random matrices \ $\bG_T^{(1)}$ \ and
 \ $\bG_T^{(2)}$ \ are invertible almost surely, and hence
 \ $\hbtheta_T$ \ given in \eqref{LSE_cont} exists almost surely.
Moreover, \ $\hbtheta_{T,n} \stoch \hbtheta_T$ \ as \ $n \to \infty$.
\end{Lem}

\noindent{\bf Proof.}
The aim of the following discussion is to show that the random matrix \ $\bG_T^{(1)}$ \ is
 almost surely strictly positive definite checking that for all
 \ $\bx \in \RR^2 \setminus \{\bzero\}$, \ we have \ $\bx^\top \bG_T^{(1)} \bx > 0$ \ almost
 surely.
Indeed, for all \ $\bx = (x_1, x_2)^\top \in \RR^2$,
 \[
   \begin{bmatrix} x_1 \\ x_2 \end{bmatrix}^\top
   \bG_T^{(1)}
   \begin{bmatrix} x_1 \\ x_2 \end{bmatrix}
   = \int_0^T
      \begin{bmatrix} x_1 \\ x_2 \end{bmatrix}^\top
      \begin{bmatrix} 1 \\ -Y_s \end{bmatrix}
      \begin{bmatrix} 1 \\ -Y_s \end{bmatrix}^\top
      \begin{bmatrix} x_1 \\ x_2 \end{bmatrix}
      \dd s
   = \int_0^T (x_1 - x_2 Y_s)^2 \, \dd s
   \geq 0 ,
 \]
 and \ $\bx^\top \bG_T^{(1)} \bx = 0$ \ if and only if \ $x_1 - x_2 Y_s = 0$ \ for all
 \ $s \in [0, T]$, \ which happens with probability 0, since
 \ $\bx = (x_1, x_2)^\top \ne \bzero$ \ and, for each \ $s \in (0, T]$, \ the distribution
 \ of \ $Y_s$ \ is absolutely continuous, since the conditional distribution of \ $Y_s$
 \ given \ $Y_0$ \ is absolutely continuous, see, eg., Ben Alaya and Kebaier
 \cite[Proof of Proposition 2]{BenKeb2} and Ikeda and Watanabe \cite[page 222]{IkeWat}.

In a similar way, the stochastic matrix \ $\bG_T^{(2)}$ \ is almost surely strictly positive
 definite, since for all \ $\bx \in \RR^3 \setminus \{\bzero\}$, \ we have
 \ $\bx^\top \bG_T^{(2)} \bx > 0$ \ almost surely.
Indeed, for all \ $\bx = (x_1, x_2, x_3)^\top \in \RR^3$,
 \[
   \begin{bmatrix} x_1 \\ x_2 \\ x_3 \end{bmatrix}^\top
   \bG_T^{(2)}
   \begin{bmatrix} x_1 \\ x_2 \\ x_3 \end{bmatrix}
   = \int_0^T
      \begin{bmatrix} x_1 \\ x_2 \\ x_3 \end{bmatrix}^\top
      \begin{bmatrix} 1 \\ -Y_s \\ -X_s \end{bmatrix}
      \begin{bmatrix} 1 \\ -Y_s \\ -X_s \end{bmatrix}^\top
      \begin{bmatrix} x_1 \\ x_2 \\ x_3 \end{bmatrix}
      \dd s
   = \int_0^T (x_1 - x_2 Y_s - x_3 X_s)^2 \, \dd s
   \geq 0 ,
 \]
 and \ $\bx^\top \bG_T^{(2)} \bx = 0$ \ if and only if \ $x_1 - x_2 Y_s - x_3 X_s = 0$ \ for
 all \ $s \in [0, T]$, \ which happens with probability 0, since
 \ $\bx = (x_1, x_2, x_3)^\top \ne \bzero$ \ and, for each
 \ $s \in (0, T]$, \ the distribution \ of \ $(Y_s, X_s)$ \ is absolutely continuous,
 because, as in the proof of part (b) in the proof of Theorem \ref{Thm_ergodic1} in Bolyog
 and Pap \cite{BolPap}, the conditional distribution of \ $(Y_s, X_s)$ \ given \ $(Y_0, X_0)$
 \ is absolutely continuous.

Next we are going to show \ $\hbtheta_{T,n} \stoch \hbtheta_T$ \ as \ $n \to \infty$.
\ The function \ $g_n$ \ introduced in \eqref{gn} admits an inverse
 \ $g_n^{-1} : \RR \times (-\infty, 1) \times \RR^2 \times (-\infty, 1) \to \RR^5$
 \ satisfying
 \[
   g_n^{-1}(c, d, \delta, \vare, \zeta) = (a, b, \alpha, \beta, \gamma)
 \]
 with
 \begin{gather*}
  b = -n \log(1 - d) , \qquad
  a = \frac{c}{\int_0^{\frac{1}{n}} \ee^{-bw} \, \dd w} , \qquad
  \gamma = -n \log(1 - \zeta) , \\
  \beta
  = \frac{\vare}{\int_0^{\frac{1}{n}} \ee^{(\gamma-b)w-\frac{\gamma}{n}} \, \dd w} , \qquad
  \alpha
  = \frac{\delta
          +a\beta
           \int_0^{\frac{1}{n}}
            \ee^{\gamma w-\frac{\gamma}{n}}
            \bigl(\int_0^w \ee^{-b(w-v)} \, \dd v\bigr) \dd w}
         {\int_0^{\frac{1}{n}} \ee^{-\gamma w} \, \dd w} .
 \end{gather*}
Convergence \eqref{LSE_cont} yields
 \ $(\hc_{T,n}, \hd_{T,n}, \hdelta_{T,n}, \hvare_{T,n}, \hzeta_{T,n}) \stoch \bzero$ \ as
 \ $n \to \infty$, \ hence \ $\hd_{T,n} \in (-\infty, 1)$ \ and
 \ $\hzeta_{T,n} \in (-\infty, 1)$ \ with probability tending to one as \ $n \to \infty$.
Consequently,
 \ $g_n^{-1}(\hc_{T,n}, \hd_{T,n}, \hdelta_{T,n}, \hvare_{T,n}, \hzeta_{T,n})
    =\hbtheta_{T,n}$
 \ with probability tending to one as \ $n \to \infty$.
\ We have
 \[
   \hb_{T,n} = -n \log(1 - \hd_{T,n}) = n \hd_{T,n} h_1(\hd_{T,n})
 \]
 with probability tending to one as \ $n \to \infty$, \ where the continuous function
 \ $h_1 : (-\infty, 1) \to \RR$ \ is given by
 \[
   h_1(x)
   := \begin{cases}
       - \frac{1}{x} \log(1 - x) & \text{if \ $x \ne 0$,} \\
       1 & \text{if \ $x = 0$.}
      \end{cases}
 \]
By \eqref{LSE_cont}, we have \ $n \hd_{T,n} \stoch \hb_T$ \ and \ $\hd_{T,n} \stoch 0$, \ thus
 we obtain \ $h_1(\hd_{T,n}) \stoch h_1(0) = 1$, \ and hence \ $\hb_{T,n} \stoch \hb_T$ \ as
 \ $n \to \infty$. 

Moreover,
 \[
   \ha_{T,n}
   = \frac{\hc_{T,n}}{\int_0^{\frac{1}{n}} \ee^{-\hb_{T,n}w} \, \dd w}
   = \frac{n\hc_{T,n}}{n\int_0^{\frac{1}{n}} \ee^{-\hb_{T,n}w} \, \dd w}
   = \frac{n\hc_{T,n}}{\int_0^1 \exp\bigl\{-n^{-1}\hb_{T,n}v\bigr\} \, \dd v}
   = \frac{n\hc_{T,n}}{h_2(n^{-1}\hb_{T,n})}
 \]
 with probability tending to one as \ $n \to \infty$, \ where the continuous function
 \ $h_2 : \RR \to \RR$ \ is given by
 \[
   h_2(x)
   := \int_0^1 \ee^{-xv} \, \dd v
   = \begin{cases}
      \frac{1-\ee^{-x}}{x} & \text{if \ $x \ne 0$,} \\
      1 & \text{if \ $x = 0$.}
     \end{cases}
 \]
We have already showed \ $\hb_{T,n} \stoch \hb_T$, \ yielding \ $n^{-1} \hb_{T,n} \stoch 0$,
 \ and hence \ $h_2(n^{-1} \hb_{T,n}) \stoch h_2(0) = 1$ \ as \ $n \to \infty$. 
\ By \eqref{LSE_cont}, we have \ $n \hc_{T,n} \stoch \ha_T$, \ thus we obtain
 \ $\ha_{T,n} \stoch \ha_T$ \ as \ $n \to \infty$. 

In a similar way,
 \[
   \hgamma_{T,n} = -n \log(1 - \hzeta_{T,n}) = n \hzeta_{T,n} h_1(\hzeta_{T,n}) 
 \]
 with probability tending to one as \ $n \to \infty$.
\ By \eqref{LSE_cont}, we have \ $n \hzeta_{T,n} \stoch \hgamma_T$ \ and
 \ $\hzeta_{T,n} \stoch 0$, \ thus we obtain \ $h_1(\hzeta_{T,n}) \stoch h_1(0) = 1$, \ and
 hence \ $\hgamma_{T,n} \stoch \hgamma_T$ \ as \ $n \to \infty$. 

Further, 
 \[
   \hbeta_{T,n}
   = \frac{\hvare_{T,n}}
          {\int_0^{\frac{1}{n}}
            \ee^{(\hgamma_{T,n}-\hb_{T,n})w-\frac{\hgamma_{T,n}}{n}} \, \dd w}
   = \frac{n\hvare_{T,n}\ee^{\frac{\hgamma_{T,n}}{n}}}{h_2(n^{-1}(\hb_{T,n}-\hgamma_{T,n}))}
 \]
 with probability tending to one as \ $n \to \infty$.
\ We have already showed \ $\hb_{T,n} \stoch \hb_T$ \ and \ $\hgamma_{T,n} \stoch \hgamma_T$,
 \ yielding \ $n^{-1} \hb_{T,n} \stoch 0$ \ and \ $n^{-1} \hgamma_{T,n} \stoch 0$, \ and hence 
 \ $\ee^{\frac{\hgamma_{T,n}}{n}} \stoch 1$ \ and
 \ $h_2(n^{-1} (\hb_{T,n} - \hgamma_{T,n})) \stoch h_2(0) = 1$ \ as \ $n \to \infty$. 
\ By \eqref{LSE_cont}, we have \ $n \hvare_{T,n} \stoch \hbeta_T$, \ thus we obtain
 \ $\hbeta_{T,n} \stoch \hbeta_T$ \ as \ $n \to \infty$. 

Finally,
 \[
   \halpha_{T,n}
   = \frac{\hdelta_{T,n}
           +\ha_{T,n}\hbeta_{T,n}
            \int_0^{\frac{1}{n}}
             \ee^{\hgamma_{T,n} w-\frac{\hgamma_{T,n}}{n}}
             \bigl(\int_0^w \ee^{-b(w-v)} \, \dd v\bigr) \dd w}
          {\int_0^{\frac{1}{n}} \ee^{-\hgamma_{T,n}w} \, \dd w}
   = \frac{n\hdelta_{T,n}+\ha_{T,n}\hbeta_{T,n}\ee^{-\frac{\hgamma_{T,n}}{n}}I_{T,n}}
          {h_2(n^{-1}\hgamma_{T,n})}
 \]
 with probability tending to one as \ $n \to \infty$, \ where
 \[
   I_{T,n}
   = n \int_0^{\frac{1}{n}}
        \ee^{\hgamma_{T,n}w} \bigl(\int_0^w \ee^{-b(w-v)} \, \dd v\bigr) \dd w
   \leq \frac{1}{n} \ee^{\frac{|\hgamma_{T,n}|}{n}} \ee^{\frac{|\hb_{T,n}|}{n}} .
 \]
We have already showed \ $\ha_{T,n} \stoch \ha_T$, \ $\hb_{T,n} \stoch \hb_T$,
 \ $\hbeta_{T,n} \stoch \hbeta_T$ \ and \ $\hgamma_{T,n} \stoch \hgamma_T$, \ yielding
 \ $n^{-1} \hb_{T,n} \stoch 0$ \ and \ $n^{-1} \hgamma_{T,n} \stoch 0$, \ and hence 
 \ $h_2(n^{-1} \hgamma_{T,n}) \stoch h_2(0) = 1$,
 \ $\ee^{\frac{\hgamma_{T,n}}{n}} \stoch 1$, \ $\ee^{\frac{|\hgamma_{T,n}|}{n}} \stoch 1$
 \ and \ $\ee^{\frac{|\hb_{T,n}|}{n}} \stoch 1$, \ implying \ $I_{T,n} \stoch 0$ \ as
 \ $n \to \infty$. 
\ By \eqref{LSE_cont}, we have \ $n \hdelta_{T,n} \stoch \halpha_T$, \ thus we
 obtain \ $\halpha_{T,n} \stoch \halpha_T$ \ as \ $n \to \infty$. 
\proofend

Using the SDE \eqref{2dim_affine} and Corollary 3.2.20 in Karatzas and Shreve \cite{KarShr},
 one can check that
 \begin{equation}\label{LSE_cont_true}
  \hbtheta_T - \btheta
  = \begin{bmatrix}
     \ha_T - a \\
     \hb_T -b \\
     \halpha_T - \alpha\\
     \hbeta_T - \beta \\
     \hgamma_T - \gamma
    \end{bmatrix}
  = \begin{bmatrix}
     (\bG_T^{(1)})^{-1} \bh_T^{(1)} \\
     (\bG_T^{(2)})^{-1} \bh_T^{(2)}
    \end{bmatrix}
  = \bG_T^{-1} \bh_T
 \end{equation}
 on the event where the random matrices \ $\bG_T^{(1)}$ \ and \ $\bG_T^{(2)}$ \ are 
 invertible, where
 \[
   \bG_T
   := \begin{bmatrix} \bG_T^{(1)} & \bzero \\ \bzero & \bG_T^{(2)} \end{bmatrix} , \qquad
   \bh_T:= \begin{bmatrix} \bh_T^{(1)} \\ \bh_T^{(2)} \end{bmatrix} ,
 \]
 with
 \[
   \bh_T^{(1)}
   := \sigma_1 \int_0^T \sqrt{Y_s} \begin{bmatrix} 1 \\ - Y_s \end{bmatrix} \dd W_s , \qquad
   \bh_T^{(2)}
   := \int_0^T
       \begin{bmatrix} 1 \\ - Y_s \\ - X_s \end{bmatrix}
       (\sigma_2 \sqrt{Y_s} \, \dd \tW_s + \sigma_3 \, \dd L_s) ,
 \]
 where 
 \begin{equation}\label{tW}
  \tW_s := \varrho W_s + \sqrt{1-\varrho^2} B_s , \qquad s \in \RR_+ ,
 \end{equation}
 is a standard Wiener process, independent of \ $L$.
\ Indeed,
 \[
   \hbtheta_T = \bG_T^{-1} \Bf_T \qquad \text{with} \qquad
   \Bf_T:= \begin{bmatrix} \Bf_T^{(1)} \\ \Bf_T^{(2)} \end{bmatrix} ,
 \]
 hence \ $\hbtheta_T - \btheta = \bG_T^{-1} (\Bf_T - \bG_T \btheta)$, \ where
 \[
   \Bf_T - \bG_T \btheta
   = \begin{bmatrix}
      \Bf_T^{(1)} - \bG_T^{(1)} \begin{bmatrix} a \\ b \end{bmatrix} \\
      \Bf_T^{(2)} - \bG_T^{(2)} \begin{bmatrix} \alpha \\ \beta \\ \gamma \end{bmatrix}
     \end{bmatrix} .
 \]
Using the first equation of SDE \eqref{2dim_affine} and Corollary 3.2.20 in Karatzas and
 Shreve \cite{KarShr}, we obtain
 \begin{align*}
  \Bf_T^{(1)} - \bG_T^{(1)} \begin{bmatrix} a \\ b \end{bmatrix}
  &= \int_0^T \begin{bmatrix} 1 \\ -Y_s \end{bmatrix} \dd Y_s
     - \int_0^T
        \begin{bmatrix} 1 \\ -Y_s \end{bmatrix}
        \begin{bmatrix} 1 \\ -Y_s \end{bmatrix}^\top
        \dd s
       \begin{bmatrix} a \\ b \end{bmatrix} \\
  &= \int_0^T
      \begin{bmatrix} 1 \\ -Y_s \end{bmatrix}
      \Biggl(\dd Y_s
             - \begin{bmatrix} 1 \\ -Y_s \end{bmatrix}^\top
                \begin{bmatrix} a \\ b \end{bmatrix}
                \dd s\Biggr)
   = \int_0^T
      \begin{bmatrix} 1 \\ -Y_s \end{bmatrix}
      \bigl(\dd Y_s - (a - b Y_s) \, \dd s\bigr) \\
  &= \int_0^T
      \begin{bmatrix} 1 \\ -Y_s \end{bmatrix}
      \sigma_1 \sqrt{Y_s} \, \dd W_s .
 \end{align*}
In a similar way,
 \begin{align*}
  &\Bf_T^{(2)} - \bG_T^{(2)} \begin{bmatrix} \alpha \\ \beta \\ \gamma \end{bmatrix}
   = \int_0^T \begin{bmatrix} 1 \\ -Y_s \\ -X_s \end{bmatrix} \dd X_s
     - \int_0^T
        \begin{bmatrix} 1 \\ -Y_s \\ -X_s \end{bmatrix}
        \begin{bmatrix} 1 \\ -Y_s \\ -X_s \end{bmatrix}^\top
        \dd s
       \begin{bmatrix} \alpha \\ \beta \\ \gamma \end{bmatrix} \\
  &= \int_0^T
      \begin{bmatrix} 1 \\ -Y_s \\ -X_s \end{bmatrix}
      \left(\dd X_s
            - \begin{bmatrix} 1 \\ -Y_s \\ -X_s \end{bmatrix}^\top
              \begin{bmatrix} \alpha \\ \beta \\ \gamma \end{bmatrix}
              \dd s\right)
   = \int_0^T
      \begin{bmatrix} 1 \\ -Y_s \\ -X_s \end{bmatrix}
       \bigl(\dd X_s - (\alpha - \beta Y_s - \gamma X_s) \, \dd s\bigr) \\
  &= \int_0^T
      \begin{bmatrix} 1 \\ -Y_s \\ -X_s \end{bmatrix}
      \bigl(\sigma_2 \sqrt{Y_s} [\varrho \dd W_s + \sqrt{1-\varrho^2} B_s]
            + \sigma_3 \dd L_s\bigr) .
 \end{align*}

\section{Consistency of CLSE}
\label{section_LSE}

First we consider the case of subcritical Heston models, i.e., when
 \ $b \in \RR_{++}$.

\begin{Thm}\label{Thm_LSE_cons_sub}
Let us consider the two-factor affine diffusion model \eqref{2dim_affine} with
 \ $a, b \in \RR_{++}$, \ $\alpha, \beta \in \RR$, \ $\gamma \in \RR_{++}$,
 \ $\sigma_1 \in \RR_{++}$, \ $\sigma_2, \sigma_3 \in \RR_+$ \ and
 \ $\varrho \in [-1, 1]$ \ with a random initial value \ $(\eta_0, \zeta_0)$
 \ independent of \ $(W_t, B_t, L_t)_{t\in\RR_+}$ \ satisfying
 \ $\PP(\eta_0 \in \RR_+) = 1$. 
\ Suppose that \ $(1 - \varrho^2) \sigma_2^2 + \sigma_3^2 > 0$.
\ Then the CLSE of \ $\btheta = (a, b, \alpha, \beta, \gamma)^\top$ \ is strongly
 consistent, i.e.,
 \ $\hbtheta_T
    = \bigl(\ha_T, \hb_T, \halpha_T, \hbeta_T, \hgamma_T\bigr)^\top
    \as \btheta = (a, b, \alpha, \beta, \gamma)^\top$
 \ as \ $T \to \infty$. 
\end{Thm}

\noindent{\bf Proof.}
By \eqref{LSE_cont_true}, we have
 \begin{equation}\label{LSE_cont_trueT}
  \hbtheta_T - \btheta = (T^{-1} \bG_T)^{-1} (T^{-1} \bh_T) 
 \end{equation}
 on the event, where the random matrix \ $\bG_T$ \ is invertible, which has propapility 1,
 see Lemma \ref{LEMMA_LSE_exist}. 
 
By Theorem \ref{Thm_ergodic2}, we obtain
 \begin{equation}\label{bG}
  T^{-1} \bG_T \as \EE(\bG_\infty) \qquad \text{as \ $T \to \infty$,}
 \end{equation}
 where
 \begin{equation}\label{bGinfty}
  \bG_\infty
  := \begin{bmatrix} \bG_\infty^{(1)} & \bzero \\ \bzero & \bG_\infty^{(2)} \end{bmatrix}
 \end{equation}
 with
 \[
   \bG_\infty^{(1)}
   := \begin{bmatrix} 1 & -Y_\infty \\ -Y_\infty & Y_\infty^2 \end{bmatrix} ,
   \qquad
   \bG_\infty^{(2)}
   := \begin{bmatrix}
       1 & -Y_\infty & -X_\infty \\
       -Y_\infty & Y_\infty^2 & Y_\infty X_\infty \\
       -X_\infty & Y_\infty X_\infty & X_\infty^2
      \end{bmatrix} ,
 \]
 where the random vector \ $(Y_\infty, X_\infty)$ \ is given by Theorem \ref{Thm_ergodic1}, 
 since, by Theorem \ref{Thm_moments}, the entries of \ $\EE(\bG_\infty)$ \ exist and finite. 

The matrix \ $\EE(\bG_\infty^{(1)})$ \ is strictly positive definite, since for all
 \ $\bx \in \RR^2 \setminus \{\bzero\}$, \ we have
 \ $\bx^\top \EE(\bG_\infty^{(1)}) \bx > 0$.
\ Indeed, for all \ $\bx = (x_1, x_2)^\top \in \RR^2 \setminus \{\bzero\}$,
 \[
   \begin{bmatrix} x_1 \\ x_2 \end{bmatrix}^\top
   \EE(\bG_\infty^{(1)})
   \begin{bmatrix} x_1 \\ x_2 \end{bmatrix}
   = \EE\bigl[(x_1 - x_2 Y_\infty)^2\bigr]
   > 0 ,
 \]
 since, by Theorem \ref{Thm_ergodic2}, the distribution \ of \ $Y_\infty$ \ is absolutely 
 continuous, hence \ $x_1 - x_2 Y_\infty \ne 0$ \ with probability 1.
In a similar way, the matrix \ $\EE(\bG_\infty^{(2)})$ \ is strictly positive definite,
 since for all \ $\bx \in \RR^3 \setminus \{\bzero\}$, \ we have
 \ $\bx^\top \EE(\bG_\infty^{(2)}) \bx > 0$.
\ Indeed, for all \ $\bx = (x_1, x_2, x_3)^\top \in \RR^3 \setminus \{\bzero\}$,
 \[
   \begin{bmatrix} x_1 \\ x_2 \\ x_3 \end{bmatrix}^\top
   \EE(\bG_\infty^{(2)})
   \begin{bmatrix} x_1 \\ x_2 \\ x_3 \end{bmatrix}
   = \EE\bigl[(x_1 - x_2 Y_\infty - x_3 Z_\infty)^2\bigr]
   > 0 ,
 \]
 since, by Theorem \ref{Thm_ergodic2}, the distribution \ of \ $(Y_\infty, X_\infty)$ \ is
 absolutely continuous, hence \ $x_1 - x_2 Y_\infty - x_3 X_\infty \ne 0$ \ with probability
 1. 
Thus the matrices \ $\EE(\bG_\infty^{(1)})$ \ and \ $\EE(\bG_\infty^{(2)})$ \ are invertible, 
 whence we conclude 
 \begin{equation}\label{bG^{-1}}
  (T^{-1} \bG_T)^{-1}
  \as \begin{bmatrix}
       [\EE(\bG_\infty^{(1)})]^{-1} & \bzero \\ 
       \bzero & [\EE(\bG_\infty^{(2)})]^{-1}
      \end{bmatrix}
  = [\EE(\bG_\infty)]^{-1}
  \qquad \text{as \ $T \to \infty$.}
 \end{equation} 
The aim of the next discussion is to show convergence
 \begin{equation}\label{bh}
  T^{-1} \bh_T \as \bzero \qquad \text{as \ $T \to \infty$.}
 \end{equation}
We have
 \begin{gather*}
  \frac{1}{T} \int_0^T \sqrt{Y_s} \, \dd W_s
  = \frac{1}{T} \int_0^T Y_s \, \dd s \cdot
    \frac{\int_0^T \sqrt{Y_s} \, \dd W_s}{\int_0^T Y_s \, \dd s}
  \as 0 \qquad \text{as \ $T \to \infty$.}
 \end{gather*}
Indeed, we have already proved
 \[
   \frac{1}{T} \int_0^T Y_s \, \dd s \as \EE(Y_\infty) = \frac{a}{b} \in \RR_{++} \qquad 
   \text{as \ $T \to \infty$,}
 \]
 and the strong law of large numbers for continuous local martingales (see, e.g., Theorem
 \ref{DDS_stoch_int}) implies
 \[
   \frac{\int_0^T \sqrt{Y_s} \, \dd W_s}{\int_0^T Y_s \, \dd s} \as 0 \qquad 
   \text{as \ $T \to \infty$,}
 \]
 since we have
 \[
   \int_0^T Y_s \, \dd s = T \cdot \frac{1}{T} \int_0^T Y_s \, \dd s \as \infty \qquad 
   \text{as \ $T \to \infty$.}
 \]
Further,
 \[
   \frac{1}{T}
   \int_0^T
    (\sigma_2 \sqrt{Y_s} \, \dd \tW_s + \sigma_3 \, \dd L_s)
   = \frac{1}{T} \int_0^T (\sigma_2^2 Y_s + \sigma_3^2) \dd s \cdot
     \frac{\int_0^T (\sigma_2 \sqrt{Y_s} \, \dd \tW_s + \sigma_3 \, \dd L_s)}
          {\int_0^T (\sigma_2^2 Y_s + \sigma_3^2) \, \dd s}
   \as 0
 \]
 as \ $T \to \infty$.
\ Indeed, we have already proved
 \[
   \frac{1}{T} \int_0^T (\sigma_2^2 Y_s + \sigma_3^2) \, \dd s
   \as \EE(\sigma_2^2 Y_\infty + \sigma_3^2)
   = \sigma_2^2 \frac{a}{b} + \sigma_3^2 \in \RR_{++} \qquad 
   \text{as \ $T \to \infty$,}
 \]
 and the strong law of large numbers for continuous local martingales (see, e.g., Theorem
 \ref{DDS_stoch_int}) implies
 \[
   \frac{\int_0^T (\sigma_2 \sqrt{Y_s} \, \dd \tW_s + \sigma_3 \, \dd L_s)}
        {\int_0^T (\sigma_2^2 Y_s + \sigma_3^2) \, \dd s}
   \as 0 \qquad 
   \text{as \ $T \to \infty$,}
 \]
 since we have
 \[
   \int_0^T (\sigma_2^2 Y_s + \sigma_3^2) \, \dd s
   = T \cdot \frac{1}{T} \int_0^T (\sigma_2^2 Y_s + \sigma_3^2) \, \dd s \as \infty
   \qquad \text{as \ $T \to \infty$.}
 \] 
One can check
 \begin{gather*}
  \frac{1}{T} \int_0^T Y_s \sqrt{Y_s} \, \dd W_s \as 0 , \\
  \frac{1}{T} \int_0^T Y_s (\sigma_2 \sqrt{Y_s} \, \dd \tW_s + \sigma_3 \, \dd L_s) \as 0 ,
  \qquad 
  \frac{1}{T} \int_0^T X_s (\sigma_2 \sqrt{Y_s} \, \dd \tW_s + \sigma_3 \, \dd L_s) \as 0
 \end{gather*}
 as \ $T \to \infty$ \ in the same way, since
 \begin{gather*}
  \frac{1}{T} \int_0^T Y_s^3 \, \dd s \as \EE(Y_\infty^3) \in \RR_{++} , \\
  \frac{1}{T} \int_0^T Y^2_s (\sigma_2^2 Y_s + \sigma_3^2) \, \dd s
  \as \EE\bigl[Y^2_s (\sigma_2^2 Y_s + \sigma_3^2)\bigr] \in \RR_{++}  , \\
  \frac{1}{T} \int_0^T X^2_s (\sigma_2^2 Y_s + \sigma_3^2) \, \dd s
  \as \EE\bigl[X^2_s (\sigma_2^2 Y_s + \sigma_3^2)\bigr] \in \RR_{++} 
 \end{gather*}
 as \ $T \to \infty$.
\ Consequently, we conclude \eqref{bh}.
Finally, by \eqref {bG^{-1}} and \eqref {bh}, we obtain the statement.
\proofend

In order to handle supercritical two-factor affine diffusion models when
 \ $b \in \RR_{--}$, \ we need the following integral version of the Toeplitz Lemma, due to
 Dietz and Kutoyants \cite{DieKut}.

\begin{Lem}\label{int_Toeplitz}
Let \ $\{\varphi_T : T \in \RR_+\}$ \ be a family of probability measures on
 \ $\RR_+$ \ such that \ $\varphi_T([0,T]) = 1$ \ for all \ $T \in \RR_+$,
 \ and \ $\lim_{T\to\infty} \varphi_T([0,K]) = 0$ \ for all \ $K \in \RR_{++}$.
\ Then for every bounded and measurable function \ $f : \RR_+ \to \RR$ \ for
 which the limit \ $f(\infty) := \lim_{t\to\infty} f(t)$ \ exists, we have
 \[
   \lim_{T\to\infty} \int_0^\infty f(t) \, \varphi_T(\dd t) = f(\infty) .
 \]
\end{Lem}

As a special case, we have the following integral version of the Kronecker Lemma, see
 K\"uchler and S{\o}rensen \cite[Lemma B.3.2]{KS}.

\begin{Lem}\label{int_Kronecker}
Let \ $a : \RR_+ \to \RR_+$ \ be a measurable function.
Put \ $b(T) := \int_0^T a(t) \, \dd t$, \ $T \in \RR_+$.
\ Suppose that \ $\lim_{T\to\infty} b(T) = \infty$.
\ Then for every bounded and measurable function \ $f : \RR_+ \to \RR$ \ for
 which the limit \ $f(\infty) := \lim_{t\to\infty} f(t)$ \ exists, we have
 \[
   \lim_{T\to\infty} \frac{1}{b(T)} \int_0^T a(t) f(t) \, \dd t = f(\infty) .
 \]
\end{Lem}

Next we present an auxiliary lemma in the supercritical case about the asymptotic behavior of
 \ $Y_t$ \ as \ $t \to \infty$.

\begin{Lem}\label{martconvY}
Let us consider the two-factor affine diffusion model \eqref{2dim_affine} with
 \ $a \in \RR_+$, \ $b \in \RR_{--}$, \ $\alpha, \beta, \gamma \in \RR$,
 \ $\sigma_1, \sigma_2, \sigma_3 \in \RR_+$ \ and \ $\varrho \in [-1, 1]$ \ with a random
 initial value \ $(\eta_0, \zeta_0)$ \ independent of \ $(W_t, B_t, L_t)_{t\in\RR_+}$
 \ satisfying \ $\PP(\eta_0 \in \RR_+) = 1$. 
\ Then there exists a random variable \ $V_Y$ \ such that
 \begin{align}\label{lim_Y}
  \ee^{bt} Y_t \as V_Y \qquad \text{as \ $t \to \infty$}
 \end{align}
 with \ $\PP(V_Y \ne 0) = 1$, \ and, for each \ $k \in \NN$,
 \begin{align}\label{lim_intY}
  \ee^{kbt} \int_0^t Y_u^k \, \dd u \as -\frac{V_Y^k}{kb} \qquad \text{as \ $t \to \infty$.}
 \end{align}
\end{Lem}

\noindent{\bf Proof.}
By \eqref{SolutionY},
 \begin{align*}
  \EE( Y_t \mid \cF_s )
  = \EE( Y_t \mid Y_s )
  = \ee^{-b(t-s)} Y_s + a \int_s^t \ee^{-b(t-u)} \, \dd u
 \end{align*}
 for all \ $s, t \in \RR_+$ \ with \ $0 \leq s \leq t$.
\ Thus
 \[
   \EE( \ee^{bt} Y_t \mid \cF^Y_s )
   = \ee^{bs} Y_s + a \int_s^t \ee^{bu} \, \dd u
   \geq \ee^{bs} Y_s
 \]
 for all \ $s, t \in \RR_+$ \ with \ $0 \leq s \leq t$, \ consequently, the process
 \ $(\ee^{bt} Y_t)_{t\in\RR_+}$ \ is a non-negative submartingale with respect to the
 filtration \ $(\cF^Y_t)_{t\in\RR_+}$.
\ Moreover, \ $b \in \RR_{--}$ \ implies
 \[
   \EE(\ee^{bt} Y_t)
   = y_0 + a \int_0^t \ee^{bu} \, \dd u
   \leq y_0 + a \int_0^\infty \ee^{bu} \, \dd u
   = y_0 - \frac{a}{b}
   < \infty , \qquad t \in \RR_+ ,
 \]
 hence, by the submartingale convergence theorem, there exists a non-negative random variable
 \ $V_Y$ \ such that \eqref{lim_Y} holds.

The distribution of \ $V_Y$ \ coincides with the distribution of \ $\tcY_{-1/b}$, \ where
 \ $(\tcY_t)_{t\in\RR_+}$ \ is a CIR process given by the SDE
 \begin{align*}
  \dd \tcY_t = a \dd t + \sigma_1 \sqrt{\tcY_t} \, \dd \cW_t ,
  \qquad t \in \RR_+ ,
 \end{align*}
 with initial value \ $\tcY_0 = y_0$, \ where \ $(\cW_t)_{t\in\RR_+}$ \ is a standard Wiener
 process, see Ben Alaya and Kebaier \cite[Proposition 3]{BenKeb1}.
Consequently, \ $\PP(V_Y \in \RR_{++}) = 1$, \ since \ $\tcY_t$, \ $t \in \RR_{++}$, \ are
 absolutely continuous random variables.
 
If \ $\omega \in \Omega$ \ such that \ $\RR_+ \ni t \mapsto Y_t(\omega)$ \ is continuous and
 \ $\ee^{bt} Y_t(\omega) \to V_Y(\omega)$ \ as \ $t \to \infty$, \ then, by the integral
 Kronecker Lemma \ref{int_Kronecker} with \ $f(t) = \ee^{kbt} Y_t(\omega)^k$ \ and
 \ $a(t) = \ee^{-kbt}$, \ $t \in \RR_+$, \ we have
 \[
   \frac{1}{\int_0^t \ee^{-kbu} \, \dd u}
   \int_0^t \ee^{-kbu} (\ee^{kbu} Y_u(\omega)^k) \, \dd u
   \to V_Y(\omega)^k \qquad \text{as \ $t \to \infty$.}
 \]
Here \ $\int_0^t \ee^{-kbu} \, \dd u = - \frac{\ee^{-kbt} - 1}{kb}$, \ $t \in \RR_+$, \ thus
 we conclude the second convergence in \eqref{lim_intY}.
\proofend.

The next theorem states strong consistency of the CLSE of \ $b$ \ in the supercritical case.

\begin{Thm}\label{Thm_MLE_cons_super}
Let us consider the two-factor affine diffusion model \eqref{2dim_affine} with
 \ $a \in \RR_+$, \ $b \in \RR_{--}$, \ $\alpha, \beta, \gamma \in \RR$,
 \ $\sigma_1 \in \RR_{++}$, \ $\sigma_2, \sigma_3 \in \RR_+$ \ and \ $\varrho \in [-1, 1]$
 \ with a random initial value \ $(\eta_0, \zeta_0)$ \ independent of
 \ $(W_t, B_t, L_t)_{t\in\RR_+}$ \ satisfying \ $\PP(\eta_0 \in \RR_+) = 1$. 
\ Then the CLSE of \ $b$ \ is strongly consistent, i.e., \ $\hb_T \as b$ \ as
 \ $T \to \infty$.
\end{Thm}

\noindent{\bf Proof.}
By Lemma \ref{LEMMA_LSE_exist}, there exists a unique CLSE \ $\hb_T$ \ of \ $b$ \ for all
 \ $T \in \RR_{++}$ \ which has the form given in \eqref{LSE_cont}.
By Ito's formula,
 \[
   \int_0^T Y_s \, \dd Y_s
   = \frac{1}{2} ( Y_T^2 - Y_0^2) - \frac{1}{2} \sigma_1^2 \int_0^T Y_s \, \dd s , \qquad
   T \in \RR_+ ,
 \]
 hence, by \eqref{lim_Y} and \eqref{lim_intY}, we have
 \begin{align*}
  \hb_T
  &= \frac{(Y_T-Y_0)\int_0^T Y_s\,\dd s-T\int_0^T Y_s\,\dd Y_s} 
          {T\int_0^T Y_s^2\,\dd s-\bigl(\int_0^T Y_s\,\dd s\bigr)^2}
   = \frac{(Y_T-Y_0)\int_0^T Y_s\,\dd s-\frac{T}{2}(Y_T^2-Y_0^2)
           +\frac{T}{2}\sigma_1^2\int_0^T Y_s\,\dd s} 
          {T\int_0^T Y_s^2\,\dd s-\bigl(\int_0^T Y_s\,\dd s\bigr)^2} \\
  &= \frac{\frac{1}{T}\bigl(\ee^{bT}Y_T-\ee^{bT}Y_0\bigr)
           \bigl(\ee^{bT}\int_0^T Y_s\,\dd s \bigr)
           -\frac{1}{2}\bigl(\ee^{2bT}Y_T^2-\ee^{2bT}Y_0^2\bigr)
           +\frac{1}{2}\sigma_1^2\ee^{bT}\bigl(\ee^{bT}\int_0^T Y_s\,\dd s\bigr)}
          {\ee^{2bT}\int_0^T Y_s^2\,\dd s
           -\frac{1}{T}\bigl(\ee^{bT}\int_0^T Y_s\,\dd s\bigr)^2} \\
  &\as \frac{0(V_Y-0)\bigl(-\frac{V_Y}{b}\bigr)-\frac{1}{2}(V_Y^2-0)
             +\frac{1}{2}\sigma_1^2 0\bigl(-\frac{V_Y}{b}\bigr)}
            {-\frac{V_Y^2}{2b}-0\bigl(-\frac{V_Y}{b}\bigr)^2}
   = b
 \end{align*}
 as \ $T \to \infty$.
\proofend.

\begin{Rem}\label{Rem_comparison5}
For critical two-factor affine diffusion models, it will turn out that the CLSE of \ $a$
 \ and \ $\alpha$ \ are not even weakly consistent, but the CLSE of \ $b$, \ $\beta$ \ and
 \ $\gamma$ \ are weakly consistent, see Theorem \ref{Thm_LSE_crit}.
\proofend
\end{Rem}

\begin{Rem}\label{Rem_comparison6}
For supercritical two-factor affine diffusion models, it will turn out that the CLSE of
 \ $a$ \ and \ $\alpha$ \ are not even weakly consistent, but the CLSE of \ $\beta$ \ and
 \ $\gamma$ \ are weakly consistent, see Theorem \ref{Thm_LSE_super}.
\proofend
\end{Rem}

\section{Asymptotic behavior of CLSE: subcritical case}
\label{section_ALSE_subcritical}

\begin{Thm}\label{Thm_LSE_sub}
Let us consider the two-factor affine diffusion model \eqref{2dim_affine} with
 \ $a, b \in \RR_{++}$, \ $\alpha, \beta \in \RR$, \ $\gamma \in \RR_{++}$,
 \ $\sigma_1 \in \RR_{++}$, \ $\sigma_2, \sigma_3 \in \RR_+$ \ and \ $\varrho \in [-1, 1]$
 \ with a random initial value \ $(\eta_0, \zeta_0)$ \ independent of
 \ $(W_t, B_t, L_t)_{t\in\RR_+}$ \ satisfying \ $\PP(\eta_0 \in \RR_+) = 1$. 
\ Suppose that \ $(1 - \varrho^2) \sigma_2^2 + \sigma_3^2 > 0$.
\ Then the CLSE of \ $\btheta = (a, b, \alpha, \beta, \gamma)^\top$ \ is asymptotically
 normal, namely,
 \begin{align}\label{LSE_sub}
  T^{\frac{1}{2}} (\hbtheta_T - \btheta)
  \distr \cN_5(\bzero, [\EE(\bG_\infty)]^{-1} \EE(\tbG_\infty) [\EE(\bG_\infty)]^{-1}) \qquad
  \text{as \ $T \to \infty$,}
 \end{align}
 where \ $\bG_\infty$ \ is given in \eqref{bGinfty} and \ $\tbG_\infty$ \ has the form
 \[
    \begin{bmatrix}
       \sigma_1^2 Y_\infty & - \sigma_1^2 Y_\infty^2 & \varrho \sigma_1 \sigma_2 Y_\infty
        & - \varrho \sigma_1 \sigma_2 Y_\infty^2
        & - \varrho \sigma_1 \sigma_2 Y_\infty X_\infty \\
       - \sigma_1^2 Y_\infty^2 & \sigma_1^2 Y_\infty^3 & - \varrho \sigma_1 \sigma_2 Y_\infty^2
        & \varrho \sigma_1 \sigma_2 Y_\infty^3
        & \varrho \sigma_1 \sigma_2 Y_\infty^2 X_\infty \\
       \varrho \sigma_1 \sigma_2 Y_\infty & - \varrho \sigma_1 \sigma_2 Y_\infty^2
        & \sigma_2^2 Y_\infty + \sigma_3^2 & - (\sigma_2^2 Y_\infty + \sigma_3^2) Y_\infty
        & - (\sigma_2^2 Y_\infty + \sigma_3^2) X_\infty \\
       - \varrho \sigma_1 \sigma_2 Y_\infty^2 & \varrho \sigma_1 \sigma_2 Y_\infty^3
        & - (\sigma_2^2 Y_\infty + \sigma_3^2) Y_\infty
        & (\sigma_2^2 Y_\infty + \sigma_3^2) Y_\infty^2
        & (\sigma_2^2 Y_\infty + \sigma_3^2) Y_\infty X_\infty \\
       - \varrho \sigma_1 \sigma_2 Y_\infty X_\infty
        & \varrho \sigma_1 \sigma_2 Y_\infty^2 X_\infty
        & - (\sigma_2^2 Y_\infty + \sigma_3^2) X_\infty
        & (\sigma_2^2 Y_\infty + \sigma_3^2) Y_\infty X_\infty
        & (\sigma_2^2 Y_\infty + \sigma_3^2) X_\infty^2
      \end{bmatrix} ,
 \]
 where the random vector \ $(Y_\infty, X_\infty)$ \ is given by Theorem \ref{Thm_ergodic1}.
\end{Thm}

\noindent{\bf Proof.}
By \eqref{LSE_cont_true}, we have
 \begin{equation}\label{ThT}
  T^{\frac{1}{2}} (\hbtheta_T - \btheta)
  = (T^{-1} \bG_T)^{-1} (T^{-\frac{1}{2}} \bh_T) 
 \end{equation}
 on the event where \ $\bG_T$ \ is invertible, which holds almost surely, see Lemma
 \ref{LEMMA_LSE_exist}. 
By \eqref{bG^{-1}}, we have \ $(T^{-1} \bG_T)^{-1} \as [\EE(\bG_\infty)]^{-1}$ \ as
 \ $T \to \infty$.
\ The process \ $(\bh_t)_{t\in\RR_+}$ \ is a 5-dimensional continuous local martingale with
 quadratic variation process \ $\langle \bh \rangle_t = \tbG_t$, \ $t \in \RR_+$, \ where
 \[
   \tbG_t
   := \int_0^t
       \begin{bmatrix}
        \sigma_1^2 Y_s & - \sigma_1^2 Y_s^2 & \varrho \sigma_1 \sigma_2 Y_s
         & - \varrho \sigma_1 \sigma_2 Y_s^2 & - \varrho \sigma_1 \sigma_2 Y_s X_s \\
        - \sigma_1^2 Y_s^2 & \sigma_1^2 Y_s^3 & - \varrho \sigma_1 \sigma_2 Y_s^2
         & \varrho \sigma_1 \sigma_2 Y_s^3 & \varrho \sigma_1 \sigma_2 Y_s^2 X_s \\
        \varrho \sigma_1 \sigma_2 Y_s & - \varrho \sigma_1 \sigma_2 Y_s^2
         & \sigma_2^2 Y_s + \sigma_3^2 & - (\sigma_2^2 Y_s + \sigma_3^2) Y_s
         & - (\sigma_2^2 Y_s + \sigma_3^2) X_s \\
        - \varrho \sigma_1 \sigma_2 Y_s^2 & \varrho \sigma_1 \sigma_2 Y_s^3
         & - (\sigma_2^2 Y_s + \sigma_3^2) Y_s & (\sigma_2^2 Y_s + \sigma_3^2) Y_s^2
         & (\sigma_2^2 Y_s + \sigma_3^2) Y_s X_s \\
        - \varrho \sigma_1 \sigma_2 Y_s X_s & \varrho \sigma_1 \sigma_2 Y_s^2 X_s
         & - (\sigma_2^2 Y_s + \sigma_3^2) X_s & (\sigma_2^2 Y_s + \sigma_3^2) Y_s X_s
         & (\sigma_2^2 Y_s + \sigma_3^2) X_s^2
       \end{bmatrix}
       \dd s .
 \]
By Theorem \ref{Thm_ergodic2}, we obtain
 \begin{equation}\label{tbG}
  T^{-1} \tbG_T \as \EE(\tbG_\infty) \qquad \text{as \ $T \to \infty$,}
 \end{equation}
 since, by Theorem \ref{Thm_moments}, the entries of \ $\EE(\tbG_\infty)$ \ exist and finite. 
Using \eqref{tbG}, Theorem \ref{THM_Zanten} yields
 \ $T^{-\frac{1}{2}} \bh_T \distr \cN_5(\bzero, \EE(\tbG_\infty))$ \ as \ $T \to \infty$.
\ Hence, by \eqref{ThT} and by Slutsky's lemma,
 \[
   T^{\frac{1}{2}} (\hbtheta_T - \btheta)
   \distr [\EE(\bG_\infty)]^{-1} \cN_5\bigl(\bzero, \EE(\tbG_\infty)\bigr)
   =\cN_5\bigl(\bzero, 
               [\EE(\bG_\infty)]^{-1} \EE(\tbG_\infty)
               \bigl([\EE(\bG_\infty)]^{-1}\bigr)^\top\bigr)
 \]
 as \ $T \to \infty$. 
\proofend

\section{Asymptotic behavior of CLSE: critical case}
\label{section_ALSE_critical}

First we present an auxiliary lemma.

\begin{Lem}\label{Lem_intf_cont}
If \ $(\cY_t, \cX_t)_{t\in\RR_+}$ \ and \ $(\tcY_t, \tcX_t)_{t\in\RR_+}$ \ are continuous
 semimartingales with \ $(\cY_t, \cX_t)_{t\in\RR_+} \distre (\tcY_t, \tcX_t)_{t\in\RR_+}$,
 \ then
 \begin{align*}
  &\biggl(\cY_1, \cX_1, \int_0^1 \cX_s \, \dd\cY_s, \int_0^1 \cY_s^k \cX_s^\ell \, \dd s
          : k, \ell \in \ZZ_+ , k + \ell \leq n\biggr) \\
  &\distre
   \biggl(\tcY_1, \tcX_1, \int_0^1 \tcX_s \, \dd\tcY_s,
          \int_0^1 \tcY_s^k \tcX_s^\ell \, \dd s
          : k, \ell \in \ZZ_+ , k + \ell \leq n\biggr)
 \end{align*}
 for each \ $n \in \NN$. 
\end{Lem}

\noindent{\bf Proof.}
By Proposition I.4.44 in Jacod and Shiryaev \cite{JSh} with the Riemann sequence of
 deterministic subdivisions \ $\left(\frac{i}{n} \land T\right)_{i\in\NN}$, \ $n \in \NN$,
 \ we have
 \[
   \sum_{i=1}^n \cX_{\frac{i-1}{n}} (\cY_{\frac{i}{n}} - \cY_{\frac{i-1}{n}})
   \stoch
   \int_0^1 \cX_s \, \dd\cY_s , \qquad
   \frac{1}{n} \sum_{i=1}^n \cY_{\frac{i}{n}}^k \cX_{\frac{i}{n}}^\ell
   \stoch
   \int_0^1 \cY_s^k \cX_s^\ell \, \dd s
 \]
 as \ $n \to \infty$ \ for each \ $k, \ell \in \ZZ_+$, \ and similar convergences hold for
 \ $(\tcY_t, \tcX_t)_{t\in\RR_+}$.
\ The assumption implies
 \[
   \sum_{i=1}^n \cX_{\frac{i-1}{n}} (\cY_{\frac{i}{n}} - \cY_{\frac{i-1}{n}})
   \distre
   \sum_{i=1}^n \tcX_{\frac{i-1}{n}} (\tcY_{\frac{i}{n}} - \tcY_{\frac{i-1}{n}}) , \qquad
   \frac{1}{n} \sum_{i=1}^n \cY_{\frac{i}{n}}^k \cX_{\frac{i}{n}}^\ell
   \distre
   \frac{1}{n} \sum_{i=1}^n \tcY_{\frac{i}{n}}^k \tcX_{\frac{i}{n}}^\ell
 \]
 for each \ $n \in \NN$ \ and \ $k, \ell \in \ZZ_+$, \ hence we obtain the statement.
\proofend

\begin{Thm}\label{Thm_LSE_crit}
Let us consider the two-factor affine diffusion model \eqref{2dim_affine} with
 \ $a \in \RR_+$, \ $b = 0$, \ $\alpha \in \RR$, \ $\beta = 0$, \ $\gamma = 0$,
 \ $\sigma_1, \sigma_2, \sigma_3 \in \RR_+$ \ and \ $\varrho \in [-1, 1]$ \ with a random
 initial value \ $(\eta_0, \zeta_0)$ \ independent of \ $(W_t, B_t, L_t)_{t\in\RR_+}$
 \ satisfying \ $\PP(\eta_0 \in \RR_+) = 1$. 
\ Suppose that \ $(1 - \varrho^2) \sigma_2^2 + \sigma_3^2 > 0$.
\ Then
 \begin{align}\label{LSE_crit}
  \begin{bmatrix}
   \ha_T - a \\
   T \hb_T \\
   \halpha_T - \alpha \\
   T \hbeta_T \\
   T \hgamma_T
  \end{bmatrix}
  \distr 
  \begin{bmatrix}
   \left(\int_0^1
          \begin{bmatrix} 1 \\ - \cY_s \end{bmatrix}
          \begin{bmatrix} 1 \\ - \cY_s \end{bmatrix}^\top
          \dd s\right)^{-1}
   \begin{bmatrix}
    \cY_1 - a \\
    -\frac{1}{2} \cY_1^2 + \bigl(a + \frac{\sigma_1^2}{2}\bigr) \int_0^1 \cY_s \, \dd s
   \end{bmatrix} \\
   \left(\int_0^1
          \begin{bmatrix} 1 \\ - \cY_s \\ - \cX_s \end{bmatrix}
          \begin{bmatrix} 1 \\ - \cY_s \\ - \cX_s \end{bmatrix}^\top
          \dd s\right)^{-1}
   \begin{bmatrix}
    \cX_1 - \alpha \\
    - \cY_1 \cX_1
     + (\alpha + \varrho \sigma_1 \sigma_2) \int_0^1 \cY_s \, \dd s
     + \int_0^1 \cX_s \, \dd\cY_s \\
    -\frac{1}{2} \cX_1^2 + \alpha \int_0^1 \cX_s \, \dd s
     + \frac{\sigma_2^2}{2} \int_0^1 \cY_s \, \dd s + \frac{\sigma_3^2}{2}
   \end{bmatrix}
  \end{bmatrix}
 \end{align}
 as \ $T \to \infty$, \ where \ $(\cY_t, \cX_t)_{t\in\RR_+}$ \ is the unique strong solution of
 the SDE
 \begin{align}\label{c2dim_affine}
  \begin{cases}
   \dd \cY_t = a \, \dd t + \sigma_1 \sqrt{\cY_t} \, \dd W_t , \\
   \dd \cX_t
   = \alpha \, \dd t
     + \sigma_2 \sqrt{\cY_t} \, (\varrho \, \dd W_t + \sqrt{1 - \varrho^2} \, \dd B_t) ,
  \end{cases}
  \qquad t \in [0, \infty) ,
 \end{align}
 with initial value \ $(\cY_0, \cX_0) = (0, 0)$.
\end{Thm}

\noindent{\bf Proof.}
By \eqref{LSE_cont_true}, we have
 \[
   \begin{bmatrix}
    \ha_T - a \\
    \hb_T
   \end{bmatrix}
   = \begin{bmatrix}
      \ha_T - a \\
      \hb_T - b
     \end{bmatrix}
   = \begin{bmatrix}
      T & - \int_0^T Y_s \, \dd s \\
      - \int_0^T Y_s \, \dd s & \int_0^T Y_s^2 \, \dd s
     \end{bmatrix}^{-1}
     \begin{bmatrix}
      \sigma_1 \int_0^T Y_s^{\frac{1}{2}} \, \dd W_s \\
      - \sigma_1 \int_0^T Y_s^{\frac{3}{2}} \, \dd W_s
     \end{bmatrix} .
 \]
We can write
 \[
   \begin{bmatrix}
    T & - \int_0^T Y_s \, \dd s \\
    - \int_0^T Y_s \, \dd s & \int_0^T Y_s^2 \, \dd s
   \end{bmatrix}
   = \begin{bmatrix}
      T^{\frac{1}{2}} & 0 \\
      0 & T^{\frac{3}{2}}
     \end{bmatrix}
     \begin{bmatrix}
      1 & - \frac{1}{T^2} \int_0^T Y_s \, \dd s \\
      - \frac{1}{T^2} \int_0^T Y_s \, \dd s & \frac{1}{T^3} \int_0^T Y_s^2 \, \dd s
     \end{bmatrix}
     \begin{bmatrix}
      T^{\frac{1}{2}} & 0 \\
      0 & T^{\frac{3}{2}}
     \end{bmatrix}
 \]
 and
 \[
   \begin{bmatrix}
    \sigma_1 \int_0^T Y_s^{\frac{1}{2}} \, \dd W_s \\
    - \sigma_1 \int_0^T Y_s^{\frac{3}{2}} \, \dd W_s
   \end{bmatrix}
   = \begin{bmatrix}
      T & 0 \\
      0 & T^2
     \end{bmatrix}
     \begin{bmatrix}
      \frac{\sigma_1}{T} \int_0^T Y_s^{\frac{1}{2}} \, \dd W_s \\
      - \frac{\sigma_1}{T^2} \int_0^T Y_s^{\frac{3}{2}} \, \dd W_s
     \end{bmatrix} .
 \]
Consequently,
 \[
   \begin{bmatrix}
    \ha_T - a \\
    T \hb_T
   \end{bmatrix}
   = \begin{bmatrix}
      1 & 0 \\
      0 & T
     \end{bmatrix}
     \begin{bmatrix}
      \ha_T - a \\
      \hb_T
     \end{bmatrix}
   = \begin{bmatrix}
      1 & - \frac{1}{T^2} \int_0^T Y_s \, \dd s \\
      - \frac{1}{T^2} \int_0^T Y_s \, \dd s & \frac{1}{T^3} \int_0^T Y_s^2 \, \dd s
     \end{bmatrix}^{-1}
     \begin{bmatrix}
      \frac{\sigma_1}{T} \int_0^T Y_s^{\frac{1}{2}} \, \dd W_s \\
      - \frac{\sigma_1}{T^2} \int_0^T Y_s^{\frac{3}{2}} \, \dd W_s
     \end{bmatrix} .
 \]
In a similar way,
 \begin{align*}
  \begin{bmatrix}
   \halpha_T - \alpha \\
   T \hbeta_T \\
   T \hgamma_T
  \end{bmatrix}
  &= \begin{bmatrix}
      1 & - \frac{1}{T^2} \int_0^T Y_s \, \dd s & - \frac{1}{T^2} \int_0^T X_s \, \dd s \\
      - \frac{1}{T^2} \int_0^T Y_s \, \dd s & \frac{1}{T^3} \int_0^T Y_s^2 \, \dd s
       & \frac{1}{T^3} \int_0^T Y_s X_s \, \dd s \\
      - \frac{1}{T^2} \int_0^T X_s \, \dd s & \frac{1}{T^3} \int_0^T Y_s X_s \, \dd s
       & \frac{1}{T^3} \int_0^T X_s^2 \, \dd s
     \end{bmatrix}^{-1} \\
  &\quad
     \times
     \begin{bmatrix}
      \frac{\sigma_2}{T} \int_0^T Y_s^{\frac{1}{2}} \, \dd \tW_s + \frac{\sigma_3}{T} L_T \\
      - \frac{\sigma_2}{T^2} \int_0^T Y_s^{\frac{3}{2}} \, \dd \tW_s
       - \frac{\sigma_3}{T^2} \int_0^T Y_s \, \dd L_s \\
      - \frac{\sigma_2}{T^2} \int_0^T Y_s^{\frac{1}{2}} X_s \, \dd \tW_s
       - \frac{\sigma_3}{T^2} \int_0^T X_s \, \dd L_s
     \end{bmatrix} .
 \end{align*}
The aim of the following discussion is to prove
 \begin{equation}\label{intYX}
  \begin{aligned}
   &\biggl(\frac{1}{T} Y_T, \frac{1}{T} X_T, \frac{1}{T^2} \int_0^T X_s \, \dd Y_s, 
           \frac{1}{T^{k+\ell+1}} \int_0^T Y_s^k X_s^\ell \, \dd s
           : k, \ell \in \ZZ_+ , k + \ell \leq 2\biggr) \\
   &\distr \biggl(\cY_1, \cX_1, \int_0^1 \cX_s \, \dd\cY_s,
                  \int_0^1 \cY_s^k \cX_s^\ell \, \dd s
                  : k, \ell \in \ZZ_+ , k + \ell \leq 2\biggr)
  \end{aligned}
 \end{equation}
 as \ $T \to \infty$.
\ By part (ii) of Remark 2.7 in Barczy et al.\ \cite{BarDorLiPap}, we have
 \[
   \bigl(\tcY_t^{(T)}, \tcX_t^{(T)}\bigr)_{t\in\RR_+}
   := \Bigl(\frac{1}{T} \cY_{Tt}, \frac{1}{T} \cX_{Tt} \Bigr)_{t\in\RR_+}
   \distre (\cY_t, \cX_t)_{t\in\RR_+} \qquad
   \text{for all \ $T \in \RR_{++}$,}
 \]
 since, by Proposition \ref{Pro_affine}, \ $(\cY_t, \cX_t)_{t\in\RR_+}$ \ is an affine
 process with infinitesimal generator
 \[
   (\cA_{(\cY,\cX)} f)(y, x)
   = a f_1'(y, x) + \alpha f_2'(y, x)
     + \frac{1}{2} y
       \bigl[\sigma_1^2 f_{1,1}''(y, x)
             + 2 \varrho \sigma_1 \sigma_2 f_{1,2}''(y, x)
             + \sigma_2^2 f_{2, 2}''(y,x) \bigr] .
 \]
Hence, by Lemma \ref{Lem_intf_cont}, we obtain
 \begin{align*}
  &\biggl(\cY_1, \cX_1, \int_0^1 \cX_s \, \dd\cY_s, \int_0^1 \cY_s^k \cX_s^\ell \, \dd s
          : k, \ell \in \ZZ_+ , k + \ell \leq 2\biggr) \\
  &\distre
   \biggl(\tcY_1^{(T)}, \tcX_1^{(T)}, \int_0^1 \tcX_s^{(T)} \, \dd\tcY_s^{(T)},
          \int_0^1 \bigl(\tcY_s^{(T)}\bigr)^k \bigl(\tcX_s^{(T)}\bigr)^\ell \, \dd s
          : k, \ell \in \ZZ_+ , k + \ell \leq 2\biggr) \\
  &= \biggl(\frac{1}{T} \cY_T, \frac{1}{T} \cX_T, \frac{1}{T^2} \int_0^T \cX_s \, \dd\cY_s,
            \frac{1}{T^{k+\ell+1}} \int_0^T \cY_s^k \cX_s^\ell \, \dd s
            : k, \ell \in \ZZ_+ , k + \ell \leq 2\biggr)
 \end{align*}
 for all \ $T \in \RR_{++}$.
\ Then, by Slutsky's lemma, in order to prove \eqref{intYX}, it suffices to show the
 convergences
 \begin{gather}\label{intYXT}
  \frac{1}{T} (Y_T - \cY_T) \stoch 0 , \qquad
  \frac{1}{T} (X_T - \cX_T) \stoch 0 , \\
  \frac{1}{T^2} \left(\int_0^T X_s \, \dd Y_s - \int_0^T \cX_s \, \dd\cY_s\right)
  \stoch 0 , \qquad
  \frac{1}{T^{k+\ell+1}} \int_0^T (Y_s^k X_s^\ell - \cY_s^k \cX_s^\ell) \, \dd s \stoch 0
  \label{intYkXlT}
 \end{gather}
 as \ $T \to \infty$ \ for all \ $k, \ell \in \ZZ_+$ \ with \ $k + \ell \leq 2$.
\ By (3.21) in Barczy et al.\ \cite{BarDorLiPap}, we have
 \begin{align}\label{E|Y-tY|}
  \EE(|Y_s - \cY_s|) \leq \EE(Y_0) , \qquad s \in \RR_+ ,
 \end{align}
 hence
 \begin{align*}
  &\EE\left(\left|\frac{1}{T} (Y_T - \cY_T)\right|\right)
   \leq \frac{1}{T} \EE(Y_0)
   \to 0 , \\
  &\EE\left(\left|\frac{1}{T^2} \int_0^T (Y_s - \cY_s) \, \dd s\right|\right)
   \leq \frac{1}{T^2} \int_0^T \EE(|Y_s - \cY_s|) \, \dd s
   \leq \frac{1}{T} \EE(Y_0)
   \to 0 ,
 \end{align*}
 as \ $T \to \infty$, \ implying \ $\frac{1}{T} (Y_T - \cY_T) \stoch 0$ \ and
 \ $\frac{1}{T^2} \int_0^T (Y_s - \cY_s) \, \dd s \stoch 0$ \ as \ $T \to \infty$, \ i.e.,
 the first convergence in \eqref{intYXT} and the second convergence in \eqref{intYkXlT} for
 \ $(k, \ell) = (1, 0)$.

As in (3.23) in Barczy et al.\ \cite{BarDorLiPap}, we have
 \ $\EE(|X_s - \cX_s|) \leq \EE(|X_0|) + \sqrt{(\sigma_2^2 \EE(Y_0) + \sigma_3^2) s}$ \ for
 all \ $s \in \RR_+$, \ hence
 \begin{align}\label{E|X-tX|}
  \sup_{s\in[0,T]} \EE(|X_s - \cX_s|) = \OO(T^{\frac{1}{2}}) \qquad 
  \text{as \ $T \to \infty$,}
 \end{align}
 thus
 \begin{align*}
  &\EE\left(\left|\frac{1}{T} (X_T - \cX_T)\right|\right)
   = \frac{1}{T} \OO(T^{\frac{1}{2}})
   \to 0 , \\
  &\EE\left(\left|\frac{1}{T^2} \int_0^T (X_s - \cX_s) \, \dd s\right|\right)
   \leq \frac{1}{T^2} \int_0^T \EE(|X_s - \cX_s|) \, \dd s
   = \frac{1}{T^2} \int_0^T \OO(T^{\frac{1}{2}}) \, \dd s
   = \frac{1}{T^2} \OO(T^{\frac{3}{2}})     
   \to 0 ,
 \end{align*}
 as \ $T \to \infty$, \ implying \ $\frac{1}{T} (X_T - \cX_T) \stoch 0$ \ and
 \ $\frac{1}{T^2} \int_0^T (X_s - \cX_s) \, \dd s \stoch 0$ \ as \ $T \to \infty$, \ i.e.,
 the second convergence in \eqref{intYXT} and the second convergence in \eqref{intYkXlT} for
 \ $(k, \ell) = (0, 1)$.

As in (3.25) in Barczy et al.\ \cite{BarDorLiPap}, we have
 \ $\EE[(Y_s - \cY_s)^2] \leq 2 \EE(Y_0^2) + 2 s \sigma_1^2 \EE(Y_0)$ \ for all
 \ $s \in \RR_+$, \ hence
 \begin{equation}\label{(Ys-cYs)2}
  \sup_{s\in[0,T]} \EE[(Y_s - \cY_s)^2] = \OO(T) \qquad \text{as \ $T \to \infty$.}
 \end{equation}
By Proposition \ref{Pro_momentsnp},
 \ $\EE(Y_s^2) = \EE(Y_0^2) + (2 a + \sigma_1^2) \bigl(\EE(Y_0) s + a \frac{s^2}{2}\bigr)$
 \ for all \ $s \in \RR_+$, \ hence
 \begin{equation}\label{Ys2}
  \sup_{s\in[0,T]} \EE(Y_s^2) = \OO(T^2) \qquad \text{as \ $T \to \infty$,}
 \end{equation}
 and \ $\sup_{s\in[0,T]} \EE(\cY_s^2) = \OO(T^2)$ \ as \ $T \to \infty$. 
\ We have
 \begin{align*}
  \EE(|Y_s^2 - \cY_s^2|)
  = \EE(|(Y_s - \cY_s) (Y_s + \cY_s)|)
  &\leq \sqrt{\EE[(Y_s - \cY_s)^2] \EE[(Y_s + \cY_s)^2]} \\
  &\leq \sqrt{2 \EE[(Y_s - \cY_s)^2] (\EE(Y_s^2) + \EE(\cY_s^2))} ,
 \end{align*}
 yielding
 \[
   \sup_{s\in[0,T]} \EE(|Y_s^2 - \cY_s^2|)
   = \sqrt{2 \OO(T) (\OO(T^2) + \OO(T^2))}
   = \OO(T^{\frac{3}{2}}) \qquad \text{as \ $T \to \infty$,}
 \]
 thus
 \[
   \EE\left(\left|\frac{1}{T^3} \int_0^T (Y_s^2 - \cY_s^2) \, \dd s\right|\right)
   \leq \frac{1}{T^3} \int_0^T \EE(|Y_s^2 - \cY_s^2|) \, \dd s
   = \frac{1}{T^3} \int_0^T \OO(T^{\frac{3}{2}}) \, \dd s
   = \frac{1}{T^3} \OO(T^{\frac{5}{2}})
   \to 0 ,
 \]
 as \ $T \to \infty$, \ implying
 \ $\frac{1}{T^3} \int_0^T (Y_s^2 - \cY_s^2) \, \dd s \stoch 0$ \ as \ $T \to \infty$.,
 \ i.e., the second convergence in \eqref{intYkXlT} for \ $(k, \ell) = (2, 0)$.
 
In a similar way,
 \ $\EE[(X_s - \cX_s)^2] \leq 2 \EE(X_0^2) + 2 s (\sigma_2^2 \EE(Y_0) + \sigma_3^2)$ \ for all
 \ $s \in \RR_+$, \ hence
 \begin{equation}\label{(Xs-cXs)2}
  \sup_{s\in[0,T]} \EE[(X_s - \cX_s)^2] = \OO(T) \qquad \text{as \ $T \to \infty$.}
 \end{equation}
By Proposition \ref{Pro_momentsnp},
 \ $\EE(X_s^2) = \EE(X_0^2) + \alpha \bigl(s \EE(X_0) + \alpha \frac{s^2}{2}\bigr)
                 + \sigma_2^2 \bigl(s \EE(Y_0) + a \frac{s^2}{2}\bigr) + \sigma_3^2 s$, \ thus
 \ $\sup_{s\in[0,T]} \EE(X_s^2) = \OO(T^2)$ \ and
 \ $\sup_{s\in[0,T]} \EE(\cX_s^2) = \OO(T^2)$ \ as \ $T \to \infty$. 
\ We have
 \[
   \EE(|X_s^2 - \cX_s^2|)
   \leq \sqrt{2 \EE[(X_s - \cX_s)^2] (\EE(X_s^2) + \EE(\cX_s^2))} ,
 \]
 yielding
 \[
   \sup_{s\in[0,T]} \EE(|X_s^2 - \cX_s^2|)
   = \sqrt{2 \OO(T) (\OO(T^2) + \OO(T^2))}
   = \OO(T^{\frac{3}{2}}) \qquad \text{as \ $T \to \infty$,}
 \]
 thus 
 \[
   \EE\left(\left|\frac{1}{T^3} \int_0^T (X_s^2 - \cX_s^2) \, \dd s\right|\right)
   \leq \frac{1}{T^3} \int_0^T \EE(|X_s^2 - \cX_s^2|) \, \dd s
   = \frac{1}{T^3} \int_0^T \OO(T^{\frac{3}{2}}) \, \dd s
   = \frac{1}{T^3} \OO(T^{\frac{5}{2}})
   \to 0 ,
 \]
 as \ $T \to \infty$, \ implying 
 \ $\frac{1}{T^3} \int_0^T (X_s^2 - \cX_s^2) \, \dd s \stoch 0$ \ as \ $T \to \infty$, \ i.e.,
 the second convergence in \eqref{intYkXlT} for \ $(k, \ell) = (0, 2)$.

Further,
 \begin{align*}
  \EE(|Y_s X_s - \cY_s \cX_s|)
  &\leq \EE(|Y_s - \cY_s| |X_s|) + \EE(\cY_s |X_s - \cX_s|) \\
  &\leq \sqrt{\EE[(Y_s - \cY_s)^2] \EE(X_s^2)} + \sqrt{\EE(\cY_s^2) \EE[(X_s - \cX_s)^2]}
 \end{align*}
 yields
 \[
   \sup_{s\in[0,T]} \EE(|Y_s X_s - \cY_s \cX_s|)
   = \sqrt{\OO(T) \OO(T^2)} + \sqrt{\OO(T^2) \OO(T))}
   = \OO(T^{\frac{3}{2}}) \qquad \text{as \ $T \to \infty$,}
 \]
 thus
 \[
   \EE\left(\left|\frac{1}{T^3} \int_0^T (Y_s X_s - \cY_s \cX_s) \, \dd s\right|\right)
   \leq \frac{1}{T^3} \int_0^T \EE(|Y_s X_s - \cY_s \cX_s|) \, \dd s
   = \frac{1}{T^3} \int_0^T \OO(T^{\frac{3}{2}}) \, \dd s
   = \frac{1}{T^3} \OO(T^{\frac{5}{2}})
   \to 0 ,
 \]
 as \ $T \to \infty$, \ implying
 \ $\frac{1}{T^3} \int_0^T (Y_s X_s - \cY_s \cX_s) \, \dd s \stoch 0$ \ as \ $T \to \infty$, 
 \ i.e., the second convergence in \eqref{intYkXlT} for \ $(k, \ell) = (1, 1)$.

Using the Cauchy--Schwarz inequality, we obtain
 \begin{align*}
  \EE\left(\left|\int_0^T X_s \, \dd Y_s - \int_0^T \cX_s \, \dd \cY_s\right|\right)
  &\leq \EE\left(\left|\int_0^T (X_s - \cX_s) \, \dd Y_s\right|\right)
        + \EE\left(\left|\int_0^T \cX_s \, \dd (Y_s - \cY_s)\right|\right) \\
  &\leq \sqrt{E_1(T)} + \sqrt{E_2(T)} 
 \end{align*}
 with
 \[
   E_1(T) := \EE\left(\left|\int_0^T (X_s - \cX_s) \, \dd Y_s\right|^2\right) , \qquad
   E_2(T) := \EE\left(\left|\int_0^T \cX_s \, \dd (Y_s - \cY_s)\right|^2\right) .
 \]
Using \ $\dd Y_s = a \, \dd s + \sigma_1 \sqrt{Y_s} \, \dd W_s$, \ we have
 \[
   E_1(T)
   = \EE\left(\left|a \int_0^T (X_s - \cX_s) \, \dd s
                    + \sigma_1 \int_0^T (X_s - \cX_s) \sqrt{Y_s} \, \dd W_s\right|^2\right)
   \leq 2 a^2 E_{1,1}(T) + 2 \sigma_1^2 E_{1,2}(T)
 \]
 with
 \[
   E_{1,1}(T) := \EE\left(\left|\int_0^T (X_s - \cX_s) \, \dd s\right|^2\right) , \qquad
   E_{1,2}(T) := \EE\left(\left|\int_0^T (X_s - \cX_s) \sqrt{Y_s} \, \dd W_s\right|^2\right) .
 \]
Applying \eqref{(Xs-cXs)2}, we obtain
 \begin{align*}
  E_{1,1}(T) 
  &= \EE\left(\int_0^T \int_0^T (X_s - \cX_s) (X_u - \cX_u) \, \dd s \, \dd u\right)
   = \int_0^T \int_0^T \EE[(X_s - \cX_s) (X_u - \cX_u)] \, \dd s \, \dd u \\
  &\leq \int_0^T \int_0^T
        \sqrt{\EE[(X_s - \cX_s)^2] \EE[(X_u - \cX_u)^2]} \, \dd s \, \dd u
   = \int_0^T \int_0^T \sqrt{\OO(T) \OO(T)} \, \dd s \, \dd u
   = \OO(T^3) .
 \end{align*}
Again by the Cauchy--Schwarz inequality, we obtain
 \[
   E_{1,2}(T)
   = \EE\left(\int_0^T (X_s - \cX_s)^2 Y_s \, \dd s\right)
   = \int_0^T \EE[(X_s - \cX_s)^2 Y_s] \, \dd s
   \leq \int_0^T \sqrt{\EE[(X_s - \cX_s)^4] \EE(Y_s^2)} \, \dd s .
 \]
Using \ $X_t = X_0 + \sigma_2 \int_0^t \sqrt{Y_s} \, \dd \tW_s + \sigma_3 L_t$ \ and
 \ $\cX_t = \sigma_2 \int_0^t \sqrt{\cY_s} \, \dd\tW_s$, \ we get
 \ $X_t - \cX_t
    = X_0 + \sigma_2 \int_0^t (\sqrt{Y_s} - \sqrt{\cY_s}) \, \dd\tW_s + \sigma_3 L_t$,
 \ and, applying Minkowski inequality and a martingale moment inequality in Karatzas and
 Shreve \cite[3.3.25]{KarShr}, we obtain
 \begin{align*}
  (\EE[(X_t - \cX_t)^4])^{\frac{1}{4}}
  &\leq [\EE(X_0^4)]^{\frac{1}{4}}
        + \sigma_2
          \left(\EE\left[\left(\int_0^t (\sqrt{Y_s} - \sqrt{\cY_s}) \, \dd\tW_s\right)^4
                   \right]\right)^{\frac{1}{4}}
        + \sigma_3 [\EE(L_t^4)]^{\frac{1}{4}} \\
  &\leq [\EE(X_0^4)]^{\frac{1}{4}}
        + \sigma_2
          \left((2 \cdot 3)^2 t
                \EE\left(\int_0^t (\sqrt{Y_s} - \sqrt{\cY_s})^4 \, \dd s\right)
          \right)^{\frac{1}{4}}
        + \sigma_3 \sqrt[4]{3} \sqrt{t} \\
  &\leq [\EE(X_0^4)]^{\frac{1}{4}}
        + \sigma_2 \left(36 t \int_0^t \EE[(Y_s - \cY_s)^2] \, \dd s\right)^{\frac{1}{4}}
        + \sigma_3 \sqrt[4]{3} \sqrt{t} .
 \end{align*}
Applying \eqref{(Ys-cYs)2}, we get
 \begin{equation}\label{sup_(Xs-cXs)4}
  \sup_{t\in[0,T]} \EE[(X_t - \cX_t)^4] = \OO(T^3) \qquad \text{as \ $T \to \infty$,}
 \end{equation}
 which, by \eqref{Ys2}, implies 
 \ $E_{1,2}(T) = \int_0^T \sqrt{\OO(T^3) \OO(T^2)} \, \dd s = \OO(T^{\frac{7}{2}})$ \ as
 \ $T \to \infty$.
\ Using \ $E_{1,1}(T) = \OO(T^3)$ \ as \ $T \to \infty$, \ we conclude
 \ $E_1(T) = \OO(T^3) + \OO(T^{\frac{7}{2}})= \OO(T^{\frac{7}{2}})$ \ as \ $T \to \infty$.

Using \ $\dd Y_s = a \, \dd s + \sigma_1 \sqrt{Y_s} \, \dd W_s$ \ and
 \ $\dd \cY_s = a \, \dd s + \sigma_1 \sqrt{\cY_s} \, \dd W_s$, \ we obtain
 \ $\dd(Y_t - \cY_t) = \sigma_1 (\sqrt{Y_t} - \sqrt{\cY_t}) \, \dd W_t$, \ thus
 \begin{align*}
  E_2(T)
  = \sigma_1^2 \EE\left(\int_0^T \cX_s^2 (\sqrt{Y_s} - \sqrt{\cY_s})^2 \, \dd s\right)
  &\leq \sigma_1^2 \int_0^T \EE[\cX_s^2 |Y_s - \cY_s|] \, \dd s \\
  &\leq \sigma_1^2 \int_0^T \sqrt{\EE(\cX_s^4) \EE[(Y_s - \cY_s)^2]} \, \dd s .
 \end{align*}
Using \ $\cX_t = \alpha t + \sigma_2 \int_0^t \sqrt{\cY_s} \, \dd \tW_s$, \ we obtain
 \begin{align*}
  [\EE(\cX_t^4)]^{\frac{1}{4}}
  &\leq |\alpha| t
        + \sigma_2
          \left(\EE\left[\left(\int_0^t \sqrt{\cY_s} \, \dd\tW_s\right)^4
                   \right]\right)^{\frac{1}{4}}
   \leq |\alpha| t
        + \sigma_2
          \left((2 \cdot 3)^2 t \EE\left(\int_0^t \cY_s^2 \, \dd s\right)
          \right)^{\frac{1}{4}} \\
  &= |\alpha| t
     + \sigma_2
       \left(36 t \int_0^t a \left(a + \frac{\sigma_1^2}{2}\right) s^2 \, \dd s
       \right)^{\frac{1}{4}}
   = \left(|\alpha| + \sigma_2 \sqrt[4]{6 a (2 a + \sigma_1^2)}\right) t ,
 \end{align*}
 hence we conclude 
 \begin{equation}\label{sup_cXs4}
  \sup_{s\in[0,T]} \EE(\cX_s^4) = \OO(T^4) \qquad \text{as \ $T \to \infty$.}
 \end{equation}
Using \eqref{(Ys-cYs)2}, we obtain
 \ $E_2(T) = \int_0^T \sqrt{\OO(T^4) \OO(T)} \, \dd s = \OO(T^{\frac{7}{2}})$ \ as
 \ $T \to \infty$.
\ Hence
 \[
   \EE\left(\left|\frac{1}{T^2} 
                  \left(\int_0^T X_s \, \dd Y_s
                        - \int_0^T \cX_s \, \dd \cY_s\right)\right|\right)
   \leq \frac{1}{T^2} \bigl(\sqrt{E_1(T)} + \sqrt{E_2(T)}\bigr)
   = \frac{1}{T^2} \OO(T^{\frac{7}{4}})
   \to 0
 \]
 as \ $T \to \infty$, \ implying
 \ $\frac{1}{T^2} \left(\int_0^T X_s \, \dd Y_s - \int_0^T X_s \, \dd Y_s\right) \stoch 0$
 \ as \ $T \to \infty$, \ i.e., the first convergence in \eqref{intYkXlT}.
Thus we conclude convergence \eqref{intYX}.
 
Applying the first equation of \eqref{2dim_affine} and using \ $b = 0$, \ we obtain
 \[
   \frac{\sigma_1}{T} \int_0^T Y_s^{\frac{1}{2}} \, \dd W_s
   = \frac{Y_T - Y_0 - a T}{T} \distr \cY_1 - a ,
   \qquad \text{as \ $T \to \infty$.}
 \]
By It\^o's formula and using \ $b = 0$,
 \[
   \dd(Y_t^2)
   = 2 Y_t \, \dd Y_t + \sigma_1^2 Y_t \, \dd t
   = 2 Y_t (a \, \dd t + \sigma_1 Y_t^{\frac{1}{2}} \, \dd W_t) + \sigma_1^2 Y_t \, \dd t
   = (2 a + \sigma_1^2) Y_t \, \dd t + 2 \sigma_1 Y_t^{\frac{3}{2}} \, \dd W_t ,
 \]
 hence
 \[
   Y_T^2 = Y_0^2 + (2 a + \sigma_1^2) \int_0^T Y_s \, \dd s
           + 2 \sigma_1 \int_0^T Y_s^{\frac{3}{2}} \, \dd W_s .
 \]
Consequently,
 \[
   - \frac{\sigma_1}{T^2} \int_0^T Y_s^{\frac{3}{2}} \, \dd W_s
   = - \frac{Y_T^2 - Y_0^2 - (2 a + \sigma_1^2) \int_0^T Y_s \, \dd s}{2T^2}
   \distr - \frac{\cY_1^2 - (2 a + \sigma_1^2) \int_0^1 \cY_s \, \dd s}{2}
 \]
 as \ $T \to \infty$.
\ In a similar way, applying the second equation of \eqref{2dim_affine} and using
 \ $\beta = 0$ \ and \ $\gamma = 0$, \ we obtain
 \[
   \frac{\sigma_2}{T} \int_0^T Y_s^{\frac{1}{2}} \, \dd \tW_s + \frac{\sigma_3}{T} L_T
   = \frac{X_T - X_0 - \alpha T}{T}
   \distr \cX_1 - \alpha , \qquad \text{as \ $T \to \infty$.}
 \]
By It\^o's formula and using \ $\beta = 0$ \ and \ $\gamma = 0$,
 \begin{align*}
  \dd(Y_t X_t)
  &= Y_t \, \dd X_t + X_t \, \dd Y_t + \varrho \sigma_1 \sigma_2 Y_t \, \dd t
   = Y_t (\alpha \, \dd t + \sigma_2 Y_t^{\frac{1}{2}} \, \dd \tW_t + \sigma_3 \, \dd L_t)
     + X_t \, \dd Y_t + \varrho \sigma_1 \sigma_2 Y_t \, \dd t \\
  &= (\alpha + \varrho \sigma_1 \sigma_2) Y_t \, \dd t
     + \sigma_2 Y_t^{\frac{3}{2}} \, \dd \tW_t
     + X_t \, \dd Y_t + \sigma_3 Y_t \, \dd L_t ,
 \end{align*}
 hence
 \[
   Y_T X_T = Y_0 X_0 + (\alpha + \varrho \sigma_1 \sigma_2) \int_0^T Y_s \, \dd s 
             + \sigma_2 \int_0^T Y_s^{\frac{3}{2}} \, \dd \tW_s + \int_0^T X_s \, \dd Y_s
             + \sigma_3 \int_0^T Y_s \, \dd L_s .
 \]
Consequently,
 \begin{align*}
  - \frac{\sigma_2}{T^2} \int_0^T Y_s^{\frac{3}{2}} \, \dd \tW_s
  - \frac{\sigma_3}{T^2} \int_0^T Y_s \, \dd L_s
  &= - \frac{Y_T X_T - Y_0 X_0 - (\alpha + \varrho \sigma_1 \sigma_2) \int_0^T Y_s \, \dd s
             - \int_0^T X_s \, \dd Y_s}
            {T^2} \\
  &\distr - \cY_1 \cX_1 + (\alpha + \varrho \sigma_1 \sigma_2) \int_0^1 \cY_s \, \dd s
          + \int_0^1 \cX_s \, \dd \cY_s
 \end{align*}
 as \ $T \to \infty$.
\ Again by It\^o's formula and using \ $\beta = 0$ \ and \ $\gamma = 0$,
 \[
   \dd(X_t^2)
   = 2 X_t \, \dd X_t + (\sigma_2^2 Y_t + \sigma_3^2) \, \dd t
   = 2 X_t (\alpha \, \dd t + \sigma_2 Y_t^{\frac{1}{2}} \, \dd \tW_t + \sigma_3 \, \dd L_t)
     + (\sigma_2^2 Y_t + \sigma_3^2) \, \dd t ,
 \]
 hence
 \[
   X_T^2 = X_0^2 + \int_0^T (2 \alpha X_s + \sigma_2^2 Y_s + \sigma_3^2) \, \dd s 
           + 2 \sigma_2 \int_0^T Y_s^{\frac{1}{2}} X_s \, \dd \tW_s
           + 2 \sigma_3 \int_0^T X_s \, \dd L_s .
 \]
Consequently,
 \begin{align*}
  - \frac{\sigma_2}{T^2} \int_0^T Y_s^{\frac{1}{2}} X_s \, \dd \tW_s
  - \frac{\sigma_3}{T^2} \int_0^T X_s \, \dd L_s
  &= - \frac{X_T^2 - X_0^2 - \int_0^T (2 \alpha X_s + \sigma_2^2 Y_s + \sigma_3^2) \, \dd s}
            {2T^2} \\
  &\distr
     - \frac{\cX_1^2 - \int_0^1 (2 \alpha \cX_s + \sigma_2^2 \cY_s + \sigma_3^2) \, \dd s}{2}
 \end{align*}
 as \ $T \to \infty$, \ and we conclude \eqref{LSE_crit}.
\proofend

\section{Asymptotic behavior of CLSE: supercritical case}
\label{section_ALSE_supercritical}

First we present an auxiliary lemma about the asymptotic behavior of \ $\EE(X_t^2)$ \ as
 \ $t \to \infty$.

\begin{Lem}\label{momentX2}
Let us consider the two-factor affine diffusion model \eqref{2dim_affine} with
 \ $a \in \RR_+$, \ $b \in \RR_{--}$, \ $\alpha, \beta \in \RR$, \ $\gamma \in (-\infty, b)$,
 \ $\sigma_1 \in \RR_{++}$, \ $\sigma_2, \sigma_3 \in \RR_+$ \ and \ $\varrho \in [-1, 1]$
 \ with a random initial value \ $(\eta_0, \zeta_0)$ \ independent of
 \ $(W_t, B_t, L_t)_{t\in\RR_+}$ \ satisfying \ $\PP(\eta_0 \in \RR_+) = 1$.
\ Then \ $\sup_{t\in\RR_+} \ee^{2\gamma t} \EE(X_t^2) < \infty$.
\end{Lem}

\noindent{\bf Proof.}
By Proposition \ref{Pro_momentsnp},
 \[
   \sup_{t\in\RR_+} \ee^{bt} \EE(Y_t)
   = \sup_{t\in\RR_+} \bigl(\EE(Y_0) + a \int_0^t \ee^{bu} \, \dd u\bigr)
   = \EE(Y_0) + a \int_0^\infty \ee^{bu} \, \dd u < \infty ,
 \]
 since \ $b < 0$. 
\ Moreover,
 \begin{align*}
  \sup_{t\in\RR_+} \ee^{\gamma t} |\EE(X_t)|
  &= \sup_{t\in\RR_+}
      \biggl|\EE(X_0) + \alpha \int_0^t \ee^{\gamma u} \, \dd u 
             - \beta \int_0^t \ee^{\gamma u} \EE(Y_u) \, \dd u\biggr| \\
  &\leq |\EE(X_0)| + |\alpha| \int_0^\infty \ee^{\gamma u} \, \dd u
        + |\beta| \biggl(\sup_{u\in\RR_+} \ee^{bu} \EE(Y_u)\biggr) 
          \int_0^\infty \ee^{(\gamma-b)u} \, \dd u
  < \infty ,
 \end{align*}
 using \ $\gamma < 0$ \ and \ $\gamma - b < 0$.
\ Again by Proposition \ref{Pro_momentsnp},
 \begin{align*}
  \sup_{t\in\RR_+} \ee^{2bt} \EE(Y_t^2)
  &= \sup_{t\in\RR_+}
      \biggl(\EE(Y_0^2) + (2a + \sigma_1^2) \int_0^t \ee^{2bu} \EE(Y_u) \, \dd u\biggr) \\
  &\leq \EE(Y_0^2)
        + (2a + \sigma_1^2) \biggl(\sup_{u\in\RR_+} \ee^{bu} \EE(Y_u)\biggr)
          \int_0^\infty \ee^{bu} \, \dd u
   < \infty ,
 \end{align*}
 using \ $b < 0$.
\ Hence
 \begin{align*}
  &\sup_{t\in\RR_+} \ee^{(b+\gamma)t} |\EE(Y_t X_t)|
   = \sup_{t\in\RR_+}
      \biggl|\EE(Y_0 X_0)
             + a \int_0^t \ee^{(b+\gamma)u} \EE(X_u) \, \dd u \\
  &\phantom{\sup_{t\in\RR_+} \ee^{(b+\gamma)t} |\EE(Y_t X_t)| = \sup_{t\in\RR_+} \biggl(}
     + (\alpha + \varrho \sigma_1 \sigma_2)
       \int_0^t \ee^{(b+\gamma)u} \EE(Y_u) \, \dd u
     - \beta \int_0^t \ee^{(b+\gamma)u} \EE(Y_u^2) \, \dd u\biggr| \\
  &\leq |\EE(Y_0 X_0)|
        + a \biggl(\sup_{u\in\RR_+} \ee^{\gamma u} |\EE(X_u)|\biggr)
          \int_0^\infty \ee^{bu} \, \dd u
        + (|\alpha| + |\varrho| \sigma_1 \sigma_2)
          \biggl(\sup_{u\in\RR_+} \ee^{bu} \EE(Y_u)\biggr)
          \int_0^\infty \ee^{\gamma u} \, \dd u \\
  &\quad
        + |\beta| \biggl(\sup_{u\in\RR_+} \ee^{2bu} \EE(Y_u^2)\biggr)
          \int_0^\infty \ee^{(\gamma-b)u} \, \dd u
   < \infty ,
 \end{align*}
 using \ $b < 0$, \ $\gamma < 0$ \ and \ $\gamma - b < 0$.
\ Consequently,
 \begin{align*}
  &\sup_{t\in\RR_+} \ee^{2\gamma t} \EE(X_t^2)
   = \sup_{t\in\RR_+}
      \biggl(\EE(X_0^2)
             + \alpha \int_0^t \ee^{2\gamma u} X_u \, \dd u
             - 2 \beta \int_0^t \ee^{2\gamma u} Y_u X_u \, \dd u \\
  &\phantom{\sup_{t\in\RR_+} \ee^{2\gamma t} \EE(X_t^2) = \sup_{t\in\RR_+} \biggl(}
             + \sigma_2^2 \int_0^t \ee^{2\gamma u} Y_u \, \dd u
             + \sigma_3^2 \int_0^t \ee^{2\gamma u} \, \dd u\biggl) \\
  &\leq \EE(X_0^2)
        + |\alpha| \biggl(\sup_{u\in\RR_+} \ee^{\gamma u} |\EE(X_u)|\biggr)
          \int_0^\infty \ee^{\gamma u} \, \dd u
        + 2 |\beta| \biggl(\sup_{u\in\RR_+} \ee^{(b+\gamma)u} |\EE(Y_u X_u)|\biggr)
          \int_0^\infty \ee^{(\gamma-b)u} \, \dd u \\
  &\quad
        + \sigma_2^2 \biggl(\sup_{u\in\RR_+} \ee^{bu} \EE(Y_u)\biggr)
          \int_0^\infty \ee^{(2\gamma-b) u} \, \dd u
        + \sigma_3^2 \int_0^\infty \ee^{2\gamma u} \, \dd u
   < \infty
 \end{align*}
 using \ $\gamma < 0$, \ $\gamma - b < 0$ \ and \ $2 \gamma - b < 0$.
\proofend

Next we present an auxiliary lemma about the asymptotic behavior of \ $X_t$ \ as
 \ $t \to \infty$.

\begin{Lem}\label{martconvYX}
Let us consider the two-factor affine diffusion model \eqref{2dim_affine} with
 \ $a \in \RR_+$, \ $b \in \RR_{--}$, \ $\alpha, \beta \in \RR$, \ $\gamma \in (-\infty, b)$,
 \ $\sigma_1 \in \RR_{++}$, \ $\sigma_2, \sigma_3 \in \RR_+$ \ and \ $\varrho \in [-1, 1]$
 \ with a random initial value \ $(\eta_0, \zeta_0)$ \ independent of
 \ $(W_t, B_t, L_t)_{t\in\RR_+}$ \ satisfying \ $\PP(\eta_0 \in \RR_+) = 1$.
\ Suppose that \ $\alpha \beta \in \RR_-$.
\ Then there exists a random variable \ $V_X$ \ such that
 \begin{align}\label{lim_X}
  \ee^{\gamma t} X_t \as V_X \qquad \text{as \ $t \to \infty$}
 \end{align}
 and, for each \ $k, \ell \in \ZZ_+$ \ with \ $k + \ell > 0$,
 \begin{align}\label{lim_intYX}
  \ee^{(kb+\ell\gamma)t} \int_0^t Y_u^k X_u^\ell \, \dd u
  \as -\frac{V_Y^k V_X^\ell}{kb+\ell\gamma} \qquad \text{as \ $t \to \infty$,}
 \end{align}
 where \ $V_Y$ \ is given in \eqref{lim_Y}.
If, in addition, \ $\sigma_3 \in \RR_{++}$ \ or
 \ $\bigl(a - \frac{\sigma_1^2}{2}\bigr) (1 - \varrho^2) \sigma_2^2 \in \RR_{++}$,
 \ then the distribution of the random variable \ $V_X$ \ is absolutely continuous.
Particularly, \ $\PP(V_X \ne 0) = 1$.
\end{Lem}

\noindent{\bf Proof.}
By \eqref{SolutionX},
 \begin{align*}
  \EE( X_t \mid \cF_s )
  = \EE( X_t \mid Y_s, X_s )
  = \ee^{-\gamma(t-s)} X_s + \int_s^t \ee^{-\gamma(t-u)} (\alpha - \beta Y_u) \, \dd u
 \end{align*}
 for all \ $s, t \in \RR_+$ \ with \ $0 \leq s \leq t$.
\ If \ $\alpha \in \RR_+$ \ and \ $\beta \in \RR_-$, \ then
 \[
   \EE( \ee^{\gamma t} X_t \mid \cF^{Y,X}_s )
   = \ee^{\gamma s} X_s + \int_s^t \ee^{\gamma u} (\alpha - \beta Y_u) \, \dd u
   \geq \ee^{\gamma s} X_s
 \]
 for all \ $s, t \in \RR_+$ \ with \ $0 \leq s \leq t$, \ consequently, the process
 \ $(\ee^{\gamma t} X_t)_{t\in\RR_+}$ \ is a submartingale with respect to the filtration
 \ $(\cF^{Y,X}_t)_{t\in\RR_+}$.
\ If \ $\alpha \in \RR_-$ \ and \ $\beta \in \RR_+$, \ then
 \[
   \EE( \ee^{\gamma t} X_t \mid \cF^{Y,X}_s )
   = \ee^{\gamma s} X_s + \int_s^t \ee^{\gamma u} (\alpha - \beta Y_u) \, \dd u
   \leq \ee^{\gamma s} X_s
 \]
 for all \ $s, t \in \RR_+$ \ with \ $0 \leq s \leq t$, \ consequently, the process
 \ $(\ee^{\gamma t} X_t)_{t\in\RR_+}$ \ is a supermartingale with respect to the filtration
 \ $(\cF^{Y,X}_t)_{t\in\RR_+}$, \ hence the process \ $(-\ee^{\gamma t} X_t)_{t\in\RR_+}$
 \ is a submartingale with respect to the filtration \ $(\cF^{Y,X}_t)_{t\in\RR_+}$.
\ In both cases, \ $\sup_{t\in\RR_+} \EE(|\ee^{\gamma t} X_t|^2) < \infty$, \ see Lemma
 \ref{momentX2}.
Hence, by the submartingale convergence theorem, there exists a random variable \ $V_X$
 \ such that \eqref{lim_X} holds.
 
If \ $\omega \in \Omega$ \ such that \ $\RR_+ \ni t \mapsto (Y_t(\omega), X_t(\omega))$ \ is
 continuous and
 \ $(\ee^{bt} Y_t(\omega), \ee^{\gamma t} X_t(\omega)) \to (V_Y(\omega), V_X(\omega))$ \ as
 \ $t \to \infty$, \ then, by the integral Kronecker Lemma \ref{int_Kronecker} with
 \ $f(t) = \ee^{(kb+\ell\gamma)t} Y_t(\omega)^k X_t(\omega)^\ell$ \ and
 \ $a(t) = \ee^{-(kb+\ell\gamma)t}$, \ $t \in \RR_+$, \ we have
 \[
   \frac{1}{\int_0^t \ee^{-(kb+\ell\gamma)u} \, \dd u}
   \int_0^t
    \ee^{-(kb+\ell\gamma)u} (\ee^{(kb+\ell\gamma)u} Y_u(\omega)^k X_u(\omega)^\ell) \, \dd u
   \to V_Y(\omega)^k V_X(\omega)^\ell \qquad \text{as \ $t \to \infty$.}
 \]
Here
 \ $\int_0^t \ee^{-(kb+\ell\gamma)u} \, \dd u
    = - \frac{\ee^{-(kb+\ell\gamma)t} - 1}{kb+\ell\gamma}$,
 \ $t \in \RR_+$, \ thus we conclude \eqref{lim_intYX}.
 
Now suppose that \ $\sigma_3 \in \RR_{++}$ \ or
 \ $\bigl(a - \frac{\sigma_1^2}{2}\bigr) (1 - \varrho^2) \sigma_2^2 \in \RR_{++}$.
\ We are going to show that the random variable \ $V_X$ \ is absolutely continuous.
Put \ $Z_t := X_t - r Y_t$, \ $t \in \RR_+$ \ with
 \ $r := \frac{\sigma_2 \varrho}{\sigma_1}$.
\ Then the process \ $(Y_t, Z_t)_{t\in\RR_+}$ \ is an affine process satisfying
 \[
   \begin{cases}
    \dd Y_t = (a - b Y_t) \, \dd t + \sigma_1 \sqrt{Y_t} \, \dd W_t , \\
    \dd Z_t = (A - B Y_t - \gamma Z_t) \, \dd t + \Sigma_2 \sqrt{Y_t} \, \dd B_t
              + \sigma_3 \, \dd L_t .
   \end{cases} \qquad t \in \RR_+ ,
 \]
 where \ $A := \alpha - r a$, \ $B := \beta - r (b - \gamma)$ \ and
 \ $\Sigma_2 := \sigma_2 \sqrt{1 - \varrho^2}$, \ see Bolyog and Pap
 \cite[Proposition 2.5]{BolPap}.
We have
 \begin{align*}
  \ee^{\gamma t} X_t
  = r \ee^{\gamma t} Y_t + \ee^{\gamma t} Z_t
  = r \ee^{\gamma t} Y_t + Z_0 + \int_0^t \ee^{\gamma u} (A - B Y_u) \, \dd u
    + \Sigma_2 \int_0^t \ee^{\gamma u} \sqrt{Y_u} \, \dd B_u
    + \sigma_3 \int_0^t \ee^{\gamma u} \, \dd L_u ,
 \end{align*}
 where we used \eqref{SolutionX} with \ $s = 0$ \ multiplied both sides by \ $\ee^{\gamma t}$.
\ Thus the conditional distribution of \ $\ee^{\gamma t} X_t$ \ given \ $(Y_u)_{u\in[0,t]}$
 \ and \ $X_0$ \ is a normal distribution with mean
 \ $r \ee^{\gamma t} Y_t + Z_0 + \int_0^t \ee^{\gamma u} (A - B Y_u) \, \dd u$ \ and with
 variance
 \ $\Sigma_2^2 \int_0^t \ee^{2\gamma u} Y_u \, \dd u
    + \sigma_3^2 \int_0^t \ee^{2\gamma u} \, \dd u$.
\ Hence
 \begin{align*}
  &\EE\bigl(\ee^{\ii \lambda \ee^{\gamma t} X_t} \bmid (Y_u)_{u\in[0,t]}, X_0\bigr) \\
  &= \exp\biggl\{\ii \lambda
                 \biggl(r \ee^{\gamma t} Y_t + Z_0
                        + \int_0^t \ee^{\gamma u} (A - B Y_u) \, \dd u\biggr)
                 - \frac{\lambda^2}{2}
                   \biggl(\Sigma_2^2 \int_0^t \ee^{2\gamma u} Y_u \, \dd u
                          + \sigma_3^2 \int_0^t \ee^{2\gamma u} \, \dd u\biggr)\biggr\} .
 \end{align*}
Consequently,
 \begin{align*}
  &\bigl|\EE\bigl(\ee^{\ii \lambda \ee^{\gamma t} X_t}\bigr)\bigr|
   = \bigl|\EE\bigl(\EE\bigl(\ee^{\ii \lambda \ee^{\gamma t} X_t}
                             \bmid (Y_u)_{u\in[0,t]}, X_0\bigr)\bigr)\bigr| \\
  &= \biggl|\EE\biggl(\exp\biggl\{\ii \lambda
                                  \biggl(r \ee^{\gamma t} Y_t + Z_0
                                         + \int_0^t \ee^{\gamma u} (A - B Y_u) \, \dd u\biggr)
                                  - \frac{\lambda^2}{2}
                                    \biggl(\Sigma_2^2 \int_0^t \ee^{2\gamma u} Y_u \, \dd u
                                           + \sigma_3^2
                                             \int_0^t
                                              \ee^{2\gamma u} \, \dd u\biggr)\biggr\}\biggr)
     \biggr| \\
  &\leq \EE\biggl(\biggl|\exp\biggl\{\ii \lambda
                                     \biggl(r \ee^{\gamma t} Y_t + Z_0
                                            + \int_0^t
                                               \ee^{\gamma u} (A - B Y_u) \, \dd u\biggr)
                                     - \frac{\lambda^2}{2}
                                       \biggl(\Sigma_2^2 \int_0^t \ee^{2\gamma u} Y_u \, \dd u
                                              + \sigma_3^2
                                                \int_0^t
                                                 \ee^{2\gamma u} \, \dd u\bigr)\biggr\}\biggr|
           \biggr) \\
  &= \EE\biggl(\exp\biggl\{- \frac{\lambda^2}{2}
                             \biggl(\Sigma_2^2 \int_0^t \ee^{2\gamma u} Y_u \, \dd u
                                    + \sigma_3^2
                                      \int_0^t
                                       \ee^{2\gamma u} \, \dd u\biggr)\biggr\}\biggr) .
 \end{align*}
Convergence \eqref{lim_X} implies \ $\ee^{\gamma t} X_t \distr V_X$ \ as \ $t \to \infty$,
 \ hence, by the continuity theorem and by the monotone convergence theorem, 
 \begin{align*}
  \bigl|\EE\bigl(\ee^{\ii \lambda V_X}\bigr)\bigr|
  = \lim_{t\to\infty} \bigl|\EE\bigl(\ee^{\ii \lambda \ee^{\gamma t} X_t}\bigr)\bigr|
  &\leq \lim_{t\to\infty}
         \EE\biggl(\exp\biggl\{- \frac{\lambda^2}{2}
                                 \biggl(\Sigma_2^2 \int_0^t \ee^{2\gamma u} Y_u \, \dd u
                                        + \sigma_3^2
                                          \int_0^t
                                           \ee^{2\gamma u} \, \dd u\biggr)\biggr\}\biggr) \\
  &= \EE\biggl(\exp\biggl\{- \frac{\lambda^2}{2}
                             \biggl(\Sigma_2^2 \int_0^\infty \ee^{2\gamma u} Y_u \, \dd u
                                    + \sigma_3^2
                                      \int_0^\infty
                                       \ee^{2\gamma u} \, \dd u\biggr)\biggr\}\biggr) .
 \end{align*}
 for all \ $\lambda \in \RR$. 
\ If \ $\sigma_3 \in \RR_{++}$, \ then we have
 \[
   \bigl|\EE\bigl(\ee^{\ii \lambda V_X}\bigr)\bigr|
   \leq \exp\left\{- \frac{\sigma_3^2}{4(-\gamma)} \lambda^2\right\}
 \]
 for all \ $\lambda \in \RR$, \ hence
 \ $\int_{-\infty}^\infty \bigl|\EE\bigl(\ee^{\ii \lambda V_X}\bigr)\bigr| \, \dd\lambda
    < \infty$,
 \ implying absolute continuity of the distribution of \ $V_X$. 
 
If \ $\bigl(a - \frac{\sigma_1^2}{2}\bigr) (1 - \varrho^2) \sigma_2^2 \in \RR_{++}$,
 \ then we have
 \[
   \bigl|\EE\bigl(\ee^{\ii \lambda V_X}\bigr)\bigr|
   \leq \EE\biggl(\exp\biggl\{- \frac{\Sigma_2^2}{2} \lambda^2
                                \int_0^\infty \ee^{2\gamma u} Y_u \, \dd u\biggr\}\biggr)
   \leq \EE\biggl(\exp\biggl\{- \frac{\Sigma_2^2\ee^{4\gamma}}{2} \lambda^2
                                \int_1^2 Y_u \, \dd u\biggr\}\biggr)
 \]
 for all \ $\lambda \in \RR$. 
\ Applying the comparison theorem (see, e.g., Karatzas and Shreve \cite[5.2.18]{KarShr}), we
 obtain \ $\PP(\text{$\cY_t \leq Y_t$ \ for all \ $t \in \RR_+$}) = 1$, \ where
 \ $(\cY_t)_{t\in\RR_+}$ \ is the unique strong solution of the SDE
 \[
   \dd \cY_t = (a - b \cY_t) \, \dd t + \sigma_1 \sqrt{\cY_t} \, \dd W_t , \qquad
   t \in [0, \infty) ,
 \]
 with initial value \ $\cY_0 = 0$.
\ Consequently, taking into account \ $\Sigma_2 = \sigma_2 \sqrt{1 - \varrho^2} > 0$, \ we
 obtain
 \begin{align*}
  &\int_{-\infty}^\infty \bigl|\EE\bigl(\ee^{\ii \lambda V_X}\bigr)\bigr| \, \dd\lambda
   \leq \int_{-\infty}^\infty 
        \EE\biggl(\exp\biggl\{- \frac{\Sigma_2^2\ee^{4\gamma}}{2} \lambda^2
                                \int_1^2 \cY_u \, \dd u\biggr\}\biggr) \dd\lambda \\
  &= \EE\biggl(\int_{-\infty}^\infty 
                \exp\biggl\{- \frac{\Sigma_2^2\ee^{4\gamma}}{2} \lambda^2
                              \int_1^2 \cY_u \, \dd u\biggr\} \dd\lambda\biggr) 
   = \EE\left(\frac{\sqrt{2\pi}}
                   {\Sigma_2\ee^{2\gamma}\sqrt{\int_1^2 \cY_u \, \dd u}}\right)
   = \frac{\sqrt{2\pi}}{\Sigma_2\ee^{2\gamma}}
     \EE\left(\frac{1}{\sqrt{\int_1^2 \cY_u \, \dd u}}\right)
   < \infty ,
 \end{align*}
 whenever
 \begin{equation}\label{reciproc}
  \EE\left(\frac{1}{\sqrt{\int_1^2 \cY_u \, \dd u}}\right) < \infty .
 \end{equation}
By the Cauchy--Schwarz inequality, we have
 \[
   1 = \biggl(\int_1^2 \sqrt{\cY_u} \cdot \frac{1}{\sqrt{\cY_s}} \, \dd u\biggr)^2
     \leq \int_1^2 \cY_u \, \dd u \int_1^2 \frac{1}{\cY_u} \, \dd u ,
 \]
 hence
 \[
   \EE\left(\frac{1}{\sqrt{\int_1^2 \cY_u \, \dd u}}\right)
   \leq \EE\left(\sqrt{\int_1^2 \frac{1}{\cY_u} \, \dd u}\right)
   \leq \sqrt{\EE\left(\int_1^2 \frac{1}{\cY_u} \, \dd u\right)}
   = \sqrt{\int_1^2 \EE\biggl(\frac{1}{\cY_u}\biggr) \dd u} .
 \]
For each \ $u \in \RR_{++}$, \ we have \ $\cY_u \distre c(u) \xi$, \ where the distribution of
 \ $\xi$ \ has a chi-square distribution with degrees of freedom \ $\frac{4a}{\sigma_1^2}$
 \ and
 \ $c(u) := \frac{\sigma_1^2}{4} \int_0^u \ee^{-bv} \, \dd v
          = \frac{\sigma_1^2(\ee^{-bu}-1)}{4(-b)}$,
 \ see Proposition \ref{Pro_momentsnp}.
Hence
 \[
   \EE\biggl(\frac{1}{\cY_u}\biggr) = \frac{1}{c(u)} \EE\biggl(\frac{1}{\xi}\biggr) ,
 \]
 where \ $\EE\bigl(\frac{1}{\xi}\bigr) < \infty$, \ since the density of \ $\xi$ \ has the
 form
 \[
   \RR \ni x
   \mapsto \frac{1}{2^{\frac{2a}{\sigma_1^2}}\Gamma\bigl(\frac{2a}{\sigma_1^2}\bigr)}
           x^{\frac{2a}{\sigma_1^2}-1} \ee^{-\frac{x}{2}} \bbone_{\RR_{++}}(x)
 \]
 and the assumption \ $a - \frac{\sigma_1^2}{2} > 0$ \ yields
 \ $\frac{2a}{\sigma_1^2} - 1 > 0$.
\ Consequently,
 \[
   \int_1^2 \EE\biggl(\frac{1}{\cY_u}\biggr) \dd u
   = \EE\biggl(\frac{1}{\xi}\biggr) \int_1^2 \frac{1}{c(u)} \, \dd u
   = \EE\biggl(\frac{1}{\xi}\biggr)
     \int_1^2 \frac{4(-b)}{\sigma_1^2(\ee^{-bu}-1)} \, \dd u
   < \infty ,
 \]
 thus we obtain \eqref{reciproc}, and hence
 \ $\int_{-\infty}^\infty \bigl|\EE\bigl(\ee^{\ii \lambda V_X}\bigr)\bigr| \dd\lambda
    < \infty$,
 \ and we conclude absolute continuity of the distribution of \ $V_X$.
\proofend

\begin{Thm}\label{Thm_LSE_super}
Let us consider the two-factor affine diffusion model \eqref{2dim_affine} with
 \ $a \in \RR_+$, \ $b \in \RR_{--}$, \ $\alpha, \beta \in \RR$, \ $\gamma \in (-\infty, b)$,
 \ $\sigma_1 \in \RR_{++}$, \ $\sigma_2, \sigma_3 \in \RR_+$ \ and \ $\varrho \in [-1, 1]$
 \ with a random initial value \ $(\eta_0, \zeta_0)$ \ independent of
 \ $(W_t, B_t, L_t)_{t\in\RR_+}$ \ satisfying \ $\PP(\eta_0 \in \RR_+) = 1$. 
\ Suppose that \ $\alpha \beta \in \RR_-$. 
\ Suppose that \ $\sigma_3 \in \RR_{++}$ \ or 
 \ $\bigl(a - \frac{\sigma_1^2}{2}\bigr) (1 - \varrho^2) \sigma_2^2 \in \RR_{++}$.
\ Then
 \begin{align}\label{LSE_supercrit}
  \begin{bmatrix}
   T \ee^{\frac{bT}{2}} (\ha_T - a) \\
   \ee^{-\frac{bT}{2}} (\hb_T - b) \\
   T \ee^{\frac{bT}{2}} (\halpha_T - \alpha) \\
   \ee^{-bT/2} (\hbeta_T - \beta) \\
   \ee^{\frac{(b-2\gamma)T}{2}} (\hgamma_T - \gamma)
  \end{bmatrix}
  \distr
  \bV^{-1} \Beta \bxi
 \end{align}
 as \ $T \to \infty$ \ with
 \[
   \bV := \begin{bmatrix}
           1 & \frac{V_Y}{b} & 0 & 0 & 0 \\[1mm]
           0 & -\frac{V_Y^2}{2b} & 0 & 0 & 0 \\[1mm]
           0 & 0 & 1 & \frac{V_Y}{b} & \frac{V_X}{\gamma} \\[1mm]
           0 & 0 & 0 & -\frac{V_Y^2}{2b} & -\frac{V_YV_X}{b+\gamma} \\[1mm]
           0 & 0 & 0 & -\frac{V_YV_X}{b+\gamma} & -\frac{V_X^2}{2\gamma}
          \end{bmatrix} ,
 \]
 where \ $V_Y$ \ and \ $V_X$ \ are given in \eqref{lim_Y} and \eqref{lim_X}, respectively,
 \ $\Beta$ \ is a $5 \times 5$ random matrix such that
 \[
   \Beta \Beta^\top
   = \begin{bmatrix}
      - \frac{\sigma_1^2 V_Y}{b} & \frac{\sigma_1^2 V_Y^2}{2b}
       & - \frac{\varrho\sigma_1\sigma_2 V_Y}{b}
       & \frac{\varrho\sigma_1\sigma_2 V_Y^2}{2b}
       & \frac{\varrho\sigma_1\sigma_2 V_Y V_X}{b+\gamma} \\[1mm]
      \frac{\sigma_1^2 V_Y^2}{2b} & - \frac{\sigma_1^2 V_Y^3}{3b}
       & \frac{\varrho\sigma_1\sigma_2 V_Y^2}{2b}
       & - \frac{\varrho\sigma_1\sigma_2 V_Y^3}{3b}
       & - \frac{\varrho\sigma_1\sigma_2 V_Y^2 V_X}{2b+\gamma} \\[1mm]
      - \frac{\varrho\sigma_1\sigma_2 V_Y}{b} & \frac{\varrho\sigma_1\sigma_2 V_Y^2}{2b}
       & - \frac{\sigma_2^2 V_Y}{b} & \frac{\sigma_2^2 V_Y^2}{2b}
       & \frac{\sigma_2^2 V_Y V_X}{b+\gamma} \\[1mm]
      \frac{\varrho\sigma_1\sigma_2 V_Y^2}{2b} & - \frac{\varrho\sigma_1\sigma_2 V_Y^3}{3b}
       & \frac{\sigma_2^2 V_Y^2}{2b} & - \frac{\sigma_2^2 V_Y^3}{3b}
       & - \frac{\sigma_2^2 V_Y^2 V_X}{2b+\gamma} \\[1mm]
      \frac{\varrho\sigma_1\sigma_2 V_Y V_X}{b+\gamma}
       & - \frac{\varrho\sigma_1\sigma_2 V_Y^2 V_X}{2b+\gamma}
       & \frac{\sigma_2^2 V_Y V_X}{b+\gamma}
       & - \frac{\sigma_2^2 V_Y^2 V_X}{2b+\gamma} & - \frac{\sigma_2^2 V_Y V_X^2}{b+2\gamma}
     \end{bmatrix} ,
 \]
 and \ $\bxi$ \ is a 5-dimensional standard normally distributed random vector independent of
 \ $(V_Y, V_X)$.
 \end{Thm}

\noindent{\bf Proof.}
We have
 \[
   \begin{bmatrix}
    T \ee^{\frac{bT}{2}} (\ha_T - a) \\
    \ee^{-\frac{bT}{2}} (\hb_T - b) \\
    T \ee^{\frac{bT}{2}} (\halpha_T - \alpha) \\
    \ee^{-\frac{bT}{2}} (\hbeta_T - \beta) \\
    \ee^{\frac{(b-2\gamma)T}{2}} (\hgamma_T - \gamma)
   \end{bmatrix}
   = \diag\Bigl(T \ee^{\frac{bT}{2}}, \ee^{-\frac{bT}{2}}, T \ee^{\frac{bT}{2}},
                \ee^{-\frac{bT}{2}}, \ee^{\frac{(b-2\gamma)T}{2}}\Bigr)
     \bigl(\hbtheta_T - \btheta\bigr) ,
 \]
 where, by \eqref{LSE_cont_true},
 \[
   \hbtheta_T - \btheta
   = \bG_T^{-1} \bh_T
   = \begin{bmatrix}
      \bG_T^{(1)} & \bzero \\
      \bzero & \bG_T^{(2)}
     \end{bmatrix}^{-1}
     \begin{bmatrix}
      \bh_T^{(1)} \\
      \bh_T^{(2)}
     \end{bmatrix} .
 \]
We are going to apply Theorem \ref{THM_Zanten} for the continuous local martingale
 \ $(\bh_T)_{T\in\RR_+}$ \ with quadratic variation process
 \ $\langle \bh \rangle_T = \tbG_T$, \ $T \in \RR_+$ \ (introduced in the proof of Theorem
 \ref{Thm_LSE_sub}). 
With scaling matrices
 \[
   \bQ(T) := \diag\Bigl(\ee^{\frac{bT}{2}}, \ee^{\frac{3bT}{2}}, \ee^{\frac{bT}{2}},
                        \ee^{\frac{3bT}{2}}, \ee^{\frac{(b+2\gamma)T}{2}}\Bigr) ,
   \qquad T \in \RR_{++} ,
 \]
 by \eqref{lim_intYX}, we have
 \[
   \bQ(T) \langle \bh \rangle_T \bQ(T)^\top
   \as \begin{bmatrix}
           - \frac{\sigma_1^2 V_Y}{b} & \frac{\sigma_1^2 V_Y^2}{2b}
            & - \frac{\varrho\sigma_1\sigma_2 V_Y}{b}
            & \frac{\varrho\sigma_1\sigma_2 V_Y^2}{2b}
            & \frac{\varrho\sigma_1\sigma_2 V_Y V_X}{b+\gamma} \\
           \frac{\sigma_1^2 V_Y^2}{2b} & - \frac{\sigma_1^2 V_Y^3}{3b}
            & \frac{\varrho\sigma_1\sigma_2 V_Y^2}{2b}
            & - \frac{\varrho\sigma_1\sigma_2 V_Y^3}{3b}
            & - \frac{\varrho\sigma_1\sigma_2 V_Y^2 V_X}{2b+\gamma} \\
           - \frac{\varrho\sigma_1\sigma_2 V_Y}{b} & \frac{\varrho\sigma_1\sigma_2 V_Y^2}{2b}
            & - \frac{\sigma_2^2 V_Y}{b} & \frac{\sigma_2^2 V_Y^2}{2b}
            & \frac{\sigma_2^2 V_Y V_X}{b+\gamma} \\
           \frac{\varrho\sigma_1\sigma_2 V_Y^2}{2b}
            & - \frac{\varrho\sigma_1\sigma_2 V_Y^3}{3b}
            & \frac{\sigma_2^2 V_Y^2}{2b} & - \frac{\sigma_2^2 V_Y^3}{3b}
            & - \frac{\sigma_2^2 V_Y^2 V_X}{2b+\gamma} \\
           \frac{\varrho\sigma_1\sigma_2 V_Y V_X}{b+\gamma}
            & - \frac{\varrho\sigma_1\sigma_2 V_Y^2 V_X}{2b+\gamma}
            & \frac{\sigma_2^2 V_Y V_X}{b+\gamma}
            & - \frac{\sigma_2^2 V_Y^2 V_X}{2b+\gamma}
            & - \frac{\sigma_2^2 V_Y V_X^2}{b+2\gamma}
          \end{bmatrix}
   = \Beta \Beta^\top
 \]
 as \ $T \to \infty$. 
\ Hence by Theorem \ref{THM_Zanten}, for each random matrix \ $\bA$ \ defined on
 \ $(\Omega, \cF, \PP)$, \ we obtain
 \begin{equation}\label{Zanten}
  (\bQ(T)\bh_T, \bA) \distr (\Beta \bxi, \bA) \qquad \text{as \ $T \to \infty$,}
 \end{equation}
 where \ $\bxi$ \ is a 5-dimensional standard normally distributed random vector independent
 of \ $(\Beta, \bA)$.
\ The aim of the following discussion is to include appropriate scaling matrices for
 \ $\bG_T$.
\ The matrices \ $\bG_T^{(1)}$ \ and \ $\bG_T^{(2)}$ \ can be written in the form
 \[
   \bG_T^{(1)}
   = \diag\bigl(T^{\frac{1}{2}}, \ee^{-bT}\bigr)
     \begin{bmatrix}
      1 & -\frac{\ee^{bT}}{\sqrt{T}} \int_0^T Y_s \, \dd s \\
      -\frac{\ee^{bT}}{\sqrt{T}} \int_0^T Y_s \, \dd s & \ee^{2bT} \int_0^T Y_s^2 \, \dd s
     \end{bmatrix}
     \diag\bigl(T^{\frac{1}{2}}, \ee^{-bT}\bigr)
 \]
 and
 \begin{align*}
  \bG_T^{(2)}
  &= \diag\bigl(T^{\frac{1}{2}}, \ee^{-bT}, \ee^{-\gamma T}\bigr)
     \begin{bmatrix}
      1 & -\frac{\ee^{bT}}{\sqrt{T}} \int_0^T Y_s \, \dd s
       & -\frac{\ee^{\gamma T}}{\sqrt{T}} \int_0^T X_s \, \dd s \\
      -\frac{\ee^{bT}}{\sqrt{T}} \int_0^T Y_s \, \dd s & \ee^{2bT} \int_0^T Y_s^2 \, \dd s
       & \ee^{(b+\gamma)T} \int_0^T Y_s X_s \, \dd s \\
      -\frac{\ee^{\gamma T}}{\sqrt{T}} \int_0^T X_s \, \dd s
       & \ee^{(b+\gamma)T} \int_0^T Y_s X_s \, \dd s
       & \ee^{2\gamma T} \int_0^T X_s^2 \, \dd s
     \end{bmatrix} \\
  &\quad\times
     \diag\bigl(T^{\frac{1}{2}}, \ee^{-bT}, \ee^{-\gamma T}\bigr) ,
 \end{align*}
 hence the matrices \ $(\bG_T^{(1)})^{-1}$ \ and \ $(\bG_T^{(2)})^{-1}$ \ can be written in
 the form
 \[
   (\bG_T^{(1)})^{-1}
   = \diag\bigl(T^{-\frac{1}{2}}, \ee^{bT}\bigr)
     \begin{bmatrix}
      1 & -\frac{\ee^{bT}}{\sqrt{T}} \int_0^T Y_s \, \dd s \\
      -\frac{\ee^{bT}}{\sqrt{T}} \int_0^T Y_s \, \dd s & \ee^{2bT} \int_0^T Y_s^2 \, \dd s
     \end{bmatrix}^{-1}
     \diag\bigl(T^{-\frac{1}{2}}, \ee^{bT}\bigr)
 \]
 and
 \begin{align*}
  (\bG_T^{(2)})^{-1}
  &= \diag\bigl(T^{-\frac{1}{2}}, \ee^{bT}, \ee^{\gamma T}\bigr)
     \begin{bmatrix}
      1 & -\frac{\ee^{bT}}{\sqrt{T}} \int_0^T Y_s \, \dd s
       & -\frac{\ee^{\gamma T}}{\sqrt{T}} \int_0^T X_s \, \dd s \\
      -\frac{\ee^{bT}}{\sqrt{T}} \int_0^T Y_s \, \dd s & \ee^{2bT} \int_0^T Y_s^2 \, \dd s
       & \ee^{(b+\gamma)T} \int_0^T Y_s X_s \, \dd s \\
      -\frac{\ee^{\gamma T}}{\sqrt{T}} \int_0^T X_s \, \dd s
       & \ee^{(b+\gamma)T} \int_0^T Y_s X_s \, \dd s
       & \ee^{2\gamma T} \int_0^T X_s^2 \, \dd s
     \end{bmatrix}^{-1} \\
  &\quad\times
     \diag\bigl(T^{-\frac{1}{2}}, \ee^{bT}, \ee^{\gamma T}\bigr) .
 \end{align*}
We have
 \begin{align*}
  &\diag\Bigl(T \ee^{\frac{bT}{2}}, \ee^{-\frac{bT}{2}}, T \ee^{\frac{bT}{2}}, 
              \ee^{-\frac{bT}{2}}, \ee^{\frac{(b-2\gamma)T}{2}}\Bigr)
   \diag\Bigl(T^{-\frac{1}{2}}, \ee^{bT}, T^{-\frac{1}{2}}, \ee^{bT}, \ee^{\gamma T}\Bigr) \\
  &= \diag\Bigl(T^{\frac{1}{2}} \ee^{\frac{bT}{2}}, \ee^{\frac{bT}{2}}, T^{\frac{1}{2}}
                \ee^{\frac{bT}{2}}, \ee^{\frac{bT}{2}}, \ee^{\frac{bT}{2}}\Bigr)
 \end{align*}
 and
 \begin{align*}
  &\diag\Bigl(T^{-\frac{1}{2}}, \ee^{bT}, T^{-\frac{1}{2}}, \ee^{bT}, \ee^{\gamma T}\Bigr)
   \bQ(T)^{-1} \\
  &= \diag\Bigl(T^{-\frac{1}{2}}, \ee^{bT}, T^{-\frac{1}{2}}, \ee^{bT}, \ee^{\gamma T}\Bigr)
     \diag\Bigl(\ee^{-\frac{bT}{2}}, \ee^{-\frac{3bT}{2}}, \ee^{-\frac{bT}{2}},
                \ee^{-\frac{3bT}{2}}, \ee^{-\frac{(b+2\gamma)T}{2}}\Bigr) \\
  &= \diag\Bigl(T^{-\frac{1}{2}} \ee^{-\frac{bT}{2}}, \ee^{-\frac{bT}{2}},
                T^{-\frac{1}{2}} \ee^{-\frac{bT}{2}}, \ee^{-\frac{bT}{2}},
                \ee^{-\frac{bT}{2}}\Bigr) .
 \end{align*}
Moreover,
 \begin{align*}
  &\diag\Bigl(T^{\frac{1}{2}} \ee^{\frac{bT}{2}}, \ee^{\frac{bT}{2}}\Bigr)
   \begin{bmatrix}
    1 & -\frac{\ee^{bT}}{\sqrt{T}} \int_0^T Y_s \, \dd s \\
    -\frac{\ee^{bT}}{\sqrt{T}} \int_0^T Y_s \, \dd s & \ee^{2bT} \int_0^T Y_s^2 \, \dd s
   \end{bmatrix}
   \diag\Bigl(T^{-\frac{1}{2}} \ee^{-\frac{bT}{2}}, \ee^{-\frac{bT}{2}}\Bigr) \\
  &= \begin{bmatrix}
      1 & -\ee^{bT} \int_0^T Y_s \, \dd s \\
      -\frac{\ee^{bT}}{T} \int_0^T Y_s \, \dd s & \ee^{2bT} \int_0^T Y_s^2 \, \dd s
     \end{bmatrix}
   =: \bJ_T^{(1)}
 \end{align*}
 and
 \begin{align*}
  &\diag\Bigl(T^{\frac{1}{2}} \ee^{\frac{bT}{2}}, \ee^{\frac{bT}{2}}, 
              \ee^{\frac{bT}{2}}\Bigr)
   \begin{bmatrix}
    1 & -\frac{\ee^{bT}}{\sqrt{T}} \int_0^T Y_s \, \dd s
     & -\frac{\ee^{\gamma T}}{\sqrt{T}} \int_0^T X_s \, \dd s \\
    -\frac{\ee^{bT}}{\sqrt{T}} \int_0^T Y_s \, \dd s & \ee^{2bT} \int_0^T Y_s^2 \, \dd s
     & \ee^{(b+\gamma)T} \int_0^T Y_s X_s \, \dd s \\
    -\frac{\ee^{\gamma T}}{\sqrt{T}} \int_0^T X_s \, \dd s
     & \ee^{(b+\gamma)T} \int_0^T Y_s X_s \, \dd s & \ee^{2\gamma T} \int_0^T X_s^2 \, \dd s
   \end{bmatrix} \\
  &\times
   \diag\Bigl(T^{-\frac{1}{2}} \ee^{-\frac{bT}{2}}, \ee^{-\frac{bT}{2}},
              \ee^{-\frac{bT}{2}}\Bigr) \\
  &= \begin{bmatrix}
      1 & -\ee^{bT} \int_0^T Y_s \, \dd s & -\ee^{\gamma T} \int_0^T X_s \, \dd s \\
      -\frac{\ee^{bT}}{T} \int_0^T Y_s \, \dd s & \ee^{2bT} \int_0^T Y_s^2 \, \dd s
       & \ee^{(b+\gamma)T} \int_0^T Y_s X_s \, \dd s \\
      -\frac{\ee^{\gamma T}}{T} \int_0^T X_s \, \dd s
       & \ee^{(b+\gamma)T} \int_0^T Y_s X_s \, \dd s
       & \ee^{2\gamma T} \int_0^T X_s^2 \, \dd s
     \end{bmatrix}
   =: \bJ_T^{(2)}
 \end{align*}
Consequently,
 \[
   \begin{bmatrix}
    T \ee^{\frac{bT}{2}} (\ha_T - a) \\
    \ee^{-\frac{bT}{2}} (\hb_T - b) \\
    T \ee^{\frac{bT}{2}} (\halpha_T - \alpha) \\
    \ee^{-\frac{bT}{2}} (\hbeta_T - \beta) \\
    \ee^{\frac{(b-2\gamma)T}{2}} (\hgamma_T - \gamma)
   \end{bmatrix}
   = \diag\bigl(\bJ_T^{(1)}, \bJ_T^{(2)}\bigr)^{-1} \bQ(T) \bh_T ,
 \]
 where, by Lemma \ref{martconvYX},
 \begin{equation}\label{convV}
  \diag\bigl(\bJ_T^{(1)}, \bJ_T^{(2)}\bigr) \stoch \bV \qquad \text{as \ $T \to \infty$.}
 \end{equation}
By \eqref{Zanten} with \ $\bA = \bV$, \ by \eqref{convV} and by Theorem 2.7 (iv) of van der
 Vaart \cite{Vaart}, we obtain
 \[
   \bigl(\bQ(T)\bh_T, \diag\bigl(\bJ_T^{(1)}, \bJ_T^{(2)}\bigr)\bigr)
   \distr (\Beta \bxi, \bV) \qquad \text{as \ $T \to \infty$.}
 \]
The random matrix \ $\bV$ \ is invertible almost surely, since
 \[
   \det(\bV) = - \frac{(b-\gamma)^2V_Y^4V_X^2}{8(b+\gamma)^2b^2\gamma} > 0
 \]
 almost surely by Lemma \ref{martconvYX}.
Consequently,
 \ $\diag\bigl(\bJ_T^{(1)}, \bJ_T^{(2)}\bigr)^{-1} \bQ(T) \bh_T \distr \bV^{-1} \Beta \bxi$
 \ as \ $T \to \infty$.
\proofend

\vspace*{10mm}

\appendix

\noindent{\bf\Large Appendix}

\section{Stationarity and exponential ergodicity}

The following result states the existence of a unique stationary distribution of the affine
 diffusion process given by the SDE \eqref{2dim_affine}, see Bolyog and Pap
 \cite[Theorem 3.1]{BolPap}.
Let \ $\CC_{-} := \{ z \in \CC : \Re(z) \leq 0\}$.

\begin{Thm}\label{Thm_ergodic1}
Let us consider the two-factor affine diffusion model \eqref{2dim_affine} with
 \ $a \in \RR_+$, \ $b \in \RR_{++}$, \ $\alpha, \beta \in \RR$, \ $\gamma \in \RR_{++}$,
 \ $\sigma_1, \sigma_2, \sigma_3 \in \RR_+$, \ $\varrho \in [-1, 1]$, \ and with a random 
 initial value \ $(\eta_0, \zeta_0)$ \ independent of \ $(W_t, B_t, L_t)_{t\in\RR_+}$
 \ satisfying \ $\PP(\eta_0 \in \RR_+) = 1$.
\ Then
 \renewcommand{\labelenumi}{{\rm(\roman{enumi})}}
 \begin{enumerate}
 \item 
  $(Y_t, X_t) \distr (Y_\infty, X_\infty)$ \ as \ $t \to \infty$, \ and we have
   \begin{align}\label{help18}
    \EE\bigl(\ee^{u_1 Y_\infty + \ii \lambda_2 X_\infty}\bigr)
    = \exp\left\{a \int_0^\infty \kappa_s(u_1, \lambda_2) \, \dd s
      + \ii \frac{\alpha}{\gamma} \lambda_2 - \frac{\sigma_3^2}{4\gamma} \lambda_2^2 \right\}
   \end{align}
   for \ $(u_1, \lambda_2) \in \CC_{-} \times \RR$, \ where \ $\kappa_t(u_1, \lambda_2)$,
   \ $t \in \RR_+$, \ is the unique solution of the (deterministic) differential equation
   \begin{align}\label{DE1}
    \begin{cases}
     \frac{\partial \kappa_t}{\partial t}(u_1, \lambda_2)
     = - b \kappa_t(u_1, \lambda_2)
       - \ii \beta \ee^{-\gamma t} \lambda_2
       + \frac{1}{2} \sigma_1^2 \kappa_t(u_1, \lambda_2)^2 \\
     \phantom{\frac{\partial \kappa_t}{\partial t}(u_1, \lambda_2) =}  
       + \ii \varrho \sigma_1 \sigma_2 \ee^{-\gamma t} \lambda_2 \kappa_t(u_1, \lambda_2)
       - \frac{1}{2} \sigma_2^2 \ee^{-2\gamma t} \lambda_2^2 , \\
     \kappa_0(u_1, \lambda_2) = u_1 ;
    \end{cases}
   \end{align}
 \item
  supposing that the random initial value \ $(\eta_0, \zeta_0)$ \ has the same distribution 
   as \ $(Y_\infty, X_\infty)$ \ given in part \textup{(i)}, \ $(Y_t ,X_t)_{t\in\RR_+}$ \ is
   strictly stationary.
\end{enumerate}
\end{Thm}

In the subcritical case, the following result states the exponential ergodicity and a strong 
 law of large numbers for the process \ $(Y_t, X_t)_{t\in\RR_+}$, \ see Bolyog and Pap
 \cite[Theorem 4.1]{BolPap}.

\begin{Thm}\label{Thm_ergodic2}
Let us consider the two-factor affine diffusion model \eqref{2dim_affine} with
 \ $a, b \in \RR_{++}$, \ $\alpha, \beta \in \RR$, \ $\gamma \in \RR_{++}$,
 \ $\sigma_1 \in \RR_{++}$, \ $\sigma_2, \sigma_3 \in \RR_+$ \ and
 \ $\varrho \in [-1, 1]$ \ with a random initial value \ $(\eta_0,\zeta_0)$
 \ independent of \ $(W_t, B_t, L_t)_{t\in\RR_+}$ \ satisfying
 \ $\PP(\eta_0 \in \RR_+) = 1$. 
\ Suppose that \ $(1 - \varrho^2) \sigma_2^2 + \sigma_3^2 > 0$.
\ Then the process \ $(Y_t, X_t)_{t\in\RR_+}$ \ is exponentially ergodic, namely, there exist
 \ $\delta \in \RR_{++}$, \ $B \in \RR_{++}$ \ and \ $\kappa \in \RR_{++}$, \ such that
 \begin{equation}\label{exp_ergodic}
  \sup_{|g|\leq V+1}
   \big|\EE\big(g(Y_t, X_t) \mid (Y_0, X_0) = (y_0, x_0)\big)
        - \EE(g(Y_\infty, X_\infty))\big|
  \leq B (V(y_0, x_0) + 1) \ee^{-\delta t}
 \end{equation}
 for all \ $t \in \RR_+$ \ and \ $(y_0, x_0) \in \RR_+ \times \RR$, \ where the supremum is
 running for Borel measurable functions \ $g : \RR_+ \times \RR \to \RR$,
 \begin{align}\label{help54}
  V(y, x) := y^2 + \kappa x^2 , \qquad (y, x) \in \RR_+ \times \RR ,
 \end{align}
 and the distribution of \ $(Y_\infty, X_\infty)$ \ is given by \eqref{help18} and
 \eqref{DE1}.
Moreover, for all Borel measurable functions \ $f : \RR^2 \to \RR$ \ with
 \ $\EE(|f(Y_\infty, X_\infty)|) < \infty$, \ we have
 \begin{equation}\label{help_ergodic}
  \PP\biggl( \lim_{T\to\infty} \frac{1}{T} \int_0^T f(Y_s, X_s) \, \dd s
             = \EE(f(Y_\infty, X_\infty)) \biggr)
  = 1 .
 \end{equation}
\end{Thm}

\section{Moments}

The next proposition gives a recursive formula for the moments of the process
 \ $(Y_t, X_t)_{t\in\RR_+}$.

\begin{Pro}\label{Pro_momentsnp}
Let us consider the two-factor affine diffusion model \eqref{2dim_affine} with
 \ $a \in \RR_+$,
 \ $b, \alpha, \beta, \gamma \in \RR$, \ $\sigma_1, \sigma_2, \sigma_3 \in \RR_+$,
 \ $\varrho \in [-1, 1]$.
\ Suppose that \ $\EE(Y_0^n |X_0|^p) < \infty$ \ for some \ $n, p \in \ZZ_+$.
\ Then for each \ $t \in \RR_+$, \ we have \ $\EE(Y_t^k |X_t|^\ell) < \infty$ \ for all
 \ $k \in \{0, \ldots, n\}$ \ and \ $\ell \in \{0, \ldots, p\}$, \ and the recursion
 \begin{align*}
  \EE(Y_t^k X_t^\ell)
  &= \ee^{-(kb+\ell\gamma)t} \EE(Y_0^k X_0^\ell)
     + \Bigl(k a + \frac{1}{2} k (k-1) \sigma_1^2\Bigr)
       \int_0^t \ee^{-(kb+\ell\gamma)(t-u)} \EE(Y_u^{k-1} X_u^\ell) \, \dd u \\
  &\quad
     + (\alpha + k \varrho \sigma_1 \sigma_2)
       \int_0^t \ee^{-(kb+\ell\gamma)(t-u)} \EE(Y_u^k X_u^{\ell-1}) \, \dd u 
     - \ell \beta
       \int_0^t \ee^{-(kb+\ell\gamma)(t-u)} \EE(Y_u^{k+1} X_u^{\ell-1}) \, \dd u \\
  &\quad
     + \frac{1}{2} \ell (\ell - 1) \sigma_2^2
       \int_0^t \ee^{-(kb+\ell\gamma)(t-u)} \EE(Y_u^{k+1} X_u^{\ell-2}) \, \dd u \\
  &\quad
     + \frac{1}{2} \ell (\ell - 1) \sigma_3^2
       \int_0^t \ee^{-(kb+\ell\gamma)(t-u)} \EE(Y_u^k X_u^{\ell-2}) \, \dd u
 \end{align*}
 for all \ $t \in \RR_+$, \ where \ $\EE(Y_t^i X_t^j) := 0$ \ if \ $i, j \in \ZZ$ \ with
 \ $i < 0$ \ or \ $j < 0$.
\ Especially,
 \begin{align*}
  &\EE(Y_t) = \ee^{-bt} \EE(Y_0) + a \int_0^t \ee^{-b(t-u)} \, \dd u , \\
  &\EE(X_t) = \ee^{-\gamma t} \EE(X_0) + \alpha \int_0^t \ee^{-\gamma(t-u)} \, \dd u
              - \beta \int_0^t \ee^{-\gamma(t-u)} \EE(Y_u) \, \dd u , \\
  &\EE(Y_t^2)
   = \ee^{-2bt} \EE(Y_0^2) + (2a + \sigma_1^2) \int_0^t \ee^{-2b(t-u)} \EE(Y_u) \, \dd u , \\
  &\EE(Y_t X_t)
   = \ee^{-(b+\gamma)t} \EE(Y_0 X_0) + a \int_0^t \ee^{-(b+\gamma)(t-u)} \EE(X_u) \, \dd u \\
  &\phantom{\EE(Y_t X_t) =}
     + (\alpha + \varrho \sigma_1 \sigma_2) \int_0^t \ee^{-(b+\gamma)(t-u)} \EE(Y_u) \, \dd u
     - \beta \int_0^t \ee^{-(b+\gamma)(t-u)} \EE(Y_u^2) \, \dd u, 
 \end{align*}
 \begin{align*}
  &\EE(X_t^2)
   = \ee^{-2\gamma t} \EE(X_0^2) + \alpha \int_0^t \ee^{-2\gamma(t-u)} \EE(X_u) \, \dd u
     - 2 \beta \int_0^t \ee^{-2\gamma(t-u)} \EE(Y_u X_u) \, \dd u \\
  &\phantom{\EE(X_t^2) =}
     + \sigma_2^2 \int_0^t \ee^{-2\gamma(t-u)} \EE(Y_u) \, \dd u
     + \sigma_3^2 \int_0^t \ee^{-2\gamma(t-u)} \, \dd u .
 \end{align*}
If \ $\sigma_1 > 0$ \ and \ $Y_0 = y_0$, \ then the Laplace transform of \ $Y_t$,
 \ $t \in \RR_{++}$, \ takes the form
 \begin{align}\label{Laplace}
  \EE(\ee^{-\lambda Y_t})
  = \biggl(1 + \frac{\sigma_1^2}{2} \lambda
               \int_0^t \ee^{-bu} \, \dd u\biggr)^{-\frac{2a}{\sigma_1^2}}
    \exp\Biggl\{-\frac{\lambda\ee^{-bt}y_0}
                      {1+\frac{\sigma_1^2}{2}\lambda
                         \int_0^t\ee^{-bu}\,\dd u}\Biggr\} ,
  \qquad \lambda \in \RR_+ ,
 \end{align}
 i.e., \ $Y_t$ \ has a non-centered chi-square distribution up to a multiplicative constant
 \ $\frac{\sigma_1^2}{4} \int_0^t \ee^{-bu} \, \dd u$, \ with degrees of freedom
 \ $\frac{4a}{\sigma_1^2}$ \ and with non-centrality parameter
 \ $\frac{4\ee^{-bt}y_0}{\sigma_1^2\int_0^t\ee^{-bu}\,\dd u}$.
 
If \ $\sigma_1 > 0$ \ and \ $(1 - \varrho^2) \sigma_2^2 + \sigma_3^2 > 0$, \ then for each
 \ $t \in \RR_{++}$, \ the distribution of \ $(Y_t, X_t)$ \ is absolutely continuous.
\end{Pro}

\noindent{\bf Proof.}
It is sufficient to prove the recursion in the case when \ $(Y_0, X_0) = (y_0, x_0)$ \ with
 an arbitrary \ $(y_0, x_0) \in \RR_{++} \times \RR$, \ since then, for arbitrary initial
 values with \ $\EE(Y_0^n |X_0|^p) < \infty$, \ the recursion follows by the law of total
 expectation.
One can show that
 \begin{align}\label{help65}
  \int_0^t \EE(Y_u^k X_u^{2\ell}) \, \dd u < \infty \qquad
  \text{for all \ $t \in \RR_+$ \ and \ $k, \ell \in \ZZ_+$,}
 \end{align}
 see Bolyog and Pap \cite[proof of Theorem 5.1]{BolPap}.
For all \ $k, \ell \in \ZZ_+$, \ using the independence of \ $W$, $B$ \ and \ $L$, \ by
 It\^{o}'s formula, we have
 \begin{align*}
  &\dd(\ee^{(kb+\ell\gamma)t} Y_t^k X_t^\ell) \\
  &= (kb+\ell\gamma) \ee^{(kb+\ell\gamma)t} Y_t^k X_t^\ell \, \dd t
     + k \ee^{(kb+\ell\gamma)t} Y_t^{k-1} X_t^\ell
     \bigl[(a - b Y_t) \, \dd t + \sigma_1 \sqrt{Y_t} \, \dd W_t\bigr] \\
  &\quad
     + \ell \ee^{(kb+\ell\gamma)t} Y_t^k X_t^{\ell-1}
       \bigl[(\alpha - \beta Y_t - \gamma X_t) \, \dd t
             + \sigma_2 \sqrt{Y_t} \, (\varrho \, \dd W_t + \sqrt{1 - \varrho^2} \, \dd B_t)
             + \sigma_3 \, \dd L_t\bigr] \\
  &\quad
     + \frac{1}{2} k (k - 1) \ee^{(kb+\ell\gamma)t} Y_t^{k-2} X_t^\ell \sigma_1^2
       Y_t \, \dd t
     + k \ell \ee^{-(kb+\ell\gamma)t} Y_t^{k-1} X_t^{\ell-1}
       \varrho \sigma_1 \sigma_2 Y_t \, \dd t \\
  &\quad
     + \frac{1}{2} \ell (\ell - 1) \ee^{(kb+\ell\gamma)t} Y_t^k X_t^{\ell-2}
       (\sigma_2^2 Y_t + \sigma_3^2) \, \dd t \\
  &= k \ee^{(kb+\ell\gamma)t} Y_t^{k-1} X_t^\ell
     (a \, \dd t + \sigma_1 \sqrt{Y_t} \, \dd W_t) \\
  &\quad
     + \ell \ee^{(kb+\ell\gamma)t} Y_t^k X_t^{\ell-1}
       \bigl[(\alpha - \beta Y_t) \, \dd t
             + \sigma_2 \sqrt{Y_t} \, (\varrho \, \dd W_t + \sqrt{1 - \varrho^2} \, \dd B_t)
             + \sigma_3 \, \dd L_t\bigr] \\
  &\quad
     + \frac{1}{2} k (k - 1) \sigma_1^2 \ee^{(kb+\ell\gamma)t} Y_t^{k-1} X_t^\ell \, \dd t
     + k \ell \varrho \sigma_1 \sigma_2 \ee^{(kb+\ell\gamma)t} Y_t^k X_t^{\ell-1} \, \dd t \\
  &\quad
     + \frac{1}{2} \ell (\ell - 1) \sigma_2^2
       \ee^{(kb+\ell\gamma)t} Y_t^{k+1} X_t^{\ell-2} \, \dd t
     + \frac{1}{2} \ell (\ell - 1) \sigma_3^2 \ee^{(kb+\ell\gamma)t} Y_t^k
       X_t^{\ell-2} \, \dd t .
 \end{align*}
Writing this in an integrated form and taking the expectation of both sides, we obtain the
 recursive formulas for \ $\EE(Y_t^k |X_t|^\ell) < \infty$ \ for all
 \ $k \in \{0, \ldots, n\}$ \ and \ $\ell \in \{0, \ldots, p\}$.
 
Formula \eqref{Laplace} can be found in Ikeda and Watanabe \cite[Example 8.2]{IkeWat}.

If \ $\sigma_1 > 0$ \ and \ $(1 - \varrho^2) \sigma_2^2 + \sigma_3^2 > 0$, \ then for each
 \ $t \in \RR_{++}$, \ the conditional distribution of \ $(Y_t, X_t)$ \ given \ $(Y_0, X_0)$
 \ is absolutely continuous, see the proof of part (b) in the proof of Theorem 3.1 in Bolyog
 and Pap \cite{BolPap}.
This clearly implies that the (unconditional) distribution of \ $(Y_t, X_t)$ \ is absolutely
 continuous.
\proofend

The next theorem gives a recursive formula for the moments of the stationary distribution of
 the process \ $(Y_t, X_t)_{t\in\RR_+}$ \ in the subcritical case, see Bolyog and Pap
 \cite[Theorem 5.1]{BolPap}.

\begin{Thm}\label{Thm_moments}
Let us consider the two-factor affine diffusion model \eqref{2dim_affine} with
 \ $a \in \RR_+$, \ $b \in \RR_{++}$, \ $\alpha, \beta \in \RR$, \ $\gamma \in \RR_{++}$,
 \ $\sigma_1, \sigma_2, \sigma_3 \in \RR_+$, \ $\varrho \in [-1, 1]$, \ and the random vector
 \ $(Y_\infty, X_\infty)$ \ given by Theorem \ref{Thm_ergodic1}.
Then all the (mixed) moments of \ $(Y_\infty, X_\infty)$ \ of any order are finite, i.e., we
 have \ $\EE(Y_\infty^n |X_\infty|^p) < \infty$ \ for all \ $n, p \in \ZZ_+$, \ and the
 recursion
 \begin{align*}
  \EE(Y_\infty^n X_\infty^p)
  = \frac{1}{nb+p\gamma}
    \biggl[&\biggl(n a + \frac{1}{2} \, n (n - 1) \sigma_1^2\biggr)
            \EE(Y_\infty^{n-1} X_\infty^p)
            - p \beta \EE(Y_\infty^{n+1} X_\infty^{p-1}) \\
           &+ p (\alpha + n \varrho \sigma_1 \sigma_2) \EE(Y_\infty^n X_\infty^{p-1})
            + \frac{1}{2} \, p (p - 1) \sigma_2^2 \EE(Y_\infty^{n+1} X_\infty^{p-2}) \\
           &+ \frac{1}{2} \, p (p - 1) \sigma_3^2 \EE(Y_\infty^n X_\infty^{p-2})\biggr] ,
 \end{align*}
 holds for all \ $n, p \in \ZZ_+$ \ with \ $n + p \geq 1$,
 \ where \ $\EE(Y_\infty^k X_\infty^\ell) := 0$ \ for \ $k, \ell \in \ZZ$ \ with \ $k < 0$
 \ or \ $\ell < 0$.
\ Especially,
 \begin{gather*}
  \EE(Y_\infty) = \frac{a}{b} , \qquad
  \EE(Y_\infty^2) = \frac{a(2a+\sigma_1^2)}{2b^2} , \qquad
  \EE(Y_\infty^3) = \frac{a(a+\sigma_1^2)(2a+\sigma_1^2)}{2b^3} , \\
  \EE(X_\infty) = \frac{b\alpha-a\beta}{b\gamma} , \qquad
  \EE(Y_\infty X_\infty)
  = \frac{a\EE(X_\infty)-\beta\EE(Y_\infty^2)+(\alpha+\varrho\sigma_1\sigma_2)\EE(Y_\infty)}
         {b+\gamma} , \\
  \EE(X_\infty^2)
  = \frac{-2\beta\EE(Y_\infty X_\infty)+2\alpha\EE(X_\infty)+\sigma_2^2\EE(Y_\infty)
          +\sigma_3^2}
         {2\gamma} , \\
  \EE(Y_\infty^2 X_\infty)
  = \frac{(2a+\sigma_1^2)\EE(Y_\infty X_\infty)-\beta\EE(Y_\infty^3)
          +(\alpha+2\varrho\sigma_1\sigma_2)\EE(Y_\infty^2)}
         {2b+\gamma} , \\
   \EE(Y_\infty X_\infty^2)
   = \frac{a\EE(X_\infty^2)-2\beta\EE(Y_\infty^2 X_\infty)
           +2(\alpha+\varrho\sigma_1\sigma_2)\EE(Y_\infty X_\infty) 
           + \sigma_2^2\EE(Y_\infty^2)+\sigma_3^2\EE(Y_\infty)}{b+2\gamma} .
 \end{gather*}
If \ $\sigma_1 > 0$, \ then the Laplace transform of \ $Y_\infty$ \ takes the form
 \begin{align}\label{Laplace_stac}
  \EE(\ee^{-\lambda Y_\infty})
  = \left(1 + \frac{\sigma_1^2}{2b} \lambda \right)^{-2a/\sigma_1^2} ,
  \qquad \lambda \in \RR_+ ,
 \end{align}
 i.e., \ $Y_\infty$ \ has gamma distribution with parameters \ $2a / \sigma_1^2$ \ and
 \ $2b / \sigma_1^2$, \ hence
 \[
   \EE(Y_\infty^\kappa)
   = \frac{\Gamma\Bigl(\frac{2a}{\sigma_1^2} + \kappa\Bigr)}
          {\Bigl(\frac{2b}{\sigma_1^2}\Bigr)^\kappa
           \Gamma\Bigl(\frac{2a}{\sigma_1^2}\Bigr)} , \qquad
   \kappa \in \biggl(-\frac{2a}{\sigma_1^2}, \infty\biggr) .
 \]
If \ $\sigma_1 > 0$ \ and \ $(1 - \varrho^2) \sigma_2^2 + \sigma_3^2 > 0$, \ then the 
 distribution of \ $(Y_\infty, X_\infty)$ \ is absolutely continuous.
\end{Thm}

\section{Statistics for diffusion coefficients}
\label{sigma_varrho}

Next, for any \ $T > 0$, \ we give a statistic for \ $\sigma_1^2$, $\sigma_2^2$, $\sigma_3^2$
 \ and \ $\varrho$ \ using continuous time observations \ $(Y_t, X_t)_{t\in[0,T]}$.
\ Due to this result, we do not consider the estimation of the parameters \ $\sigma_1$,
 $\sigma_2$, $\sigma_3$ \ and \ $\varrho$, \ they are supposed to be known.

Let us consider the two-factor affine diffusion model \eqref{2dim_affine} with
 \ $a \in \RR_+$, \ $b, \alpha, \beta, \gamma \in \RR$, \ $\sigma_1 \in \RR_{++}$,
 \ $\sigma_2, \sigma_3 \in \RR_+$, \ $\varrho \in [-1, 1]$.
\ Suppose that we have \ $\PP(Y_0 \in \RR_{++})$ \ or \ $a \in \RR_{++}$.
\ Then for all \ $T \in\RR_{++}$, \ we have
 \begin{equation}\label{++}
  \PP\biggl(\int_0^T Y_u \, \dd u \in \RR_{++}\biggr) = 1 .
 \end{equation}
Indeed, if \ $\omega \in \Omega$ \ is such that \ $[0, t] \ni u \mapsto Y_u(\omega)$ \ is
 continuous and \ $Y_v(\omega) \in \RR_+$ \ for all \ $v \in\RR_+$, \ then we have
 \ $\int_0^t Y_s(\omega) \, \dd s = 0$ \ if and only if \ $Y_s(\omega) = 0$ \ for all
 \ $s \in [0, t]$.
\ Using the method of the proof of Theorem 3.1 in Barczy et.\ al \cite{BarDorLiPap}, we get
 \eqref{++}.
The (predictable) quadratic variation process of \ $Y$, $X$, \ and the (predictable)
 quadratic covariation process of \ $Y$ \ and \ $X$ \ are
 \[
   \langle Y \rangle_t = \sigma_1^2 \int_0^t Y_u \, \dd u , \qquad
   \langle X \rangle_t = \sigma_2^2 \int_0^t Y_u \, \dd u + \sigma_3^2 t , \qquad
   \langle Y, X \rangle_t = \varrho \sigma_1 \sigma_2 \int_0^t Y_u \, \dd u , \qquad
   t \in \RR_+ .
 \]
\ If, in addition, \ $a \in (\sigma_1^2, \infty)$, \ then for each \ $T \in\RR_{++}$, \ we
 have
 \begin{gather*}
  \begin{bmatrix}
   \sigma_2^2 \\
   \sigma_3^2
  \end{bmatrix}
  = \begin{bmatrix}
     \int_0^T Y_u \, \dd u & T \\
     \int_0^{T/2} Y_u \, \dd u & T/2
    \end{bmatrix}^{-1}
    \begin{bmatrix}
     \langle X \rangle_T \\
     \langle X \rangle_{T/2}
    \end{bmatrix}
  =: \begin{bmatrix}
      \hsigma_2^2(T) \\
      \hsigma_3^2(T)
     \end{bmatrix} , \\
  \sigma_1^2 = \frac{\langle Y \rangle_T}{\int_0^T Y_u \, \dd u} =: \hsigma_1^2(T) , \qquad
  \varrho
  = \frac{\langle Y, X \rangle_T}{\sigma_1\sigma_2\int_0^T Y_u \, \dd u}
  = \frac{\langle Y, X \rangle_T}
         {\bigl(\hsigma_1^2(T)\hsigma_2^2(T)\bigr)^{\frac{1}{2}}\int_0^T Y_u \, \dd u}
  =: \hvarrho(T) ,
 \end{gather*}
 since the matrix
 \[
   \begin{bmatrix}
    \int_0^T Y_u \, \dd u & T \\
    \int_0^{T/2} Y_u \, \dd u & T/2
   \end{bmatrix}
 \]
 is invertible almost surely.
Indeed,
 \begin{align*}
  &\PP\biggl(\frac{T}{2} \int_0^T Y_u \, \dd u - T \int_0^{T/2} Y_u \, \dd u = 0\biggr)
   = \PP\biggl(\int_{T/2}^T Y_u \, \dd u = \int_0^{T/2} Y_u \, \dd u\biggr) \\
  &= \EE\biggl(\PP\biggl(\int_{T/2}^T Y_u \, \dd u = \int_0^{T/2} Y_u \, \dd u 
                         \,\bigg|\, (Y_u)_{u\in[0,T/2]}\biggr)\biggr) \\
  &= \EE\biggl(\PP\biggl(\int_{T/2}^T Y_u \, \dd u = I
                         \,\bigg|\, (Y_u)_{u\in[0,T/2]}\biggr)
               \bigg|_{I = \int_0^{T/2} Y_u \, \dd u}\biggr) \\
  &= \EE\biggl(\PP\biggl(\int_{T/2}^T Y_u \, \dd u = I
                         \,\bigg|\, Y_{T/2} = y\biggr)
               \bigg|_{I = \int_0^{T/2} Y_u \, \dd u, \, y = Y_{T/2}}\biggr) \\
  &= \EE\biggl(\PP\biggl(\int_0^{T/2} Y_u \, \dd u = I
                         \,\bigg|\, Y_0 = y\biggr)
               \bigg|_{I = \int_0^{T/2} Y_u \, \dd u, \, y = Y_{T/2}}\biggr)
   = 0 ,
 \end{align*}
 since
 \[
   \PP\biggl(\int_0^{T/2} Y_u \, \dd u = I \,\bigg|\, Y_0 = y\biggr) = 0
 \]
 for each \ $I \in \RR_+$ \ and \ $y \in \RR_+$, \ because the additional condition
 \ $a \in (\sigma_1^2, \infty)$ \ yields that the distribution of \ $\int_0^{T/2} Y_u \, \dd u$
 \ is absolutely continuous, see Filipovi\'{c} et al.\ \cite[Theorem 4.3]{FilMaySch}.
(The absolute continuity of \ $\int_0^{T/2} Y_u \, \dd u$ \ might hold without the additional
 condition \ $a \in (\sigma_1^2, \infty)$.)

We note that \ $(\hsigma_1^2(T), \hsigma_2^2(T), \hsigma_3^2(T), \hvarrho(T))$ \ is a
 statistic, i.e., there exists a measurable function \ $\Xi : D([0,T], \RR) \to \RR^4$ \ such
 that
 \ $(\hsigma_1^2(T), \hsigma_2^2(T), \hsigma_3^2(T), \hvarrho(T)) = \Xi((Y_u)_{u\in[0,T]})$,
 \ where \ $D([0, T], \RR)$ \ denotes the space of real-valued c\`adl\`ag functions defined
 on \ $[0, T]$, \ since
 \begin{gather}\label{Sigma1}
   \begin{bmatrix}
    \frac{1}{n} \sum_{i=1}^\nT Y_{\frac{i-1}{n}} & T \\
    \frac{1}{n} \sum_{i=1}^{\nT/2} Y_{\frac{i-1}{n}} & T/2
   \end{bmatrix}^{-1}
   \begin{bmatrix}
    \sum_{i=1}^\nT \bigl(X_{\frac{i}{n}} - X_{\frac{i-1}{n}}\bigr)^2 \\
    \sum_{i=1}^{\nT/2} \bigl(X_{\frac{i}{n}} - X_{\frac{i-1}{n}}\bigr)^2
   \end{bmatrix}
   \stoch 
   \begin{bmatrix}
    \hsigma_2^2(T) \\
    \hsigma_3^2(T)
   \end{bmatrix} , \\
   \frac{\sum_{i=1}^\nT \bigl(Y_{\frac{i}{n}} - Y_{\frac{i-1}{n}}\bigr)^2}
        {\frac{1}{n} \sum_{i=1}^\nT Y_{\frac{i-1}{n}}}
   \stoch \hsigma_1^2(T) , \qquad
   \frac{\sum_{i=1}^\nT
          \bigl(Y_{\frac{i}{n}} - Y_{\frac{i-1}{n}}\bigr)
          \bigl(X_{\frac{i}{n}} - X_{\frac{i-1}{n}}\bigr)}
        {\frac{1}{n} \sum_{i=1}^\nT Y_{\frac{i-1}{n}}}
   \stoch \hvarrho(T) \label{Sigma2}
 \end{gather}
 as \ $n \to \infty$, \ where the convergences in \eqref{Sigma1} and \eqref{Sigma2} hold
 almost surely along a suitable subsequence, the members of the sequences in \eqref{Sigma1}
 and \eqref{Sigma2} are measurable functions of \ $(Y_u, X_u)_{u\in[0,T]}$, \ and one can
 use Theorems 4.2.2 and 4.2.8 in Dudley \cite{Dud}.
Next we prove \eqref{Sigma1} and \eqref{Sigma2}.
By Theorems I.4.47 a) and I.4.52 in Jacod and Shiryaev \cite{JSh},
 \begin{gather*}
  \sum_{i=1}^\nT \bigl(X_{\frac{i}{n}} - X_{\frac{i-1}{n}}\bigr)^2
  \stoch
  [X]_T = \langle X \rangle_T , \qquad
  \sum_{i=1}^{\nT/2} \bigl(X_{\frac{i}{n}} - X_{\frac{i-1}{n}}\bigr)^2
  \stoch
  [X]_{T/2} = \langle X \rangle_{T/2} , \\
  \sum_{i=1}^\nT \bigl(Y_{\frac{i}{n}} - Y_{\frac{i-1}{n}}\bigr)^2
  \stoch
  [Y]_T = \langle Y \rangle_T , \qquad
  \sum_{i=1}^\nT
   \bigl(Y_{\frac{i}{n}} - Y_{\frac{i-1}{n}}\bigr)
   \bigl(X_{\frac{i}{n}} - X_{\frac{i-1}{n}}\bigr)
  \stoch
  [Y, X]_T = \langle Y, X \rangle_T
 \end{gather*}
 as \ $n \to \infty$.
\ Moreover, for all \ $T \in \RR_+$, \ we have
 \[
   \frac{1}{n} \sum_{i=1}^\nT Y_{\frac{i-1}{n}}
   \stoch \int_0^T Y_u \, \dd u , \qquad
   \frac{1}{n} \sum_{i=1}^{\nT/2} Y_{\frac{i-1}{n}}
   \stoch \int_0^{T/2} Y_u \, \dd u
 \]
 as \ $n \to \infty$, \ see Proposition I.4.44 in Jacod and Shiryaev \cite{JSh}.
Hence \eqref{Sigma1} and \eqref{Sigma2} follow by Slutsky's lemma.

\section{Limit theorems for continuous local martingales}

In what follows we recall some limit theorems for continuous local martingales.
We use these limit theorems for studying the asymptotic behaviour of the MLE of
 \ $\btheta = (a, b, \alpha, \beta, \gamma)^\top$.
\ First we recall a strong law of large numbers for continuous local martingales.

\begin{Thm}{\bf (Liptser and Shiryaev \cite[Lemma 17.4]{LipShiII})}
\label{DDS_stoch_int}
Let \ $\bigl( \Omega, \cF, (\cF_t)_{t\in\RR_+}, \PP \bigr)$ \ be a filtered probability space
 satisfying the usual conditions.
Let \ $(M_t)_{t\in\RR_+}$ \ be a square-integrable continuous local martingale with respect
 to the filtration \ $(\cF_t)_{t\in\RR_+}$ \ such that \ $\PP(M_0 = 0) = 1$.
\ Let \ $(\xi_t)_{t\in\RR_+}$ \ be a progressively measurable process such that
 \[
   \PP\left( \int_0^t \xi_u^2 \, \dd \langle M \rangle_u < \infty \right) = 1 ,
   \qquad t \in \RR_+ ,
 \]
 and
 \begin{align}\label{SEGED_STRONG_CONSISTENCY2}
  \int_0^t \xi_u^2 \, \dd \langle M \rangle_u \as \infty \qquad
  \text{as \ $t \to \infty$,}
 \end{align}
 where \ $(\langle M \rangle_t)_{t\in\RR_+}$ \ denotes the quadratic variation process of
 \ $M$.
\ Then
 \begin{align}\label{SEGED_STOCH_INT_SLLN}
  \frac{\int_0^t \xi_u \, \dd M_u}
       {\int_0^t \xi_u^2 \, \dd \langle M \rangle_u} \as 0 \qquad
  \text{as \ $t \to \infty$.}
 \end{align}
If \ $(M_t)_{t\in\RR_+}$ \ is a standard Wiener process, the progressive measurability of
 \ $(\xi_t)_{t\in\RR_+}$ \ can be relaxed to measurability and adaptedness to the filtration
 \ $(\cF_t)_{t\in\RR_+}$.
\end{Thm}

The next theorem is about the asymptotic behaviour of continuous multivariate local
 martingales, see van Zanten \cite[Theorem 4.1]{Zan}.

\begin{Thm}{\bf (van Zanten \cite[Theorem 4.1]{Zan})}\label{THM_Zanten}
Let \ $\bigl( \Omega, \cF, (\cF_t)_{t\in\RR_+}, \PP \bigr)$ \ be a filtered probability space
 satisfying the usual conditions.
Let \ $(\bM_t)_{t\in\RR_+}$ \ be a $d$-dimensional square-integrable continuous local
 martingale with respect to the filtration \ $(\cF_t)_{t\in\RR_+}$ \ such that
 \ $\PP(\bM_0 = \bzero) = 1$.
\ Suppose that there exists a function \ $\bQ : [t_0, \infty) \to \RR^{d \times d}$ \ with
 some \ $t_0 \in \RR_+$ \ such that \ $\bQ(t)$ \ is an invertible (non-random) matrix for all
 \ $t \in \RR_+$, \ $\lim_{t\to\infty} \|\bQ(t)\| = 0$ \ and
 \[
   \bQ(t) \langle \bM \rangle_t \, \bQ(t)^\top \stoch \bfeta \bfeta^\top
   \qquad \text{as \ $t \to \infty$,}
 \]
 where \ $\bfeta$ \ is a \ $d \times d$ random matrix.
Then, for each \ $\RR^{k\times\ell}$-valued random matrix \ $\bA$ \ defined on
 \ $(\Omega, \cF, \PP)$, \ we have
 \[
   (\bQ(t) \bM_t, \bA) \distr (\bfeta \bZ, \bA) \qquad
   \text{as \ $t \to \infty$,}
 \]
 where \ $\bZ$ \ is a \ $d$-dimensional standard normally distributed random vector
 independent of \ $(\bfeta, \bA)$.
\end{Thm}

\end{document}